\documentclass[hidelinks,onefignum,onetabnum]{siamart220329}
\usepackage{amssymb}
\usepackage{amsmath}
\usepackage{enumitem}
\usepackage{bm}
\usepackage{color}
\usepackage{exscale}
\usepackage{relsize}
\usepackage{graphicx}
\DeclareSymbolFont{yhlargesymbols}{OMX}{yhex}{m}{n}
\DeclareSymbolFont{ugmL}{OMX}{mdugm}{m}{n}
\DeclareMathAccent{\wideparen}{\mathord}{yhlargesymbols}{"F3}
\usepackage{wrapfig}
\usepackage{multicol}
\usepackage{verbatim}
\newtheorem{rmk}{Remark}
\usepackage{appendix}

\usepackage{lipsum}
\usepackage{amsfonts}
\usepackage{graphicx}
\usepackage{epstopdf}
\usepackage{algorithmic}
\ifpdf
  \DeclareGraphicsExtensions{.eps,.pdf,.png,.jpg}
\else
  \DeclareGraphicsExtensions{.eps}
\fi

\newsiamremark{remark}{Remark}
\newsiamremark{hypothesis}{Hypothesis}
\crefname{hypothesis}{Hypothesis}{Hypotheses}
\newsiamthm{claim}{Claim}

\headers{Multi-scale method for thermo-electro-mechanical coupling problems}{H. Dong}

\title{High-order multi-scale method and its convergence analysis for nonlinear thermo-electro-mechanical coupling problems of composite structures\thanks{Submitted to the editors' DATE.
\funding{This work was supported by the National Natural Science Foundation of China (No.\hspace{1mm}12471387), Xidian University Specially Funded Project for Interdisciplinary Exploration (No.\hspace{1mm}TZJH2024008), the Fundamental Research Funds for the Central Universities (No.\hspace{1mm}QTZX25082), the Innovation Capability Support Program of Shaanxi Province (No.\hspace{1mm}2024RS-CXTD-88).
}}}

\author{Hao Dong\thanks{Corresponding author. School of Mathematics and Statistics, Xidian University, Xi'an 710071, China
  (\email{donghao@mail.nwpu.edu.cn}).}
}

\usepackage{amsopn}

\ifpdf
\hypersetup{
  pdftitle={High-order multi-scale method for thermo-electro-mechanical coupling problems},
  pdfauthor={H. Dong}
}
\fi

\begin{document}

\maketitle

\begin{abstract}
This study proposes a high-order multi-scale method tailored for time-dependent nonlinear thermo-electro-mechanical coupling problems of composite structures with highly spatial heterogeneity, which incorporate temperature-dependent material properties and Joule heating effect. By employing the multi-scale asymptotic approach and the Taylor series technique, a high-accuracy multi-scale asymptotic model featuring novel high-order correction terms is established for nonlinear multi-physics simulation of periodic solid structures. A local point-wise error analysis is derived to theoretically and physically illustrate the local balance preserving of heat quantity, electric charge and stress,thereby enabling high-accuracy multi-scale computation. Moreover, a global error estimation is obtained that provides an explicit convergence rate for high-order multi-scale solutions. Furthermore, an efficient numerical algorithm featuring with off-line and on-line stages is presented meticulously, accompanied by a corresponding error analysis. Numerical experiments are conducted to showcase the competitive advantages of the proposed method for simulating the time-dependent nonlinear thermo-electro-mechanical coupling problems with highly oscillatory and discontinuous coefficients, demonstrating superior numerical accuracy and reduced computational cost.
\end{abstract}



\section{Introduction}
As the integrated circuit industry enters the post-Moore's Law era, the electronic packaging structures serve as a pivotal bridge connecting chips and external systems, and its importance is gaining increasing prominence. The next-generation integrated circuits have shifted from the traditional technique of continuously reducing structural sizes to advanced three-dimensional (3D) packaging technologies, which is regarded as a crucial scheme for extending Moore's Law. Electronic packaging structures involve a variety of materials with distinct properties, while their components span several orders of magnitude in spatial size. Hence, electronic packaging structure is a complicated three-dimensional composite system with pronounced multi-scale spatial features \cite{R1}. In real-world applications, the rising power density necessitates that electronic packaging structures withstand extremely thermo-electro-mechanical coupling environments, whose localized heat flux densities can reach up to 1kW/cm$^2$. These composite structures will exhibit significant and dominant nonlinear effects, including temperature-dependent material properties and Joule heating effects \cite{R2,R3}, which cannot be neglected in modeling, simulation, and application, and lead to enormous challenges by use of classical numerical methods. Therefore, high-accuracy and efficient multi-scale computation is of paramount importance for simulating and predicting nonlinear thermo-electro-mechanical coupling behaviors of composite structures in electronic packaging industry.

Within the continuum mechanics (CM) framework, the following partial differential equations (PDEs)-based multi-scale model for time-dependent nonlinear thermo-electro-mechanical coupling problems of composite structures with periodic micro-structures can be formulated on the basis of the first law of thermodynamics, the Fourier's law and the conservation law of electric charge \cite{R3,R4} over domain $\Omega\subset\mathbb{R}^n(n=2,3)$, where $\Omega$ is a bounded convex domain with the Lipschitz continuous boundary $\partial\Omega$.
\begin{equation}
\left\{
\begin{aligned}
&\rho^{\varepsilon} ({\bm{x}},{T^{\varepsilon}})c^{\varepsilon}({\bm{x}},{T^{\varepsilon}})\frac{{\partial {T^{\varepsilon}(\bm{x},t)}}}{{\partial t}}-\frac{\partial }{{\partial {x_i}}}\Big( {{k_{ij}^{\varepsilon}}({\bm{x}},{T^{\varepsilon}})\frac{{\partial {T^{\varepsilon}(\bm{x},t)}}}{{\partial {x_j}}}}\Big)\\
&\quad\quad={{\lambda_{ij}^{\varepsilon}}({\bm{x}},{T^{\varepsilon}})\frac{\partial \mathit{\Phi}^{\varepsilon}(\bm{x},t)}{\partial {x_i}}}\frac{\partial \mathit{\Phi}^{\varepsilon}(\bm{x},t)}{\partial {x_j}}-T^{\varepsilon}(\bm{x},t)\beta_{ij}^{\varepsilon}({\bm{x}},{T^{\varepsilon}})\frac{\partial}{\partial t}\Big(\frac{\partial U_i^{\varepsilon}(\bm{x},t)}{\partial x_j}\Big)\\
&\quad\quad+f_T(\bm{x},t),\;\text{in}\;\Omega\times(0,\mathcal T),\\
&- \frac{\partial }{{\partial {x_i}}}\Big({{\lambda_{ij}^{\varepsilon}}({\bm{x}},{T^{\varepsilon} })\frac{{\partial \mathit{\Phi}^{\varepsilon}(\bm{x},t)}}{{\partial {x_j}}}} \Big) = f_\mathit{\Phi}(\bm{x},t),\;\text{in}\;\Omega\times(0,\mathcal T),\\
&\rho^{\varepsilon} ({\bm{x}},{T^{\varepsilon}})\frac{{\partial^2 {U_i^{\varepsilon}(\bm{x},t)}}}{{\partial t^2}}-\frac{\partial }{{\partial {x_j}}}\Big( {{c_{ijkl}^{\varepsilon}}({\bm{x}},{T^{\varepsilon}})\frac{{\partial {U_k^{\varepsilon}(\bm{x},t)}}}{{\partial {x_l}}}}\Big)\\
&\quad\quad=-\frac{\partial }{{\partial {x_j}}}\Big(\beta_{ij}^{\varepsilon}({\bm{x}},{T^{\varepsilon}})(T^{\varepsilon}(\bm{x},t)-\widetilde T)\Big)+f_i(\bm{x},t),\;\text{in}\;\Omega\times(0,\mathcal T),\\
&T^{{\varepsilon}}(\bm{x},t) = \widehat T(\bm{x},t),\;\mathit{\Phi}^{{\varepsilon}}(\bm{x},t)=\widehat{\mathit{\Phi}}(\bm{x},t),\;{\bm{U}}^{{\varepsilon}}(\bm{x},t)=\widehat {\bm{U}}(\bm{x},t),\;\text{on}\;\partial\Omega\!\times\!(0,\mathcal T),\\
&T^{{\varepsilon}}({\bm{x}},0)=\widetilde T,\;{\bm{U}}^{{\varepsilon}}({\bm{x}},0)=\widetilde{\bm{U}}^0({\bm{x}}),\;\frac{\partial {\bm{U}}^{{\varepsilon}}({\bm{x}},0)}{\partial t}=\widetilde{\bm{U}}^1({\bm{x}}),\;\text{in}\;\Omega.
\end{aligned} \right.
\end{equation}
In the multi-physics coupling model (1.1), the temperature field $T^{\varepsilon}(\bm{x},t)$, electric potential field $\mathit{\Phi}^{{\varepsilon}}(\bm{x},t)$ and displacement field ${\bm{U}}^{{\varepsilon}}(\bm{x},t)$ are targeted and undetermined, in which $\varepsilon$ represents the characteristic periodic length. $\widehat T(\bm{x},t)$, $\widehat{\mathit{\Phi}}(\bm{x},t)$ and $\widehat {\bm{U}}(\bm{x},t)$ are the prescribed temperature, electric potential and displacement on the domain boundary $\partial\Omega$.
$\widetilde T$, $\widetilde{\bm{U}}^0({\bm{x}})$ and $\widetilde{\bm{U}}^1({\bm{x}})$ respectively represents the initial temperature, displacement and velocity of domain $\Omega$. $f_T(\bm{x},t)$, $f_\mathit{\Phi}(\bm{x},t)$ and $f_i(\bm{x},t)$ denote the internal heat source, electric charge density and body force respectively. Here $\rho^{\varepsilon}({\bm{x}},{T^{\varepsilon}})$, $c^{\varepsilon}({\bm{x}},{T^{\varepsilon}})$, $\{{k_{ij}^{\varepsilon}}({\bm{x}},{T^{\varepsilon}})\}$,  $\{{\lambda_{ij}^{\varepsilon}}({\bm{x}},{T^{\varepsilon}})\}$, $\{{\beta_{ij}^{\varepsilon}}({\bm{x}},{T^{\varepsilon}})\}$ and $\{{c_{ijkl}^{\varepsilon}}({\bm{x}},{T^{\varepsilon}})\}$ are, respectively, the mass density, the specific heat, the thermal conductivity tensor, the electric conductivity tensor, the thermal modulus and the elastic modulus, which all are temperature-dependent. The multi-physics coupling model (1.1) governs the thermal, electronic and mechanical coupling interactions of composite structures in the field of electronic packaging, whose inherent multi-scale nature stems from the high-frequency oscillating coefficients arising from periodic micro-heterogeneities when $0<\varepsilon\ll1$. This necessitates an extremely fine discretization to achieve high-resolution simulations, thereby rendering the computational burden prohibitive, particularly for time-dependent scenarios. Furthermore, the Laplace transform technique \cite{R5,R6,R7}, typically employed for time-dependent linear systems, is inapplicable for the theoretical and numerical investigation of the nonlinear model (1.1) due to its temperature-dependent material properties (nonlinear equation coefficients) and Joule heating effects (nonlinear product terms of gradient).

Over the past few decades, extensive researches have been conducted on the theoretical analysis and numerical simulation of thermo-electro-mechanical coupling problems. In 2006 and 2009, researchers established the existence theory for solutions to nonlinear thermo-electro-mechanical coupling problems \cite{R8,R9}. In 2010, Fern$\rm{\acute{a}}$ndez and Kuttler established a fully discrete numerical scheme for the nonlinear thermo-electro-mechanical coupling problem based on finite element method and a forward Euler scheme, and also obtained a linear convergence estimate under appropriate regularity assumptions \cite{R10}. In 2017, M$\rm{\mathring{a}}$lqvist and Stillfjord presented a fully discrete numerical scheme for the nonlinear thermo-electro-mechanical coupling problem based on finite element method in space and a semi-implicit Euler scheme in time, while proving optimal convergence orders, i.e. second-order in space and first-order in time \cite{R11}. In 2023, Jiang et al. established a Galerkin mixed finite element method for the nonlinear thermo-electro-mechanical coupling problem, and employed space-time error splitting technique to obtain optimal error estimate \cite{R12}. However, the above-mentioned computational methods are only applicable for the computation of homogeneous materials and don't consider the impact of microscopic material heterogeneities on nonlinear multi-physics behaviors of composite structures.

The multi-scale nature of composite structures presents inherent challenges in balancing accuracy and efficiency. Moreover, in the context of large-scale composite structures, classical numerical methods suffer from convergence issues caused by the prohibitive mesh refinement. To address these, researchers have developed a variety of multi-scale methodologies, including Asymptotic homogenization method (AHM) \cite{R13,R14,R15}, Multi-scale finite element method (MsFEM) \cite{R16,R17}, Heterogeneous multi-scale method (HMM) \cite{R18}, Variational multi-scale method (VMS) \cite{R19}, Multi-scale eigenelement method (MEM) \cite{R20}, Multi-scale finite volume method (MSFV) \cite{R21}, Localized orthogonal decomposition method (LOD) \cite{R22} and Multi-scale spectral generalized finite element method (MS-GFEM) \cite{R23}, etc. Driven by real-world engineering applications, Cui and his research team systematically established a class of high-order multi-scale approaches for precisely and efficiently simulating the thermal, mechanical and multi-physics behaviors of composite structures, as shown in references \cite{R5,R24,R25,R26,R27} for further details. It is important to recognize, however, that the aforementioned multi-scale approaches are predominantly tailored to linear problems of composite structures in practical applications. As the range of applications for composites broadens, their service environments have become increasingly complex and extreme. To accurately characterize the intricate service behavior of these materials, nonlinear mathematical models have been developed. Therefore, research attention of scientists and engineers has naturally shifted towards multi-scale modeling and high-accuracy simulation for nonlinear problems of heterogeneous solids. In reference \cite{R17}, Efendiev et al. developed multi-scale finite element method for accurately simulating the elliptic and parabolic heat conduction problems of heterogeneous materials with temperature-dependent equation coefficients. Moreover, generalized multi-scale finite element methods are developed by Efendiev et al. for nonlinear elliptic equations in reference \cite{R28}. Furthermore, based on the asymptotic homogenization method, researchers established effective approaches to simulate static and transient nonlinear heat conduction problems of inhomogeneous solids with temperature-dependent material coefficients in the existing studies \cite{R29,R30,R31,R32}. Besides, the G-convergence and homogenized theory of nonlinear partial differential operators was systematically given in reference \cite{R33}. For the nonlinear multi-physics problems of composite structures, Dong et al. first established high-order multi-scale method and obtained corresponding convergence analysis for nonlinear thermo-mechanical and thermo-electronic coupling problems of composite structures in \cite{R34,R35,R36}. In summary, few studies have addressed nonlinear multi-physics problems of composite structures within the multi-scale modeling and computation. Nevertheless, widespread engineering demands strongly necessitate continued research into this challenging issue, especially for nonlinear thermo-electro-mechanical coupling problems in integrated circuit industry.

The remainder of this study is organized as follows. By virtue of multi-scale asymptotic expansion and the Taylor series techniques, Section 2 establishes the high-order multi-scale computational model for nonlinear thermo-electro-mechanical coupling simulation of composite structures. After introducing the high-order multi-scale computational model, both local and global error analyses are conducted in Section 3 for the proposed multi-scale approximate solutions. With the high-order multi-scale model obtained, in Section 4 a two-stage numerical algorithm with off-line micro-scale computation, and on-line macro-scale and multi-scale computation is provided to effectively compute time-dependent nonlinear thermo-electro-mechanical coupling problems of composite structures. And also, the corresponding error analysis of the two-stage algorithm is obtained in Section 5. Numerical examples and results are presented in Section 6 to validate the numerical accuracy and efficiency of the proposed high-order multi-scale approach. Eventually, concluding remarks and potential directions are proclaimed in Section 7. Throughout this study, Einstein summation convention is utilized to streamline repetitive indices.

\section{High-order multi-scale computational model of nonlinear thermo-electro-mechanical coupling problems}
\subsection{The statement of multi-scale nonlinear coupling system}
Following the theoretical framework of AHM, we introduce the macroscopic coordinates $\bm{x}$ and corresponding microscopic coordinates $\bm{y}={\bm{x}}/{\varepsilon}=({x_1}/{\varepsilon},\cdots,{x_n}/{\varepsilon})=(y_1,\cdots,y_n)$ of periodic unit cell (PUC) $\mathbf{Y}=[0,1]^n$. Thereby, the material parameters $\rho^{\varepsilon}({\bm{x}},{T^{\varepsilon}})$, $c^{\varepsilon}({\bm{x}},{T^{\varepsilon}})$, ${k_{ij}^{\varepsilon}}({\bm{x}},{T^{\varepsilon}})$,  ${\lambda_{ij}^{\varepsilon}}({\bm{x}},{T^{\varepsilon}})$, ${\beta_{ij}^{\varepsilon}}({\bm{x}},{T^{\varepsilon}})$ and ${c_{ijkl}^{\varepsilon}}({\bm{x}},{T^{\varepsilon}})$ in multi-scale nonlinear system (1.1) can be reformulated as new forms $\rho({\bm{y}},{T^{\varepsilon}})$, $c({\bm{y}},{T^{\varepsilon}})$, ${k_{ij}}({\bm{y}},{T^{\varepsilon}})$, ${\lambda_{ij}}({\bm{y}},{T^{\varepsilon}})$, ${\beta_{ij}}({\bm{y}},{T^{\varepsilon}})$ and ${c_{ijkl}}({\bm{y}},{T^{\varepsilon}})$, which necessitate these material parameters be 1-periodic functions regarding micro-scale variable ${\bm{y}}$.

Inspired by the well-posedness theory in references \cite{R8,R9,R10,R11,R12}, some assumptions for multi-scale nonlinear equations (1.1) are provided as follows.
\begin{enumerate}
\item[(A$_1$)]
${k_{ij}^{\varepsilon}}$, ${\lambda_{ij}^{\varepsilon}}$, ${\beta_{ij}^{\varepsilon}}$ and ${c_{ijkl}^{\varepsilon}}$ are symmetric, and there exist two positive constants $\kappa_0$ and $\kappa_1$ irrespective of $\varepsilon$ for the following uniform elliptic conditions
\begin{displaymath}
\begin{aligned}
&k_{ij}^{\varepsilon}=k_{ji}^{\varepsilon},\;\kappa_0|\bm{\xi}|^2\leq {k_{ij}^{\varepsilon}}({\bm{x}},{T^{\varepsilon}})\xi_i\xi_j\le\kappa_1|\bm{\xi}|^2,\\
&\lambda_{ij}^{\varepsilon}=\lambda_{ji}^{\varepsilon},\;\kappa_0|\bm{\xi}|^2\leq {\lambda_{ij}^{\varepsilon}}({\bm{x}},{T^{\varepsilon}})\xi_i\xi_j\le\kappa_1|\bm{\xi}|^2,\\
&\beta_{ij}^{\varepsilon}=\beta_{ji}^{\varepsilon},\;\kappa_0|\bm{\xi}|^2\leq {\beta_{ij}^{\varepsilon}}({\bm{x}},{T^{\varepsilon}})\xi_i\xi_j\le\kappa_1|\bm{\xi}|^2,\\
&c_{ijkl}^{\varepsilon}=c_{jikl}^{\varepsilon}=c_{klij}^{\varepsilon},\;\kappa_0|\bm{\zeta}|^2\leq {c_{ijkl}^{\varepsilon}}({\bm{x}},{T^{\varepsilon}})\zeta_{ij}\zeta_{kl}\le\kappa_1|\bm{\zeta}|^2,
\end{aligned}
\end{displaymath}
where $\bm{\xi}=(\xi_1,\cdots,\xi_n)$ is an optional vector with
real elements in $\mathbb{R}^n$, $\bm{\zeta}=\{\zeta_{ij}\}$ is an arbitrary symmetric matrix with real elements in $\mathbb{R}^{n\times n}$, and ${\bm{x}}$ is an optional point within $\Omega$.
\item[(A$_2$)]
$\rho^{\varepsilon}$, $c^{\varepsilon}$, ${k_{ij}^{\varepsilon}}$, ${\lambda_{ij}^{\varepsilon}}$, ${\beta_{ij}^{\varepsilon}}$ and ${c_{ijkl}^{\varepsilon}}\in L^\infty (\Omega)$; $0<\rho_{*}\leq\rho^{\varepsilon} ({\bm{x}},{T^{\varepsilon}})<\rho^{*}$, $0<c_{*}\leq c^{\varepsilon} ({\bm{x}},{T^{\varepsilon}})<c^{*}$, where $\rho_{*}$, $\rho^{*}$, $c_{*}$ and $c^{*}$ are constants irrespective of $\varepsilon$.
\item[(A$_3$)]
$f_{T}({\bm{x}},t)\in L^2(\Omega\times(0,T))$, $f_\mathit{\Phi}({\bm{x}},t)\in L^2(\Omega\times(0,T))$, $f_i({\bm{x}},t)\in L^2(\Omega\times(0,T))$, ${\widehat T}({\bm{x}},t)\in H^{\frac{7}{2},\frac{7}{4}}(\partial\Omega_1\times(0,T))$, $\widehat{\mathit{\Phi}}(\bm{x},t)\in H^{\frac{7}{2},\frac{7}{4}}(\partial\Omega_2\times(0,T))$, $\widehat {\bm{U}}(\bm{x},t)\in (H^{\frac{7}{2},\frac{7}{4}}(\partial\Omega_3\times(0,T)))^n$, $\widetilde{\bm{U}}^0({\bm{x}})\in{(H^1(\Omega))^n}$, $\widetilde{\bm{U}}^1({\bm{x}})\in{(L^2(\Omega))^n}$.
\end{enumerate}
\subsection{High-order multi-scale analysis for nonlinear multi-physics problem}
Before establishing high-order multi-scale computational model, we firstly posits that the physical fields ${T^{\varepsilon} }(\bm{x},t)$, $\mathit{\Phi}^{\varepsilon}(\bm{x},t)$ and ${U_i^{\varepsilon} }(\bm{x},t)$ are approximated via two-scale asymptotic expansion in the forms inspired by \cite{R14,R15}.
\begin{equation}
\left\{\begin{array}{l}
{T^{\varepsilon} }(\bm{x},t)={T_{0}}(\bm{x},\bm{y},t) + \varepsilon {T_{1}}(\bm{x},\bm{y},t) + {\varepsilon ^2}{T_{2}}(\bm{x},\bm{y},t) + {\rm O}({\varepsilon ^3}),\\
\mathit{\Phi}^{\varepsilon}(\bm{x},t)=\mathit{\Phi}_{0}(\bm{x},\bm{y},t) + \varepsilon \mathit{\Phi}_{1}(\bm{x},\bm{y},t) + {\varepsilon ^2}\mathit{\Phi}_{2}(\bm{x},\bm{y},t) + {\rm O}({\varepsilon ^3}),\\
{U_i^{\varepsilon} }(\bm{x},t)={U_{i0}}(\bm{x},\bm{y},t) + \varepsilon {U_{i1}}(\bm{x},\bm{y},t) + {\varepsilon ^2}{U_{i2}}(\bm{x},\bm{y},t) + {\rm O}({\varepsilon ^3}).
\end{array}\right.
\end{equation}
Among the formula (2.1), the expansion terms with subscript $i$ for $i=0,1,2$ represent the $i$-order expansion terms. Notably, first-order and second-order asymptotic terms are also referred to as low-order and high-order asymptotic terms respectively.

Next, the pivotal strategy in tackling the temperature-dependent material parameters is introduced for multi-scale nonlinear system (1.1). Employing the temperature-dependent material parameter $c^\varepsilon ({\bm{x}},T^{\varepsilon})$ as an example, it can be expanded around leading temperature term $T_0$ by displacing ${x_0} = \bm{y}$, ${y_0} = {T_0}$ and $\delta = \varepsilon {T_1} + {\varepsilon^2}{T_2} + O({\varepsilon^3})$ in multi-index Taylor's formula $\displaystyle f({x_0},{y_0} + \delta)=f({x_0},{y_0}) + \mathbf{D}^{(0,1)}{f}({x_0},{y_0})\delta + \frac{1}{2}\mathbf{D}^{(0,2)}{f}({x_0},{y_0}){\delta^2} + O({\delta^3})$ with $\mathbf{D}^{(0,1)}{f}({x_0},{y_0})={f_y}({x_0},{y_0})$ and $\mathbf{D}^{(0,2)}{f}({x_0},{y_0})={f_{yy}}({x_0},{y_0})$ as below
\begin{equation}
\begin{aligned}
&\quad\; c^\varepsilon ({\bm{x}},T^{\varepsilon})= {c}({\bm{y}},T^{\varepsilon}) = {c}({\bm{y}},{T_{0}} + \varepsilon{T_{1}} + {\varepsilon^2}{T_{2}} + {\rm O}({\varepsilon^3}))\\
& = {c}({\bm{y}},{T_{0}}) + \mathbf{D}^{(0,1)}{c}({\bm{y}},{T_{0}})\big[ {\varepsilon{T_{1}} + {\varepsilon^2}{T_{2}} + {\rm O}({\varepsilon^3})} \big]\\
&+ \frac{1}{2}\mathbf{D}^{(0,2)}{c}({\bm{y}},{T_{0}}){\big[ {\varepsilon{T_{1}} + {\varepsilon^2}{T_{2}} + {\rm O}({\varepsilon^3})} \big]^2} + {\rm O}\Big( \big[ {\varepsilon{T_{1}} + {\varepsilon^2}{T_{2}} + {\rm O}({\varepsilon^3})} \big]^3 \Big)\\
& = {c}({\bm{y}},{T_{0}}) + \varepsilon{T_{1}}\mathbf{D}^{(0,1)}{c}({\bm{y}},{T_{0}})\\
&+ {\varepsilon^2}\big[ {{T_{2}}\mathbf{D}^{(0,1)}{c}({\bm{y}},{T_{0}}) + \frac{1}{2}{{( {{T_{1}}} )}^2}\mathbf{D}^{(0,2)}{c}({\bm{y}},{T_{0}})} \big] + {\rm O}({\varepsilon^3})\\
&=:c^{(0)}({\bm{y}},{T_{0}})+ \varepsilon c^{(1)}({\bm{x}},{\bm{y}},{T_{0}})+ {\varepsilon^2}c^{(2)}({\bm{x}},{\bm{y}},{T_{0}})+ {\rm O}({\varepsilon^3}).
\end{aligned}
\end{equation}
Next, the remaining temperature-dependent material property parameters $\rho^{\varepsilon}({\bm{x}},{T^{\varepsilon}})$, $k_{ij}^{\varepsilon}({\bm{x}},{T^{\varepsilon}})$, ${\lambda_{ij}^{\varepsilon}}({\bm{x}},{T^{\varepsilon}})$, ${\beta_{ij}^{\varepsilon}}({\bm{x}},{T^{\varepsilon}})$ and ${c_{ijkl}^{\varepsilon}}({\bm{x}},{T^{\varepsilon}})$ successively possess the expansion forms as follows
\begin{equation}
\begin{aligned}
&{\rho^{\varepsilon}}({\bm{x}},{T^{\varepsilon}})= {\rho^{(0)}}({\bm{y}},{T_{0}}) + \varepsilon{\rho^{(1)}}({\bm{x}},{\bm{y}},{T_{0}}) + {\varepsilon^2}{\rho^{(2)}}({\bm{x}},{\bm{y}},{T_{0}}) + {\rm O}({\varepsilon^3}),\\
&k_{ij}^\varepsilon ({\bm{x}},T^{\varepsilon})=k_{ij}^{(0)}({\bm{y}},{T_{0}})+ \varepsilon k_{ij}^{(1)}({\bm{x}},{\bm{y}},{T_{0}})+ {\varepsilon^2}k_{ij}^{(2)}({\bm{x}},{\bm{y}},{T_{0}})+ {\rm O}({\varepsilon^3}),\\
&\lambda_{ij}^{\varepsilon}({\bm{x}},{T^{\varepsilon}})= \lambda_{ij}^{(0)}({\bm{y}},{T_{0}}) + \varepsilon\lambda_{ij}^{(1)}({\bm{x}},{\bm{y}},{T_{0}}) + {\varepsilon^2}\lambda_{ij}^{(2)}({\bm{x}},{\bm{y}},{T_{0}}) + {\rm O}({\varepsilon ^3}),\\
&\beta_{ij}^{\varepsilon}({\bm{x}},{T^{\varepsilon}})= \beta_{ij}^{(0)}({\bm{y}},{T_{0}}) + \varepsilon\beta_{ij}^{(1)}({\bm{x}},{\bm{y}},{T_{0}}) + {\varepsilon^2}\beta_{ij}^{(2)}({\bm{x}},{\bm{y}},{T_{0}}) + {\rm O}({\varepsilon ^3}),\\
&c_{ijkl}^{\varepsilon}({\bm{x}},{T^{\varepsilon}})= c_{ijkl}^{(0)}({\bm{y}},{T_{0}}) + \varepsilon c_{ijkl}^{(1)}({\bm{x}},{\bm{y}},{T_{0}}) + {\varepsilon^2}c_{ijkl}^{(2)}({\bm{x}},{\bm{y}},{T_{0}}) + {\rm O}({\varepsilon ^3}).
\end{aligned}
\end{equation}

Under the assumption of periodic heterogeneity of composite structures, the chain rule yields the following result for multi-scale asymptotic analysis.
\begin{equation}
\displaystyle\frac{\partial \Theta^{\varepsilon} (\bm{x},t)}{\partial x_i}=\frac{\partial \Theta (\bm{x},\bm{y},t)}{\partial x_i}+\frac{1}{\varepsilon}\frac{\partial \Theta (\bm{x},\bm{y},t)}{\partial y_i},\;\;(i=1,\cdots,n),
\end{equation}
where $\Theta^{\varepsilon}(\bm{x},t)$ represents a function exhibiting multi-scale features, specifically defined as $\Theta^{\varepsilon}(\bm{x},t)=\Theta(\bm{x},\bm{y},t)$; see reference \cite{R14}.

Then substituting (2.1)-(2.3) into the nonlinear initial-boundary problem (1.1) with two-scale features and applying the chain rule (2.4), a sequence of equations is derived by collecting terms of the same order $\varepsilon^{-2}$ as follows.
\begin{equation}
\begin{aligned}
{\rm O}({\varepsilon^{ - 2}}):\frac{\partial }{{\partial {y_i}}}\Big( {{k_{ij}^{(0)}}\frac{\partial T_{0}}{\partial {y_j}}}\Big)=-{\lambda_{ij}^{(0)}}\frac{\partial \mathit{\Phi}_{0}}{\partial {y_i}}\frac{\partial \mathit{\Phi}_{0}}{\partial {y_j}},
\end{aligned}
\end{equation}
\begin{equation}
\begin{aligned}
{\rm O}({\varepsilon^{ - 2}}):\frac{\partial }{{\partial {y_i}}}\Big({{\lambda_{ij}^{(0)}}\frac{\partial \mathit{\Phi}_{0}}{\partial {y_j}}}\Big)=0,
\end{aligned}
\end{equation}
\begin{equation}
\begin{aligned}
{\rm O}({\varepsilon^{ - 2}}):\frac{\partial }{{\partial {y_j}}}\Big( {{c_{ijkl}^{(0)}}\frac{{\partial {U_{k0}}}}{{\partial {y_l}}}}\Big)=0.
\end{aligned}
\end{equation}
Analogously, following equations are determined by equating the coefficient $\varepsilon^{-1}$.
\begin{equation}
\begin{aligned}
{\rm O}({\varepsilon^{ - 1}}):
&\frac{\partial }{{\partial {y_i}}}\Big({{k_{ij}^{(0)}}\frac{\partial T_{0}}{\partial {x_j}}}\Big)\!+\!\frac{\partial }{{\partial {x_i}}}\Big({{k_{ij}^{(0)}}\frac{\partial T_{0}}{\partial {y_j}}}\Big)\!+\!\frac{\partial }{{\partial {y_i}}}\Big({{k_{ij}^{(0)}}\frac{\partial T_{1}}{\partial {y_j}}}\Big)\!+\!\frac{\partial }{{\partial {y_i}}}\Big({{k_{ij}^{(1)}}\frac{\partial T_{0}}{\partial {y_j}}}\Big)\\
&=-{\lambda_{ij}^{(0)}}\frac{\partial \mathit{\Phi}_{0}}{\partial {x_i}}\frac{\partial \mathit{\Phi}_{0}}{\partial {y_j}}-{\lambda_{ij}^{(0)}}\frac{\partial \mathit{\Phi}_{0}}{\partial {y_i}}\frac{\partial \mathit{\Phi}_{0}}{\partial {x_j}}-{\lambda_{ij}^{(1)}}\frac{\partial \mathit{\Phi}_{0}}{\partial {y_i}}\frac{\partial \mathit{\Phi}_{0}}{\partial {y_j}}\\
&-{\lambda_{ij}^{(0)}}\frac{\partial \mathit{\Phi}_{0}}{\partial {y_i}}\frac{\partial \mathit{\Phi}_{1}}{\partial {y_j}}-{\lambda_{ij}^{(0)}}\frac{\partial \mathit{\Phi}_{1}}{\partial {y_i}}\frac{\partial \mathit{\Phi}_{0}}{\partial {y_j}}+T_0\beta_{ij}^{(0)}\frac{\partial}{\partial t}\Big(\frac{\partial U_{i0}}{\partial y_j}\Big),
\end{aligned}
\end{equation}
\begin{equation}
\begin{aligned}
{\rm O}({\varepsilon^{ - 1}}):
&\frac{\partial }{{\partial {y_i}}}\Big( {{\lambda_{ij}^{(0)}}\frac{\partial \mathit{\Phi}_{0}}{\partial {x_j}}}\Big)+\frac{\partial }{{\partial {x_i}}}\Big({{\lambda_{ij}^{(0)}}\frac{\partial \mathit{\Phi}_{0}}{\partial {y_j}}}\Big)\\
&+\frac{\partial }{{\partial {y_i}}}\Big({{\lambda_{ij}^{(0)}}\frac{\partial \mathit{\Phi}_{1}}{\partial {y_j}}}\Big)+\frac{\partial }{{\partial {y_i}}}\Big({{\lambda_{ij}^{(1)}}\frac{\partial \mathit{\Phi}_{0}}{\partial {y_j}}}\Big)=0,
\end{aligned}
\end{equation}
\begin{equation}
\begin{aligned}
{\rm O}({\varepsilon^{ - 1}}):
&\frac{\partial }{{\partial {y_j}}}\Big({c_{ijkl}^{(0)}}\frac{\partial {U_{k0}}}{{\partial {x_l}}}\Big)+\frac{\partial }{{\partial {x_j}}}\Big( {{c_{ijkl}^{(0)}}\frac{{\partial {U_{k0}}}}{{\partial {y_l}}}}\Big)+\frac{\partial }{{\partial {y_j}}}\Big({c_{ijkl}^{(0)}}\frac{\partial {U_{k1}}}{{\partial {y_l}}}\Big)\\
&+\frac{\partial }{{\partial {y_j}}}\Big({c_{ijkl}^{(1)}}\frac{\partial {U_{k0}}}{{\partial {y_l}}}\Big)=\frac{\partial }{{\partial {y_j}}}\Big(\beta_{ij}^{(0)}(T_0-\widetilde T)\Big).
\end{aligned}
\end{equation}
Meanwhile, taking into account the coefficient $\varepsilon^{0}$ yields the following equations.
\begin{equation}
\begin{aligned}
{\rm O}({\varepsilon^0}):
&\rho^{(0)}c^{(0)}\frac{{\partial {T_{0}}}}{{\partial t}}-\frac{\partial }{{\partial {x_i}}}\Big( {{k_{ij}^{(0)}}\frac{{\partial T_{0} }}{{\partial {x_j}}}} \Big)-\frac{\partial }{{\partial {x_i}}}\Big( {{k_{ij}^{(0)}}\frac{{\partial T_{1} }}{{\partial {y_j}}}} \Big)-\frac{\partial }{{\partial {y_i}}}\Big({{k_{ij}^{(0)}}\frac{{\partial T_{1} }}{{\partial {x_j}}}} \Big)\\
&-\frac{\partial }{{\partial {y_i}}}\Big({{k_{ij}^{(0)}}\frac{{\partial T_{2} }}{{\partial {y_j}}}} \Big)-\frac{\partial }{{\partial {x_i}}}\Big({{k_{ij}^{(1)}}\frac{{\partial T_{0} }}{{\partial {y_j}}}} \Big)-\frac{\partial }{{\partial {y_i}}}\Big({{k_{ij}^{(1)}}\frac{{\partial T_{0} }}{{\partial {x_j}}}} \Big)\\
&-\frac{\partial }{{\partial {y_i}}}\Big({{k_{ij}^{(1)}}\frac{{\partial T_{1} }}{{\partial {y_j}}}} \Big)-\frac{\partial }{{\partial {y_i}}}\Big({{k_{ij}^{(2)}}\frac{{\partial T_{0} }}{{\partial {y_j}}}} \Big)={\lambda_{ij}^{(0)}}\frac{\partial \mathit{\Phi}_{0}}{\partial {x_i}}\frac{\partial \mathit{\Phi}_{0}}{\partial {x_j}}+{\lambda_{ij}^{(0)}}\frac{\partial \mathit{\Phi}_{0}}{\partial {x_i}}\frac{\partial \mathit{\Phi}_{1}}{\partial {y_j}}\\
&+{\lambda_{ij}^{(0)}}\frac{\partial \mathit{\Phi}_{1}}{\partial {x_i}}\frac{\partial \mathit{\Phi}_{0}}{\partial {y_j}}+{\lambda_{ij}^{(0)}}\frac{\partial \mathit{\Phi}_{0}}{\partial {y_i}}\frac{\partial \mathit{\Phi}_{1}}{\partial {x_j}}+{\lambda_{ij}^{(0)}}\frac{\partial \mathit{\Phi}_{0}}{\partial {y_i}}\frac{\partial \mathit{\Phi}_{2}}{\partial {y_j}}+{\lambda_{ij}^{(0)}}\frac{\partial \mathit{\Phi}_{1}}{\partial {y_i}}\frac{\partial \mathit{\Phi}_{0}}{\partial {x_j}}\\
&+{\lambda_{ij}^{(0)}}\frac{\partial \mathit{\Phi}_{1}}{\partial {y_i}}\frac{\partial \mathit{\Phi}_{1}}{\partial {y_j}}+{\lambda_{ij}^{(0)}}\frac{\partial \mathit{\Phi}_{2}}{\partial {y_i}}\frac{\partial \mathit{\Phi}_{0}}{\partial {y_j}}+{\lambda_{ij}^{(1)}}\frac{\partial \mathit{\Phi}_{0}}{\partial {x_i}}\frac{\partial \mathit{\Phi}_{0}}{\partial {y_j}}+{\lambda_{ij}^{(1)}}\frac{\partial \mathit{\Phi}_{0}}{\partial {y_i}}\frac{\partial \mathit{\Phi}_{0}}{\partial {x_j}}\\
&+{\lambda_{ij}^{(1)}}\frac{\partial \mathit{\Phi}_{0}}{\partial {y_i}}\frac{\partial \mathit{\Phi}_{1}}{\partial {y_j}}+{\lambda_{ij}^{(1)}}\frac{\partial \mathit{\Phi}_{1}}{\partial {y_i}}\frac{\partial \mathit{\Phi}_{0}}{\partial {y_j}}+{\lambda_{ij}^{(2)}}\frac{\partial \mathit{\Phi}_{0}}{\partial {y_i}}\frac{\partial \mathit{\Phi}_{0}}{\partial {y_j}}-T_0\beta_{ij}^{(0)}\frac{\partial}{\partial t}\Big(\frac{\partial U_{i0}}{\partial x_j}\Big)\\
&-T_0\beta_{ij}^{(0)}\frac{\partial}{\partial t}\Big(\frac{\partial U_{i1}}{\partial y_j}\Big)-T_0\beta_{ij}^{(1)}\frac{\partial}{\partial t}\Big(\frac{\partial U_{i0}}{\partial y_j}\Big)-T_1\beta_{ij}^{(0)}\frac{\partial}{\partial t}\Big(\frac{\partial U_{i0}}{\partial y_j}\Big)+f_T,
\end{aligned}
\end{equation}
\begin{equation}
\begin{aligned}
{\rm O}({\varepsilon^0}):&\frac{\partial }{{\partial {x_i}}}\Big({{\lambda_{ij}^{(0)}}\frac{{\partial \mathit{\Phi}_{0} }}{{\partial {x_j}}}} \Big)+\frac{\partial }{{\partial {x_i}}}\Big( {{\lambda_{ij}^{(0)}}\frac{{\partial \mathit{\Phi}_{1} }}{{\partial {y_j}}}} \Big)+\frac{\partial }{{\partial {y_i}}}\Big({{\lambda_{ij}^{(0)}}\frac{{\partial \mathit{\Phi}_{1} }}{{\partial {x_j}}}} \Big)\\
&+\frac{\partial }{{\partial {y_i}}}\Big({{\lambda_{ij}^{(0)}}\frac{{\partial \mathit{\Phi}_{2} }}{{\partial {y_j}}}} \Big)+\frac{\partial }{{\partial {x_i}}}\Big({{\lambda_{ij}^{(1)}}\frac{{\partial \mathit{\Phi}_{0} }}{{\partial {y_j}}}} \Big)+\frac{\partial }{{\partial {y_i}}}\Big({{\lambda_{ij}^{(1)}}\frac{{\partial \mathit{\Phi}_{0} }}{{\partial {x_j}}}} \Big)\\
&+\frac{\partial }{{\partial {y_i}}}\Big( {{\lambda_{ij}^{(1)}}\frac{{\partial \mathit{\Phi}_{1} }}{{\partial {y_j}}}} \Big)+\frac{\partial }{{\partial {y_i}}}\Big( {{\lambda_{ij}^{(2)}}\frac{{\partial \mathit{\Phi}_{0} }}{{\partial {y_j}}}} \Big)=-f_\mathit{\Phi},
\end{aligned}
\end{equation}
\begin{equation}
\begin{aligned}
{\rm O}({\varepsilon^0}):&\rho^{(0)}\frac{{\partial^2 {U_{i0}}}}{{\partial t^2}}-\frac{\partial }{{\partial {x_j}}}\Big( {{c_{ijkl}^{(0)}}\frac{{\partial {U_{k0}}}}{{\partial {x_l}}}}\Big)-\frac{\partial }{{\partial {x_j}}}\Big( {{c_{ijkl}^{(0)}}\frac{{\partial {U_{k1}}}}{{\partial {y_l}}}}\Big)\\
&-\frac{\partial }{{\partial {y_j}}}\Big( {{c_{ijkl}^{(0)}}\frac{{\partial {U_{k1}}}}{{\partial {x_l}}}}\Big)-\frac{\partial }{{\partial {y_j}}}\Big( {{c_{ijkl}^{(0)}}\frac{{\partial {U_{k2}}}}{{\partial {y_l}}}}\Big)-\frac{\partial }{{\partial {x_j}}}\Big( {{c_{ijkl}^{(1)}}\frac{{\partial {U_{k0}}}}{{\partial {y_l}}}}\Big)\\
&-\frac{\partial }{{\partial {y_j}}}\Big( {{c_{ijkl}^{(1)}}\frac{{\partial {U_{k0}}}}{{\partial {x_l}}}}\Big)-\frac{\partial }{{\partial {y_j}}}\Big( {{c_{ijkl}^{(1)}}\frac{{\partial {U_{k1}}}}{{\partial {y_l}}}}\Big)-\frac{\partial }{{\partial {y_j}}}\Big( {{c_{ijkl}^{(2)}}\frac{{\partial {U_{k0}}}}{{\partial {y_l}}}}\Big)\\
&=-\frac{\partial }{{\partial {x_j}}}\Big(\beta_{ij}^{(0)}(T_0-\widetilde T)\Big)-\frac{\partial }{{\partial {y_j}}}\Big(\beta_{ij}^{(1)}(T_0-\widetilde T)\Big)-\frac{\partial }{{\partial {y_j}}}\Big(\beta_{ij}^{(0)}T_1\Big)+f_i.
\end{aligned}
\end{equation}

Next, applying the standard multi-scale asymptotic expansion procedure to the equations of order $O({\varepsilon^{-2}})$, $O({\varepsilon^{- 1}})$ and $O({\varepsilon^{0}})$, the low-order multi-scale (LOMS) solutions for temperature $T^{\varepsilon}$, electric potential $\mathit{\Phi}^{\varepsilon}$ and displacement $U_i^{\varepsilon}$ are derived as follows
\begin{equation}
{\begin{aligned}
{T^{[1\varepsilon]} }(\bm{x},t)&=
T_{0}(\bm{x},t)+\varepsilon{\mathcal{M}_{\alpha_1}}(\bm{y},T_{0})\frac{\partial T_{0}(\bm{x},t)}{\partial x_{\alpha_1}},
\end{aligned}}
\end{equation}
\begin{equation}
{\begin{aligned}
\mathit{\Phi}^{[1\varepsilon]}(\bm{x},t)&=
\mathit{\Phi}_{0}(\bm{x},t)+\varepsilon \mathcal{H}_{\alpha_1}(\bm{y},T_{0})\frac{\partial \mathit{\Phi}_{0}(\bm{x},t)}{\partial x_{\alpha_1}},
\end{aligned}}
\end{equation}
\begin{equation}
{\begin{aligned}
U_i^{[1\varepsilon]}(\bm{x},t)&=
U_{i0}(\bm{x},t)\!+\!\varepsilon\Big[ \mathcal{N}_{im}^{\alpha_1}(\bm{y},T_{0})\frac{\partial U_{m0}(\bm{x},t)}{\partial x_{\alpha_1}}\!+\!\mathcal{P}_{i}(\bm{y},T_{0})(T_0(\bm{x},t)\!-\!\widetilde T)\Big],
\end{aligned}}
\end{equation}
where $\mathcal{M}_{\alpha_1}$, $\mathcal{H}_{\alpha_1}$, $\mathcal{N}_{im}^{\alpha_1}$ and $\mathcal{P}_{i}$ denote the first-order cell functions defined on PUC $\mathbf{Y}$, all of which depend on the macroscopic temperature field $T_0$.

To effectively improve the computational accuracy, the novel high-order multi-scale (HOMS) solutions for temperature $T^{\varepsilon}$, electric potential $\mathit{\Phi}^{\varepsilon}$ and displacement $U_i^{\varepsilon}$ are further established as follows
\begin{equation}
{\begin{aligned}
{T^{[2\varepsilon]} }(\bm{x},t)&\!\!=\!\!
T_{0}(\bm{x},\!t)+\varepsilon{\mathcal{M}_{\alpha_1}}(\bm{y},T_{0})\frac{\partial T_{0}(\bm{x},t)}{\partial x_{\alpha_1}}+\varepsilon^2\big[\mathcal{Q}({\bm{y}},{T_{0}})\frac{{\partial {T_{0}}(\bm{x},t)}}{{\partial t}}\\
&\!\!+\!\!{\mathcal{M}_{\alpha_1\alpha_2}}({\bm{y}},{T_{0}})\frac{\partial^2 T_{0}(\bm{x},t)}{\partial x_{\alpha_1}\partial x_{\alpha_2}}+\!\mathcal{R}_{\alpha_1}({\bm{y}},{T_{0}})\frac{\partial T_{0}(\bm{x},t)}{\partial x_{\alpha_1}}\\
&\!\!+\!\!{\mathcal{O}_{\alpha_1\alpha_2}}({\bm{y}},\!{T_{0}})\frac{\partial T_{0}(\bm{x},t)}{\partial x_{\alpha_1}}\frac{\partial T_{0}(\bm{x},t)}{\partial x_{\alpha_2}}+{\mathcal{G}_{\alpha_1\alpha_2}}({\bm{y}},\!{T_{0}})\frac{\partial \mathit{\Phi}_{0}(\bm{x},t)}{\partial x_{\alpha_1}}\frac{\partial \mathit{\Phi}_{0}(\bm{x},t)}{\partial x_{\alpha_2}}\\
&\!\!+{\mathcal{J}_{\alpha_1\alpha_2}}({\bm{y}},{T_{0}})\frac{\partial^2 U_{\alpha_10}(\bm{x},t)}{\partial x_{\alpha_2}\partial t}\big],
\end{aligned}}
\end{equation}
\begin{equation}
{\begin{aligned}
\mathit{\Phi}^{[2\varepsilon]}(\bm{x},t)&=
\mathit{\Phi}_{0}(\bm{x},t)+\varepsilon \mathcal{H}_{\alpha_1}(\bm{y},T_{0})\frac{\partial \mathit{\Phi}_{0}(\bm{x},t)}{\partial x_{\alpha_1}}+ \varepsilon^2\big[\mathcal{H}_{\alpha_1\alpha_2}({\bm{y}},{T_{0}})\frac{\partial^2 \mathit{\Phi}_{0}(\bm{x},t)}{\partial x_{\alpha_1}\partial x_{\alpha_2}}\\
&+ \mathcal{Z}_{\alpha_1}({\bm{y}},{T_{0}})\frac{\partial \mathit{\Phi}_{0}(\bm{x},t)}{\partial x_{\alpha_1}}+\mathcal{W}_{\alpha_1\alpha_2}({\bm{y}},{T_{0}})\frac{\partial T_{0}(\bm{x},t)}{\partial x_{\alpha_1}}\frac{\partial \mathit{\Phi}_{0}(\bm{x},t)}{\partial x_{\alpha_2}}\big],
\end{aligned}}
\end{equation}
\begin{equation}
{\begin{aligned}
U_i^{[2\varepsilon]}(\bm{x},t)&=
U_{i0}(\bm{x},t)\!+\!\varepsilon\Big[\mathcal{N}_{im}^{\alpha_1}(\bm{y},T_{0})\frac{\partial U_{m0}(\bm{x},t)}{\partial x_{\alpha_1}}\!+\!\mathcal{P}_{i}(\bm{y},T_{0})(T_0(\bm{x},t)\!-\!\widetilde T)\Big]\\
&+\varepsilon^2\big[\mathcal{N}_{im}^{\alpha_1\alpha_2}(\bm{y},T_{0})\frac{\partial^2 U_{m0}(\bm{x},t)}{\partial x_{\alpha_1}\partial x_{\alpha_2}}+\mathcal{F}_{i}^{\alpha_1}(\bm{y},T_{0})\frac{{\partial^2 {U_{\alpha_10}(\bm{x},t)}}}{{\partial t^2}}\\
&+\mathcal{X}_i^{\alpha_1}({\bm{y}},{T_{0}})\frac{\partial T_{0}(\bm{x},t)}{\partial x_{\alpha_1}}+\mathcal{A}_{im}^{\alpha_1}(\bm{y},T_{0})\frac{\partial U_{m0}(\bm{x},t)}{\partial x_{\alpha_1}}\\
&+\mathcal{B}_{i}(\bm{y},T_{0})(T_0(\bm{x},t)-\widetilde T)+\mathcal{C}_{i}^{\alpha_1}(\bm{y},T_{0})\frac{\partial T_{0}(\bm{x},t)}{\partial x_{\alpha_1}}(T_0(\bm{x},t)-\widetilde T)\\
&+{\mathcal{D}_{im}^{\alpha_1\alpha_2}}({\bm{y}},{T_{0}})\frac{\partial T_{0}(\bm{x},t)}{\partial x_{\alpha_1}}\frac{\partial U_{m0}(\bm{x},t)}{\partial x_{\alpha_2}}\big],
\end{aligned}}
\end{equation}
where $\mathcal{Q}$, $\mathcal{M}_{\alpha_1\alpha_2}$, $\mathcal{R}_{\alpha_1}$, $\mathcal{O}_{\alpha_1\alpha_2}$, $\mathcal{G}_{\alpha_1\alpha_2}$, $\mathcal{J}_{\alpha_1\alpha_2}$, $\mathcal{H}_{\alpha_1\alpha_2}$, $\mathcal{Z}_{\alpha_1}$, $\mathcal{W}_{\alpha_1\alpha_2}$, $\mathcal{N}_{im}^{\alpha_1\alpha_2}$, $\mathcal{F}_{i}^{\alpha_1}$, $\mathcal{X}_{i}^{\alpha_1}$, $\mathcal{A}_{im}^{\alpha_1}$, $\mathcal{B}_{i}$, $\mathcal{C}_{i}^{\alpha_1}$ and $\mathcal{D}_{im}^{\alpha_1\alpha_2}$ denote the second-order cell functions defined on PUC $\mathbf{Y}$, all of which depend on the macroscopic temperature field $T_0$. Furthermore, it is worth emphasizing that only the HOMS solutions (2.17) of temperature field can characterize the mutual coupling impact of electric potential field and displacement field on temperature field owing to the presence of correction terms $\displaystyle{\mathcal{G}_{\alpha_1\alpha_2}}\frac{\partial \mathit{\Phi}_{0}}{\partial x_{\alpha_1}}\frac{\partial \mathit{\Phi}_{0}}{\partial x_{\alpha_2}}$ and $\displaystyle{\mathcal{J}_{\alpha_1\alpha_2}}\frac{\partial^2 U_{\alpha_10}}{\partial x_{\alpha_2}\partial t}$, while only the HOMS solutions (2.18) of electric potential field can characterize the mutual coupling impact of temperature field on electric potential field owing to the presence of correction terms $\displaystyle\mathcal{W}_{\alpha_1\alpha_2}\frac{\partial T_{0}}{\partial x_{\alpha_1}}\frac{\partial \mathit{\Phi}_{0}}{\partial x_{\alpha_2}}$. This is one essential trigger to develop high-order multi-scale method for high-accuracy nonlinear  thermo-electro-mechanical coupling simulation of composite structures.

Inspired by classical multi-scale asymptotic expansion theory, the detailed definitions of auxiliary cell functions can be obtained. Firstly, leveraging the $O({\varepsilon^{-1}})$-order equations (2.8)-(2.10), the following boundary-value problems subject to periodic boundary condition are derived for solving first-order cell functions, which are designated as the first-order cell problems.
\begin{equation}
\left\{
\begin{aligned}
&\frac{\partial}{\partial y_i}\big[ { k_{ij}^{(0)}{\frac{\partial \mathcal{M}_{\alpha_1}}{\partial y_j}}} \big]= -\frac{\partial k_{i{\alpha_1}}^{(0)}}{\partial y_i},\;\;\;\bm{y}\in \mathbf{Y}, \\
&\mathcal{M}_{\alpha_1}(\bm{y},T_{0})\;\mathrm{is}\;1-\mathrm{periodic}\;\mathrm{in}\;{\bm{y}}{\rm{, }}\;\;\;{\int_{\mathbf{Y}}}\mathcal{M}_{\alpha_1}d\mathbf{Y}=0.
\end{aligned} \right.
\end{equation}
\begin{equation}
\left\{
\begin{aligned}
&\frac{\partial}{\partial y_i}\big[ { \lambda_{ij}^{(0)}{\frac{\partial \mathcal{H}_{\alpha_1}}{\partial y_j}}} \big]= -\frac{\partial \lambda_{i{\alpha_1}}^{(0)}}{\partial y_i},\;\;\;\bm{y}\in \mathbf{Y}, \\
&\mathcal{H}_{\alpha_1}(\bm{y},T_{0})\;\mathrm{is}\;1-\mathrm{periodic}\;\mathrm{in}\;{\bm{y}}{\rm{, }}\;\;\;{\int_{\mathbf{Y}}}\mathcal{H}_{\alpha_1}d\mathbf{Y}=0.
\end{aligned} \right.
\end{equation}
\begin{equation}
\left\{
\begin{aligned}
&\frac{\partial}{\partial y_j}\big[ { c_{ijkl}^{(0)}{\frac{\partial \mathcal{N}_{km}^{\alpha_1}}{\partial y_l}}} \big]= -\frac{\partial c_{ijm{\alpha_1}}^{(0)}}{\partial y_j},\;\;\;\bm{y}\in \mathbf{Y}, \\
&\mathcal{N}_{km}^{\alpha_1}(\bm{y},T_{0})\;\mathrm{is}\;1-\mathrm{periodic}\;\mathrm{in}\;{\bm{y}}{\rm{, }}\;\;\;{\int_{\mathbf{Y}}}\mathcal{N}_{km}^{\alpha_1}d\mathbf{Y}=0.
\end{aligned} \right.
\end{equation}
\begin{equation}
\left\{
\begin{aligned}
&\frac{\partial}{\partial y_j}\big[{c_{ijkl}^{(0)}{\frac{\partial \mathcal{P}_{k}}{\partial y_l}}} \big]= \frac{\partial \beta_{ij}^{(0)}}{\partial y_j},\;\;\;\bm{y}\in \mathbf{Y}, \\
&\mathcal{P}_{k}(\bm{y},T_{0})\;\mathrm{is}\;1-\mathrm{periodic}\;\mathrm{in}\;{\bm{y}}{\rm{, }}\;\;\;{\int_{\mathbf{Y}}}\mathcal{P}_{k}d\mathbf{Y}=0.
\end{aligned} \right.
\end{equation}

Subsequently, integrating the $O({\varepsilon^{0}})$-order equations (2.11)-(2.13) over microscopic unit cell $\mathcal{Y}$ (with volume integration) and applying the Gauss theorem to these equations, the macroscopic homogenized equations associated with the multi-scale nonlinear coupling problem (1.1) are derived, which are presented below
\begin{equation}
\left\{ \begin{aligned}
&\widehat S(T_{0})\frac{{\partial {T_{0}}}}{{\partial t}}- \frac{\partial }{{\partial {x_i}}}\Big( {{\widehat k_{ij}}(T_{0})\frac{{\partial {T_{0}}}}{{\partial {x_j}}}}\Big)\\
&\quad\quad= {\widehat \lambda_{ij}^*}(T_{0})\frac{\partial \mathit{\Phi}_{0}}{\partial {x_i}}\frac{\partial \mathit{\Phi}_{0}}{\partial {x_j}}-T_0\widehat \beta_{ij}^*(T_{0})\frac{\partial}{\partial t}\Big(\frac{\partial U_{i0}}{\partial x_j}\Big)+f_T,\;\text{in}\;\Omega\times(0,\mathcal T),\\
&- \frac{\partial }{{\partial {x_i}}}\Big( {{\widehat \lambda _{ij}}(T_{0})\frac{{\partial \mathit{\Phi}_0}}{{\partial {x_j}}}}\Big)=f_\mathit{\Phi},\;\text{in}\;\Omega\times(0,\mathcal T),\\
&\widehat \rho(T_0)\frac{{\partial^2 {U_{i0}}}}{{\partial t^2}}-\frac{\partial }{{\partial {x_j}}}\Big( {{\widehat c_{ijkl}}(T_0)\frac{{\partial {U_{k0}}}}{{\partial {x_l}}}}\Big)\\
&\quad\quad=-\frac{\partial }{{\partial {x_j}}}\Big(\widehat {\beta}_{ij}(T_0)(T_0-\widetilde T)\Big)+f_i,\;\text{in}\;\Omega\times(0,\mathcal T),\\
&T_0(\bm{x},t) = \widehat T(\bm{x},t),\;\mathit{\Phi}_0(\bm{x},t)=\widehat{\mathit{\Phi}}(\bm{x},t),\;{\bm{U}}_0(\bm{x},t)=\widehat {\bm{U}}(\bm{x},t),\;\text{on}\;\partial\Omega\!\times\!(0,\mathcal T),\\
&T_0({\bm{x}},0)=\widetilde T,\;{\bm{U}}_0({\bm{x}},0)=\widetilde{\bm{U}}^0({\bm{x}}),\;\frac{\partial {\bm{U}}_0({\bm{x}},0)}{\partial t}=\widetilde{\bm{U}}^1({\bm{x}}),\;\text{in}\;\Omega.
\end{aligned} \right.
\end{equation}
Here, the macroscopic homogenized material parameters in (2.24) are computed via the following formulas.
\begin{equation}
\begin{aligned}
&\widehat S(T_{0}) = \frac{1}{|\mathbf{Y}|}{\int_{\mathbf{Y}}}\big({\rho^{(0)}{c}^{(0)}}+T_0\beta_{ij}^{(0)}\frac{\partial \mathcal{P}_{i}}{\partial y_j}\big)d\mathbf{Y},\\
&{\widehat k_{ij}}(T_{0}) = \frac{1}{|\mathbf{Y}|}{\int_{\mathbf{Y}}}\big({k_{ij}^{(0)} + k_{i\alpha_1}^{(0)}{\frac{\partial \mathcal{M}_j}{\partial y_{\alpha_1}}}}\big)d\mathbf{Y},\\
&{\widehat \lambda _{ij}}(T_{0}) = \frac{1}{|\mathbf{Y}|}{\int_{\mathbf{Y}}}\big({\lambda_{ij}^{(0)}+ \lambda_{i\alpha_1}^{(0)}{\frac{\partial \mathcal{H}_j}{\partial y_{\alpha_1}}}}\big)d\mathbf{Y},\\
&\widehat \rho(T_{0}) = \frac{1}{|\mathbf{Y}|}{\int_{\mathbf{Y}}}{\rho^{(0)}}d\mathbf{Y},\\
&{\widehat c_{ijkl}}(T_{0}) = \frac{1}{|\mathbf{Y}|}{\int_{\mathbf{Y}}}\big({c_{ijkl}^{(0)} + c_{ij\alpha_1\alpha_2}^{(0)}{\frac{\partial \mathcal{N}_{\alpha_1k}^{l}}{\partial y_{\alpha_2}}}}\big)d\mathbf{Y},\\
&{\widehat \beta_{ij}}(T_{0}) = \frac{1}{|\mathbf{Y}|}{\int_{\mathbf{Y}}}\big({\beta_{ij}^{(0)}- c_{ijkl}^{(0)}{\frac{\partial \mathcal{P}_k}{\partial y_l}}}\big)d\mathbf{Y},\\
&{\widehat \lambda _{ij}^*}(T_{0}) = \frac{1}{|\mathbf{Y}|}{\int_{\mathbf{Y}}}\big({\lambda_{ij}^{(0)}+ \lambda_{i\alpha_1}^{(0)}{\frac{\partial \mathcal{H}_j}{\partial y_{\alpha_1}}}}+ {\lambda_{\alpha_1j}^{(0)}{\frac{\partial \mathcal{H}_i}{\partial y_{\alpha_1}}}}+ \lambda_{\alpha_1\alpha_2}^{(0)}{\frac{\partial \mathcal{H}_i}{\partial y_{\alpha_1}}{\frac{\partial \mathcal{H}_j}{\partial y_{\alpha_2}}}}\big)d\mathbf{Y},\\
&{\widehat \beta_{ij}^*}(T_{0}) = \frac{1}{|\mathbf{Y}|}{\int_{\mathbf{Y}}}\big({\beta_{ij}^{(0)}+ \beta_{\alpha_1\alpha_2}^{(0)}{\frac{\partial \mathcal{N}_{\alpha_1i}^j}{\partial y_{\alpha_2}}}}\big)d\mathbf{Y},
\end{aligned}
\end{equation}
where $|\mathbf{Y}|$ denotes the Lebesgue measure of microscopic unit cell $\mathbf{Y}$.
\begin{rmk}
All macroscopic homogenized material parameters of the investigated composites obviously vary with the macroscopic temperature solution $T_{0}$ because of the quasi-periodic properties of first-order cell functions with respect to $T_{0}$. This feature yields a significant difference relative to linear composites.
\end{rmk}
\begin{rmk}
Although ${\widehat \lambda _{ij}}(T_{0})$ and ${\widehat \lambda _{ij}^*}(T_{0})$, as well as ${\widehat \beta_{ij}}(T_{0})$ and ${\widehat \lambda _{ij}^*}(T_{0})$ are defined by formally distinct computational formulas, it can be proved that ${\widehat \lambda _{ij}}(T_{0})={\widehat \lambda _{ij}^*}(T_{0})$ and ${\widehat \beta_{ij}}(T_{0})={\widehat \lambda _{ij}^*}(T_{0})$ hold for any fixed $T_{0}$, with details provided in Appendix A. This result is of considerable interest to physics community.
\end{rmk}
\begin{rmk}
Following the approach in \cite{R14,R15}, $\bar{S}_0\leq \widehat S(T_0)\le\bar{S}_1$, $\bar{\kappa}_0|\bm{\xi}|^2\leq \widehat k_{ij}(T_0)\xi_i\xi_j\le\bar{\kappa}_1|\bm{\xi}|^2$, $\bar{\kappa}_0|\bm{\xi}|^2\leq \widehat \lambda_{ij}(T_0)\xi_i\xi_j\le\bar{\kappa}_1|\bm{\xi}|^2$, $\bar{\kappa}_0|\bm{\xi}|^2\leq \widehat \beta_{ij}(T_0)\xi_i\xi_j\le\bar{\kappa}_1|\bm{\xi}|^2$ and $\bar\kappa_0|\bm{\zeta}|^2\leq {\widehat c_{ijkl}}(T_{0})\zeta_{ij}\zeta_{kl}\le\bar\kappa_1|\bm{\zeta}|^2$ can be proved, where $\bar{S}_0$, $\bar{S}_1$, $\bar{\kappa}_0$ and $\bar{\kappa}_1$ are positive constants.
\end{rmk}

Furthermore, substituting the terms $f_T(\bm{x},t)$, $f_\mathit{\Phi}(\bm{x},t)$ and $f_i$ into the $O({\varepsilon^{0}})$-order equations (2.11)-(2.13) with their counterparts in the macroscopic homogenized equations (2.24), the following boundary-value problems subject to the periodic boundary condition are formulated for solving second-order cell functions, which are designated as the second-order cell problems.
\begin{equation}
\left\{
\begin{aligned}
&\frac{\partial}{\partial y_i}\big[ { k_{ij}^{(0)}{\frac{\partial \mathcal{Q}}{\partial y_j}}} \big] = { \rho^{(0)} {c}^{(0)}-\widehat S }+T_0\beta_{ij}^{(0)}\frac{\partial \mathcal{P}_{i}}{\partial y_j},\;\;\;\bm{y}\in \mathbf{Y},\\
&\mathcal{Q}({\bm{y}},{T_{0}})\;\mathrm{is}\;1-\mathrm{periodic}\;\mathrm{in}\;{\bm{y}}{\rm{, }}\;\;\;{\int_{\mathbf{Y}}}\mathcal{Q}d\mathbf{Y}=0.
\end{aligned} \right.
\end{equation}
\begin{equation}
\left\{
\begin{aligned}
&\frac{\partial}{\partial y_i}\big[ { k_{ij}^{(0)}{\frac{\partial {\mathcal{M}_{\alpha_1\alpha_2}}}{\partial y_j}}} \big]= { {\widehat k}_{\alpha_1\alpha_2}\!\!-\! \! {k_{\alpha_1\alpha_2}^{(0)}}\!-\! \frac{\partial}{\partial y_i}\big( {k_{i\alpha_1}^{(0)}{\mathcal{M}_{\alpha_2}}} \big)\!\!-{k_{\alpha_1j}^{(0)}\frac{\partial \mathcal{M}_{\alpha_2}}{\partial y_j}}},\;\;\;\bm{y}\in \mathbf{Y},\\
&{\mathcal{M}_{\alpha_1\alpha_2}}({\bm{y}},{T_{0}})\;\mathrm{is}\;1-\mathrm{periodic}\;\mathrm{in}\;{\bm{y}}{\rm{, }}\;\;\;{\int_{\mathbf{Y}}}\mathcal{M}_{\alpha_1\alpha_2}d\mathbf{Y}=0.
\end{aligned} \right.
\end{equation}
\begin{equation}
\left\{
\begin{aligned}
&\frac{\partial}{\partial y_i}\big[ { k_{ij}^{(0)}{\frac{\partial \mathcal{R}_{\alpha_1}}{\partial y_j}}} \big]\!\!=\!\!{ \frac{\partial{\widehat k}_{i\alpha_1}}{\partial x_i} \!\!-\!\!  {\frac{\partial{k}_{i\alpha_1}^{(0)}}{\partial x_i}}\!\!-\!\!\frac{\partial}{\partial y_i}\big( {k_{ij}^{(0)}\frac{\partial \mathcal{M}_{\alpha_1}}{\partial x_{j}}} \big)\!\!-\!\!\frac{\partial}{\partial x_i}\big( {k_{ij}^{(0)}\frac{\partial \mathcal{M}_{\alpha_1}}{\partial y_{j}}} \big)},\;\bm{y}\!\in\! \mathbf{Y},\\
&{\mathcal{R}_{\alpha_1}}({\bm{y}},{T_{0}})\;\mathrm{is}\;1-\mathrm{periodic}\;\mathrm{in}\;{\bm{y}}{\rm{, }}\;{\int_{\mathbf{Y}}}\mathcal{R}_{\alpha_1}d\mathbf{Y}=0.
\end{aligned} \right.
\end{equation}
\begin{equation}
\left\{
\begin{aligned}
&\frac{\partial}{\partial y_i}\big[ { k_{ij}^{(0)}{\frac{\partial \mathcal{O}_{\alpha_1\alpha_2}}{\partial y_j}}} \big]\!\!=\!\!-\frac{\partial}{\partial y_i}\Big({{\mathcal{M}_{\alpha_1}}\mathbf{D}^{(0,1)}{k_{i\alpha_2}^{(0)}}}\!+ \!{{\mathcal{M}_{\alpha_1}}\mathbf{D}^{(0,1)}{k_{ij}^{(0)}}\frac{\partial \mathcal{M}_{\alpha_2}}{\partial y_j}}\Big),\;\bm{y}\!\in\!\mathbf{Y},\\
&{\mathcal{O}_{\alpha_1\alpha_2}}({\bm{y}},{T_{0}})\;\mathrm{is}\;1-\mathrm{periodic}\;\mathrm{in}\;{\bm{y}}{\rm{, }}\;{\int_{\mathbf{Y}}}\mathcal{O}_{\alpha_1\alpha_2}d\mathbf{Y}=0.
\end{aligned} \right.
\end{equation}
\begin{equation}
\left\{
\begin{aligned}
&\frac{\partial}{\partial y_i}\big[ { k_{ij}^{(0)}{\frac{\partial {\mathcal{G}_{\alpha_1\alpha_2}}}{\partial y_j}}} \big] = {\widehat \lambda}_{\alpha_1\alpha_2}^* - \lambda_{\alpha_1\alpha_2}^{(0)}
-\lambda_{i\alpha_2}^{(0)}\frac{\partial \mathcal{H}_{\alpha_1}}{\partial y_i}\\
&\quad\quad\quad\quad\quad\quad\quad\;\;- \lambda_{\alpha_1j}^{(0)}\frac{\partial \mathcal{H}_{\alpha_2}}{\partial y_{j}}-\lambda_{ij}^{(0)}{\frac{\partial \mathcal{H}_{\alpha_1}}{\partial y_{i}}{\frac{\partial \mathcal{H}_{\alpha_2}}{\partial y_j}}},\;\;\;\bm{y}\in \mathbf{Y},\\
&{\mathcal{G}_{\alpha_1\alpha_2}}({\bm{y}},{T_{0}})\;\mathrm{is}\;1-\mathrm{periodic}\;\mathrm{in}\;{\bm{y}}{\rm{, }}\;\;\;{\int_{\mathbf{Y}}}\mathcal{G}_{\alpha_1\alpha_2}d\mathbf{Y}=0.
\end{aligned} \right.
\end{equation}
\begin{equation}
\left\{
\begin{aligned}
&\frac{\partial}{\partial y_i}\big[ { k_{ij}^{(0)}{\frac{\partial {\mathcal{J}_{\alpha_1\alpha_2}}}{\partial y_j}}} \big] = T_0\beta_{\alpha_1\alpha_2}^{(0)}- T_0{\widehat \beta}_{\alpha_1\alpha_2}^*
+ T_0\beta_{ij}^{(0)}{\frac{\partial \mathcal{N}_{i\alpha_1}^{\alpha_2}}{\partial y_j}},\;\;\;\bm{y}\in \mathbf{Y},\\
&{\mathcal{J}_{\alpha_1\alpha_2}}({\bm{y}},{T_{0}})\;\mathrm{is}\;1-\mathrm{periodic}\;\mathrm{in}\;{\bm{y}}{\rm{, }}\;\;\;{\int_{\mathbf{Y}}}\mathcal{J}_{\alpha_1\alpha_2}d\mathbf{Y}=0.
\end{aligned} \right.
\end{equation}
\begin{equation}
\left\{
\begin{aligned}
&\frac{\partial}{\partial y_i}\big[ { \lambda_{ij}^{(0)}{\frac{\partial {\mathcal{H}_{\alpha_1\alpha_2}}}{\partial y_j}}} \big]\!\!=\!\!{ {\widehat \lambda}_{\alpha_1\alpha_2}\!-\!  {\lambda_{\alpha_1\alpha_2}^{(0)}}\!-\!\frac{\partial}{\partial y_i}\big( {\lambda_{i\alpha_1}^{(0)}{\mathcal{H}_{\alpha_2}}} \big)\!-\!{\lambda_{\alpha_1j}^{(0)}\frac{\partial \mathcal{H}_{\alpha_2}}{\partial y_j}}},\;\;\;\bm{y}\in \mathbf{Y},\\
&{\mathcal{H}_{\alpha_1\alpha_2}}({\bm{y}},{T_{0}})\;\mathrm{is}\;1-\mathrm{periodic}\;\mathrm{in}\;{\bm{y}}{\rm{, }}\;\;\;{\int_{\mathbf{Y}}}\mathcal{H}_{\alpha_1\alpha_2}d\mathbf{Y}=0.
\end{aligned} \right.
\end{equation}
\begin{equation}
\left\{
\begin{aligned}
&\frac{\partial}{\partial y_i}\big[ { \lambda_{ij}^{(0)}{\frac{\partial \mathcal{Z}_{\alpha_1}}{\partial y_j}}} \big]\!\!=\!\!{ \frac{\partial{\widehat \lambda}_{i\alpha_1}}{\partial x_i}\!\!-\!\!  {\frac{\partial{\lambda}_{i\alpha_1}^{(0)}}{\partial x_i}}\!\!-\!\!\frac{\partial}{\partial y_i}\big( {\lambda_{ij}^{(0)}\frac{\partial \mathcal{H}_{\alpha_1}}{\partial x_{j}}} \big)\!\!-\!\!\frac{\partial}{\partial x_i}\big( {\lambda_{ij}^{(0)}\frac{\partial \mathcal{H}_{\alpha_1}}{\partial y_{j}}} \big)},\;\;\;\bm{y}\!\in\!\mathbf{Y},\\
&{\mathcal{Z}_{\alpha_1}}({\bm{y}},{T_{0}})\;\mathrm{is}\;1-\mathrm{periodic}\;\mathrm{in}\;{\bm{y}}{\rm{, }}\;\;\;{\int_{\mathbf{Y}}}\mathcal{Z}_{\alpha_1}d\mathbf{Y}=0.
\end{aligned} \right.
\end{equation}
\begin{equation}
\left\{
\begin{aligned}
&\frac{\partial}{\partial y_i}\big[ {\lambda_{ij}^{(0)}{\frac{\partial \mathcal{W}_{\alpha_1\alpha_2}}{\partial y_j}}} \big]\!\!=\!\!-\frac{\partial}{\partial y_i}\Big({{\mathcal{M}_{\alpha_1}}\mathbf{D}^{(0,1)}{\lambda_{i\alpha_2}^{(0)}}}\!+\! {{\mathcal{M}_{\alpha_1}}\mathbf{D}^{(0,1)}{\lambda_{ij}^{(0)}}\frac{\partial \mathcal{H}_{\alpha_2}}{\partial y_j}}\Big),\;\bm{y}\!\in\!\mathbf{Y},\\
&{\mathcal{W}_{\alpha_1\alpha_2}}({\bm{y}},{T_{0}})\;\mathrm{is}\;1-\mathrm{periodic}\;\mathrm{in}\;{\bm{y}}{\rm{, }}\;{\int_{\mathbf{Y}}}\mathcal{W}_{\alpha_1\alpha_2}d\mathbf{Y}=0.
\end{aligned} \right.
\end{equation}
\begin{equation}
\left\{
\begin{aligned}
&\frac{\partial}{\partial y_j}\big[ {c_{ijkl}^{(0)}{\frac{\partial {\mathcal{N}_{km}^{\alpha_1\alpha_2}}}{\partial y_l}}} \big]= { {\widehat c}_{i\alpha_1m\alpha_2}- {c_{i\alpha_1m\alpha_2}^{(0)}}-\frac{\partial}{\partial y_j}\big( {c_{ijk\alpha_1}^{(0)}{\mathcal{N}_{km}^{\alpha_2}}} \big)}\\
&\quad\quad\quad\quad\quad\quad\quad\quad\;{-{c_{i\alpha_1kj}^{(0)}\frac{\partial \mathcal{N}_{km}^{\alpha_2}}{\partial y_j}}},\;\;\;\bm{y}\in \mathbf{Y},\\
&{\mathcal{N}_{km}^{\alpha_1\alpha_2}}({\bm{y}},{T_{0}})\;\mathrm{is}\;1-\mathrm{periodic}\;\mathrm{in}\;{\bm{y}}{\rm{, }}\;\;\;{\int_{\mathbf{Y}}}\mathcal{N}_{km}^{\alpha_1\alpha_2}d\mathbf{Y}=0.
\end{aligned} \right.
\end{equation}
\begin{equation}
\left\{
\begin{aligned}
&\frac{\partial}{\partial y_j}\big[ {c_{ijkl}^{(0)}{\frac{\partial {\mathcal{F}_{k}^{\alpha_1}}}{\partial y_l}}} \big]= \delta_{i\alpha_1}\big(\rho^{(0)}-\widehat \rho\big),\;\;\;\bm{y}\in \mathbf{Y},\\
&{\mathcal{F}_{k}^{\alpha_1}}({\bm{y}},{T_{0}})\;\mathrm{is}\;1-\mathrm{periodic}\;\mathrm{in}\;{\bm{y}}{\rm{, }}\;\;\;{\int_{\mathbf{Y}}}\mathcal{F}_{k}^{\alpha_1}d\mathbf{Y}=0.
\end{aligned} \right.
\end{equation}
\begin{equation}
\left\{
\begin{aligned}
&\frac{\partial}{\partial y_j}\big[ {c_{ijkl}^{(0)}{\frac{\partial \mathcal{X}_{k}^{\alpha_1}}{\partial y_l}}} \big] =  {\beta _{i\alpha_1}^{(0)} - {{\widehat \beta }_{i\alpha_1}} -c_{i\alpha_1kl}^{(0)}{\frac{\partial \mathcal{P}_k}{\partial y_l}}}{-\frac{\partial}{\partial y_j}( {c_{ijk\alpha_1}^{(0)}{\mathcal{P}_k}})}\\
&\quad\quad\quad\quad\quad\quad\quad\;{+{\frac{\partial}{\partial y_j}}( {\beta _{ij}^{(0)}{\mathcal{M}_{\alpha_1}}} )},\;\;\;\bm{y}\in \mathbf{Y},\\
&\mathcal{X}_{k}^{\alpha_1}({\bm{y}},{T_{0}})\;\mathrm{is}\;1-\mathrm{periodic}\;\mathrm{in}\;{\bm{y}}{\rm{, }}\;\;\;{\int_{\mathbf{Y}}}\mathcal{X}_{k}^{\alpha_1}d\mathbf{Y}=0.
\end{aligned} \right.
\end{equation}
\begin{equation}
\left\{
\begin{aligned}
&\frac{\partial}{\partial y_j}\big[ {c_{ijkl}^{(0)}{\frac{\partial \mathcal{A}_{km}^{\alpha_1}}{\partial y_l}}} \big] ={ \frac{\partial{\widehat c}_{ijm\alpha_1}}{\partial x_j} -  {\frac{\partial{c}_{ijm\alpha_1}^{(0)}}{\partial x_j}}- \frac{\partial}{\partial y_j}\big( {c_{ijkl}^{(0)}\frac{\partial \mathcal{N}^{\alpha_1}_{km}}{\partial x_{l}}} \big)}\\
&\quad\quad\quad\quad\quad\quad\quad{-\frac{\partial}{\partial x_j}\big({c_{ijkl}^{(0)}\frac{\partial \mathcal{N}^{\alpha_1}_{km}}{\partial y_{l}}} \big)},\;\;\;\bm{y}\in \mathbf{Y},\\
&\mathcal{A}_{km}^{\alpha_1}({\bm{y}},{T_{0}})\;\mathrm{is}\;1-\mathrm{periodic}\;\mathrm{in}\;{\bm{y}}{\rm{, }}\;\;\;{\int_{\mathbf{Y}}}\mathcal{A}_{km}^{\alpha_1}d\mathbf{Y}=0.
\end{aligned} \right.
\end{equation}
\begin{equation}
\left\{
\begin{aligned}
&\frac{\partial}{\partial y_j}\big[ {c_{ijkl}^{(0)}{\frac{\partial \mathcal{B}_{k}}{\partial y_l}}} \big] \!=\!{\frac{\partial{\beta}_{ij}^{(0)}}{\partial x_j}}\!-\!{ \frac{\partial{\widehat \beta}_{ij}}{\partial x_j}\!-\!\frac{\partial}{\partial y_j}\big( {c_{ijkl}^{(0)}\frac{\partial \mathcal{P}_{k}}{\partial x_{l}}} \big)\!-\!\frac{\partial}{\partial x_j}\big( {c_{ijkl}^{(0)}\frac{\partial \mathcal{P}_{k}}{\partial y_{l}}} \big)},\;\;\;\bm{y}\in \mathbf{Y},\\
&\mathcal{B}_{k}({\bm{y}},{T_{0}})\;\mathrm{is}\;1-\mathrm{periodic}\;\mathrm{in}\;{\bm{y}}{\rm{, }}\;\;\;{\int_{\mathbf{Y}}}\mathcal{B}_{k}d\mathbf{Y}=0.
\end{aligned} \right.
\end{equation}
\begin{equation}
\left\{
\begin{aligned}
&\frac{\partial}{\partial y_j}\big[ {c_{ijkl}^{(0)}{\frac{\partial \mathcal{C}_{k}^{\alpha_1}}{\partial y_l}}} \big]  = \frac{\partial}{\partial y_j}\Big( {{\mathcal{M}_{\alpha_1}}\mathbf{D}^{(0,1)}\beta _{i j}^{(0)}- {\mathcal{M}_{\alpha_1}}\mathbf{D}^{(0,1)}c_{i jkl}^{(0)}\frac{\partial \mathcal{P}_k}{\partial y_l}} \Big),\;\;\;\bm{y}\in \mathbf{Y},\\
&\mathcal{C}_{k}^{\alpha_1}({\bm{y}},{T_{0}})\;\mathrm{is}\;1-\mathrm{periodic}\;\mathrm{in}\;{\bm{y}}{\rm{, }}\;\;\;{\int_{\mathbf{Y}}}\mathcal{C}_{k}^{\alpha_1}d\mathbf{Y}=0.
\end{aligned} \right.
\end{equation}
\begin{equation}
\left\{
\begin{aligned}
&\!\!\frac{\partial}{\partial y_j}\big[\!{c_{ijkl}^{(0)}{\frac{\partial \mathcal{D}_{km}^{\alpha_1\alpha_2}}{\partial y_l}}}\!\big]\!\!=\!\!-\!\frac{\partial}{\partial y_j}\Big(\! {{\mathcal{M}_{\alpha_1}}\!\mathbf{D}^{(0,1)}c_{ijm\alpha_2}^{(0)}\!+\! {\mathcal{M}_{\alpha_1}}\!\mathbf{D}^{(0,1)}c_{ijkl}^{(0)}\frac{\partial \mathcal{N}^{\alpha_2}_{km}}{\partial y_l}}\!\Big),\;\bm{y}\!\in \!\mathbf{Y},\\
&\!\!\mathcal{D}_{km}^{\alpha_1\alpha_2}({\bm{y}},{T_{0}})\;\mathrm{is}\;1-\mathrm{periodic}\;\mathrm{in}\;{\bm{y}}{\rm{, }}\;{\int_{\mathbf{Y}}}\mathcal{D}_{km}^{\alpha_1\alpha_2}d\mathbf{Y}=0.
\end{aligned} \right.
\end{equation}
\begin{rmk}
By virtue of the hypotheses (A$_1$)-(A$_2$) and Lax-Milgram theorem, the existence and uniqueness of solutions for auxiliary cell problems (2.20)-(2.23) and (2.26)-(2.41) are established for any fixed macroscopic temperature $T_0$.
\end{rmk}
\begin{rmk}
Based on the definitions (2.20)-(2.23) and (2.26)-(2.41) of cell functions, it can be proved that all cell functions are continuous with respect to macroscopic temperature $T_0$, with details provided in Appendix B.
\end{rmk}
\section{The error analyses of multi-scale asymptotic solutions}
\subsection{The proof of local balance preserving by local error analysis}
The fundamental motivation for developing high-order multi-scale solutions lies in their capability in local balance preserving, which shall be demonstrated by the local pointwise error analysis. First, the residual functions corresponding to LOMS and HOMS approximate solutions are defined as follows.
\begin{equation}
\left\{
\begin{aligned}
&{\mathring{T}^{[1\varepsilon]} }(\bm{x},t)={T^{\varepsilon} }(\bm{x},t)-{T^{[1\varepsilon]} }(\bm{x},t),\;{\mathring{T}^{[2\varepsilon]} }(\bm{x},t)={T^{\varepsilon} }(\bm{x},t)-{ T^{[2\varepsilon]}}(\bm{x},t),\\
&{\mathring{\mathit{\Phi}}^{[1\varepsilon]} }(\bm{x},t)={\mathit{\Phi}^{\varepsilon} }(\bm{x},t)-{\mathit{\Phi}^{[2\varepsilon]} }(\bm{x},t),\;{\mathring{\mathit{\Phi}}^{[2\varepsilon]} }(\bm{x},t)={\mathit{\Phi}^{\varepsilon} }(\bm{x},t)-{\mathit{\Phi}^{[2\varepsilon]}}(\bm{x},t),\\
&{\mathring{U}_{i}^{[1\varepsilon]}}(\bm{x},t)={U_i^{\varepsilon} }(\bm{x},t)-{U_i^{[1\varepsilon]} }(\bm{x},t),\;{\mathring{U}_{i}^{[2\varepsilon]}}(\bm{x},t)={U_i^{\varepsilon} }(\bm{x},t)-{ U_i^{[2\varepsilon]}}(\bm{x},t).
\end{aligned}\right.
\end{equation}

Then substituting the aforementioned residual functions ${\mathring{T}^{[1\varepsilon]}}$, ${\mathring{\mathit{\Phi}}^{[1\varepsilon]}}$ and ${\mathring{U}_{i}^{[1\varepsilon]}}$ into (1.1), the residual equations for the LOMS solutions can be derived as below.
\begin{equation}
\left\{
\begin{aligned}
&\rho^{\varepsilon}c^{\varepsilon}\frac{{\partial {{\mathring{T}^{[1\varepsilon]} }}}}{{\partial t}}- \frac{\partial }{{\partial {x_i}}}\Big( {{k_{ij}^{\varepsilon}}\frac{{\partial {{\mathring{T}^{[1\varepsilon]} } }}}{{\partial {x_j}}}}\Big)-{{\lambda_{ij}^{\varepsilon}}\frac{\partial {\mathring{\mathit{\Phi}}^{[1\varepsilon]}}}{\partial {x_i}}}\frac{\partial {\mathring{\mathit{\Phi}}^{[1\varepsilon]}}}{\partial {x_j}}+T^{\varepsilon}\beta_{ij}^{\varepsilon}\frac{\partial}{\partial t}\Big(\frac{\partial \mathring{U}_{i}^{[1\varepsilon]}}{\partial x_j}\Big)\\
&\quad\quad\quad\quad\quad\quad\quad\quad\;=H_0(\bm{x},\bm{y},t)+\varepsilon H_1(\bm{x},\bm{y},t),\;\text{in}\;\Omega\times(0,\mathcal T),\\
&-\frac{\partial }{{\partial {x_i}}}\Big( {{\lambda_{ij}^{\varepsilon}}\frac{{\partial {\mathring{\mathit{\Phi}}^{[1\varepsilon]}}}}{{\partial {x_j}}}} \Big) = E_0(\bm{x},\bm{y},t)+\varepsilon E_1(\bm{x},\bm{y},t),\;\text{in}\;\Omega\times(0,\mathcal T),\\
&\rho^{\varepsilon} \frac{{\partial^2 \mathring{U}_{i}^{[1\varepsilon]}}}{{\partial t^2}}-\frac{\partial }{{\partial {x_j}}}\Big( {{c_{ijkl}^{\varepsilon}}\frac{{\partial \mathring{U}_{k}^{[1\varepsilon]}}}{{\partial {x_l}}}}-\beta_{ij}^{\varepsilon}\mathring{T}^{[1\varepsilon]}\Big)\\
&\quad\quad\quad\quad\quad\quad\quad\quad\;=F_{i0}(\bm{x},\bm{y},t)+\varepsilon F_{i1}(\bm{x},\bm{y},t),\;\text{in}\;\Omega\times(0,\mathcal T),
\end{aligned} \right.
\end{equation}
where the detailed expressions of functions ${H}_{0}$, ${H}_{1}$, ${E}_{0}$, ${E}_{1}$, $F_{i0}$ and $F_{i1}$ are uncomplicated to derive and not be exhibited in this study.

Similarly, putting the residual functions ${\mathring{T}^{[2\varepsilon]}}$, ${\mathring{\mathit{\Phi}}^{[2\varepsilon]}}$ and ${\mathring{U}_{i}^{[2\varepsilon]}}$ into (1.1) leads to the following residual equations for the HOMS solutions.
\begin{equation}
\left\{
\begin{aligned}
&\rho^{\varepsilon}c^{\varepsilon}\frac{{\partial {{\mathring{T}^{[2\varepsilon]} }}}}{{\partial t}}- \frac{\partial }{{\partial {x_i}}}\Big( {{k_{ij}^{\varepsilon}}\frac{{\partial {{\mathring{T}^{[2\varepsilon]} } }}}{{\partial {x_j}}}}\Big)-{{\lambda_{ij}^{\varepsilon}}\frac{\partial {\mathring{\mathit{\Phi}}^{[2\varepsilon]}}}{\partial {x_i}}}\frac{\partial {\mathring{\mathit{\Phi}}^{[2\varepsilon]}}}{\partial {x_j}}+T^{\varepsilon}\beta_{ij}^{\varepsilon}\frac{\partial}{\partial t}\Big(\frac{\partial \mathring{U}_{i}^{[2\varepsilon]}}{\partial x_j}\Big)\\
&\quad\quad\quad\quad\;\;=\varepsilon H_2(\bm{x},\bm{y},t),\;\text{in}\;\Omega\times(0,\mathcal T),\\
&-\frac{\partial }{{\partial {x_i}}}\Big( {{\lambda_{ij}^{\varepsilon}}\frac{{\partial {\mathring{\mathit{\Phi}}^{[2\varepsilon]}}}}{{\partial {x_j}}}} \Big) = \varepsilon E_2(\bm{x},\bm{y},t),\;\text{in}\;\Omega\times(0,\mathcal T),\\
&\rho^{\varepsilon} \frac{{\partial^2 \mathring{U}_{i}^{[2\varepsilon]}}}{{\partial t^2}}-\frac{\partial }{{\partial {x_j}}}\Big( {{c_{ijkl}^{\varepsilon}}\frac{{\partial \mathring{U}_{k}^{[2\varepsilon]}}}{{\partial {x_l}}}}-\beta_{ij}^{\varepsilon}\mathring{T}^{[2\varepsilon]}\Big)=\varepsilon F_{i2}(\bm{x},\bm{y},t),\;\text{in}\;\Omega\times(0,\mathcal T),
\end{aligned} \right.
\end{equation}
where the detailed expressions of functions ${H}_{2}$, ${E}_{2}$ and $F_{i2}$ are also trivial to derive and not be displayed in this study.

Observing the above-mentioned pointwise error results, it reveals that LOMS solutions fail to preserve local balance of heat quantity, electric charge and stress; specifically, the $\varepsilon$-independent terms ${H}_{0}$, ${E}_{0}$ and $F_{i0}$ in (3.2) do not vanish as the microstructural parameter $\varepsilon$ tends to zero. In contrast, augmented by the high-order correction terms, the HOMS solutions shall yield $O(\varepsilon)$-order pointwise errors, thereby preserving the local heat quantity balance of thermal equation, local electric charge balance of electronic equation, and local stress balance of mechanical equations inherent in the original governing equations (1.1). This capability serves as the primary motivation for developing the HOMS solutions that manifest high-fidelity computational performance for the composite structures with highly spatial heterogeneities.
\subsection{The proof of convergence by global error estimate}
To ensure the computational reliability of global multi-physics simulations, the global error estimate in the integral sense is obtained for the high-order multi-scale solutions. Firstly, the following additional assumptions are stated for sequel proof.
\begin{enumerate}
\item[(B$_1$)]
Assume that $\Omega$ is a union of integral PUCs, namely $\bar{\Omega}=\cup_{\mathbf{z}\in I_{\varepsilon}}\varepsilon(\mathbf{z}+\bar{\mathbf{Y}})$, where the index set $I_{\varepsilon}=\{\mathbf{z}=(z_1,\cdots,z_n)\in Z^n,\varepsilon(\mathbf{z}+\bar{\mathbf{Y}})\subset \bar{\Omega}\}$. In addition, let $E_\mathbf{z}=\varepsilon(\mathbf{z}+\mathbf{Y})$ be the translational unit cell and $\partial E_\mathbf{z}$ be its boundary.
\item[(B$_2$)]
Assume that $\displaystyle\frac{\partial \rho^{\varepsilon}({\bm{x}},T^{\varepsilon})}{\partial t}$, $\displaystyle\frac{\partial c^{\varepsilon}({\bm{x}},T^{\varepsilon})}{\partial t}$ and $\displaystyle\frac{\partial c_{ijkl}^{\varepsilon}({\bm{x}},T^{\varepsilon})}{\partial t}\in L^\infty(\Omega\times(0,\mathcal T))$.
\item[(B$_3$)]
The nonlinear coupling term $T^{\varepsilon}(\bm{x},t)\beta_{ij}^{\varepsilon}({\bm{x}},{T^{\varepsilon}})$ of thermal equation in multi-scale multi-physics equations (1.1) is simplified as the linear term $\widetilde T\beta_{ij}^{\varepsilon}({\bm{x}},{T^{\varepsilon}})$.
\item[(B$_4$)]
Let $k_{ij}(\bm{y},{T^\varepsilon})=k_{ii}(\bm{y},{T^\varepsilon})\delta_{ij}$, $\lambda_{ij}(\bm{y},{T^\varepsilon})=\lambda_{ii}(\bm{y},{T^\varepsilon})\delta_{ij}$ and $\beta_{ij}(\bm{y},{T^\varepsilon})=\beta_{ii}(\bm{y},{T^\varepsilon})\delta_{ij}$, where $\delta_{ij}$ is a Kronecker symbol. Besides, let $\Delta_1,\cdots,\Delta_n$ be the middle hyperplanes of PUC $\mathbf{Y}=[0,1]^n$. Then suppose that material parameters $\rho(\bm{y},{T^\varepsilon})$, $c(\bm{y},{T^\varepsilon})$, $k_{ii}(\bm{y},{T^\varepsilon})$, $\lambda_{ii}(\bm{y},{T^\varepsilon})$, $\beta_{ii}(\bm{y},{T^\varepsilon})$ and $c_{ijkl}(\bm{y},{T^\varepsilon})$ are symmetric with respect to $\Delta_1,\cdots,\Delta_n$ for stationary ${T^\varepsilon }\in [T_{min},T_{max}+C_*]$.
\item[(B$_5$)]
Apply the homogeneous Dirichlet boundary condition to replace the periodic boundary condition for all auxiliary cell functions \cite{R5,R24,R37}.
\end{enumerate}

Based on the aforementioned assumptions, we further derive the following initial-boundary conditions for the residual equations (3.3) of the HOMS solutions, which will be employed for global error estimation
\begin{equation}
\left\{
\begin{aligned}
&\mathring{T}^{[2\varepsilon]}(\bm{x},t)=-\varepsilon{\mathcal{M}_{\alpha_1}}\frac{\partial T_{0}}{\partial x_{\alpha_1}}-\varepsilon^2\big[\mathcal{Q}\frac{{\partial {T_{0}}}}{{\partial t}}+{\mathcal{M}_{\alpha_1\alpha_2}}\frac{\partial^2 T_{0}}{\partial x_{\alpha_1}\partial x_{\alpha_2}}+\mathcal{R}_{\alpha_1}\frac{\partial T_{0}}{\partial x_{\alpha_1}}\\
&\!+\!{\mathcal{O}_{\alpha_1\alpha_2}}\frac{\partial T_{0}}{\partial x_{\alpha_1}}\!\frac{\partial T_{0}}{\partial x_{\alpha_2}}\!+\!{\mathcal{G}_{\alpha_1\alpha_2}}\frac{\partial \mathit{\Phi}_{0}}{\partial x_{\alpha_1}}\!\frac{\partial \mathit{\Phi}_{0}}{\partial x_{\alpha_2}}\!+\!{\mathcal{J}_{\alpha_1\alpha_2}}\frac{\partial^2 U_{\alpha_10}}{\partial x_{\alpha_2}\partial t}\big]=0,\;\text{on}\;\partial {\Omega}\times(0,\mathcal T),\\
&\mathring{\mathit{\Phi}}^{[2\varepsilon]}(\bm{x},t)=-\varepsilon \mathcal{H}_{\alpha_1}\frac{\partial \mathit{\Phi}_{0}}{\partial x_{\alpha_1}}- \varepsilon^2\big[\mathcal{H}_{\alpha_1\alpha_2}\frac{\partial^2 \mathit{\Phi}_{0}}{\partial x_{\alpha_1}\partial x_{\alpha_2}}+\mathcal{Z}_{\alpha_1}\frac{\partial \mathit{\Phi}_{0}}{\partial x_{\alpha_1}}\\
&+\mathcal{W}_{\alpha_1\alpha_2}\frac{\partial T_{0}}{\partial x_{\alpha_1}}\frac{\partial \mathit{\Phi}_{0}}{\partial x_{\alpha_2}}\big]=0,\;\text{on}\;\partial\Omega\times(0,\mathcal T),\\
&\mathring{U}_{i}^{[2\varepsilon]}(\bm{x},t)=-\varepsilon\Big[\mathcal{N}_{im}^{\alpha_1}\frac{\partial U_{m0}}{\partial x_{\alpha_1}}+\mathcal{P}_{i}(T_0-\widetilde T)\Big]-\varepsilon^2\big[\mathcal{N}_{im}^{\alpha_1\alpha_2}\frac{\partial^2 U_{m0}}{\partial x_{\alpha_1}\partial x_{\alpha_2}}\\
&+\mathcal{F}_{i}^{\alpha_1}\frac{{\partial^2 {U_{\alpha_10}}}}{{\partial t^2}}+\mathcal{X}_i^{\alpha_1}\frac{\partial T_{0}}{\partial x_{\alpha_1}}+\mathcal{A}_{im}^{\alpha_1}\frac{\partial U_{m0}}{\partial x_{\alpha_1}}+\mathcal{B}_{i}(T_0-\widetilde T)\\
&+\mathcal{C}_{i}^{\alpha_1}\frac{\partial T_{0}}{\partial x_{\alpha_1}}(T_0-\widetilde T)+{\mathcal{D}_{im}^{\alpha_1\alpha_2}}\frac{\partial T_{0}}{\partial x_{\alpha_1}}\frac{\partial U_{m0}}{\partial x_{\alpha_2}}\big]=0,\;\text{on}\;\partial\Omega\times(0,\mathcal T),\\
&\mathring{T}^{[2\varepsilon]}({\bm{x}},0)=:\varepsilon \widetilde\chi_2(\bm{x}),\;\mathring{U}_{i}^{[2\varepsilon]}({\bm{x}},0)=:\varepsilon \widetilde\eta_{2i}(\bm{x}),\;\frac{\partial \mathring{U}_{i}^{[2\varepsilon]}({\bm{x}},0)}{\partial t}=:\varepsilon \widetilde\zeta_{2i}(\bm{x}),\;\text{in}\;\Omega,
\end{aligned} \right.
\end{equation}

Furthermore, based on the above assumptions, the following lemma is obtained, which will be utilized to implement the global error estimate of the HOMS solutions.
\begin{lemma}
Define the differential operators $\displaystyle\sigma_{T\mathbf{Y}}(\chi)=n_i k_{ij}(\bm{y},{T^\varepsilon })\frac{\partial\chi}{\partial y_j}$ for temperature field, $\displaystyle\sigma_{\mathit{\Phi}\mathbf{Y}}(\eta)=n_i k_{ij}(\bm{y},{T^\varepsilon})\frac{\partial\eta}{\partial y_j}$ for electric potential field and $\displaystyle\sigma_{i\mathbf{Y}}(\bm{\phi})=n_j c_{ijkl}(\bm{y},{T^\varepsilon })\frac{\partial\phi_k}{\partial y_l}$ for displacement field, respectively. Under the assumptions (A$_1$)-(A$_2$) and (B$_4$)-(B$_5$), the normal derivatives of all microscopic cell functions are proved to be continuous on the boundary of PUC $\mathbf{Y}$, following the identical approach in references \cite{R5,R24,R37}. Furthermore, $\sigma_{T\mathbf{Y}}(T^{[2\varepsilon]})$, $\sigma_{\mathit{\Phi}\mathbf{Y}}(\mathit{\Phi}^{[2\varepsilon]})$ and $\sigma_{i\mathbf{Y}}(\bm{U}^{[2\varepsilon]})$ are proved to be continuous on the boundary of PUC $\mathbf{Y}$.
\end{lemma}

Next, the global error estimate with explicit convergence order for the HOMS solutions of the time-varying multi-scale nonlinear multi-physics equations (1.1) is established and presented in the following theorem.
\begin{theorem}
Let $T^{\varepsilon}$, $\mathit{\Phi}^{{\varepsilon}}$ and ${\bm{U}}^{{\varepsilon}}$ be the weak solutions of multi-scale nonlinear multi-physics equations (1.1), while $T_0$, $\mathit{\Phi}_0$ and ${\bm{U}}_0$ denote the weak solutions of corresponding homogenized equations (2.24). Let $T^{[2\varepsilon]}$, $\mathit{\Phi}^{{[2\varepsilon]}}$ and ${\bm{U}}^{{[2\varepsilon]}}$ be the HOMS solutions given by formulas (2.17)-(2.19). Under the above hypotheses (A$_1$)-(A$_3$) and (B$_1$)-(B$_5$), the following global error estimates are derived.
\begin{equation}
{\begin{aligned}
&\big\|{\mathit{\Phi}^{\varepsilon}}-{\mathit{\Phi}^{[2\varepsilon]}}\big\|_{L^\infty(0,\mathcal T;H_0^1(\Omega))}+\big\|{T^{\varepsilon}}-{T^{[2\varepsilon]}}\big\|_{L^\infty(0,\mathcal T;L^2(\Omega))}\\
&+
\big\|{T^{\varepsilon}}-{T^{[2\varepsilon]}}\big\|_{L^2(0,\mathcal T;H^1_0(\Omega))}+\Big\|\frac{\partial{\bm{U}}^{{\varepsilon}}}{\partial t}-\frac{\partial{\bm{U}}^{{[2\varepsilon]}}}{\partial t}\Big\|_{L^\infty(0,\mathcal T;(L^2(\Omega))^n)}\\
&+\big\|{\bm{U}}^{{\varepsilon}}-{\bm{U}}^{{[2\varepsilon]}}\big\|_{L^\infty(0,\mathcal T;(H^1_0(\Omega))^n)}\leq C(\Omega,\mathcal T)\varepsilon.
\end{aligned}}
\end{equation}
\end{theorem}
Here, $C(\Omega,\mathcal T)$ is a positive constant depending on $\Omega$ and $\mathcal T$, but independent of $\varepsilon$.\\
$\mathbf{Proof:}$
The error estimate (3.5) is established on the residual equations (3.3) and associated initial-boundary conditions (3.4). Firstly, multiplying the electronic residual equation in (3.3) by $\mathit{\Phi}^{{[2\varepsilon]}}(\bm{x},t)$ and integrating on $\Omega$, the following equality is derived upon further applying Green's formula and substituting the corresponding boundary condition in (3.4).
\begin{equation}
\begin{aligned}
&\int_{\Omega}{\lambda_{ij}^{\varepsilon}({\bm{x}},{T^\varepsilon })\frac{{\partial {\mathring{\mathit{\Phi}}^{[2\varepsilon]}}}}{\partial {x_j}}}\frac{{\partial {\mathring{\mathit{\Phi}}^{[2\varepsilon]}}}}{\partial {x_i}}d\Omega\\
&=\int_{\Omega}\varepsilon E_2(\bm{x},\bm{y},t){\mathring{\mathit{\Phi}}^{[2\varepsilon]}}d\Omega+\sum\limits _{\mathbf{z}\in I_\xi}\int_{\partial E_\mathbf{z}}\sigma_{\mathit{\Phi}\mathbf{Y}}({\mathring{\mathit{\Phi}}^{[2\varepsilon]}}){\mathring{\mathit{\Phi}}^{[2\varepsilon]}}d\Gamma_{\bm{y}}.
\end{aligned}
\end{equation}
Invoking lemma 3.1, we shall derive the following result for the boundary integral terms on $\partial E_\mathbf{z}$ in (3.6).
\begin{equation}
\begin{aligned}
&\sum\limits _{\mathbf{z}\in I_\xi}\int_{\partial E_\mathbf{z}}\sigma_{\mathit{\Phi}\mathbf{Y}}({\mathring{\mathit{\Phi}}^{[2\varepsilon]}}){\mathring{\mathit{\Phi}}^{[2\varepsilon]}}d\Gamma_{\bm{y}}=\sum\limits _{\mathbf{z}\in I_\xi}\int_{\partial E_\mathbf{z}}\sigma_{\mathit{\Phi}\mathbf{Y}}({{\mathit{\Phi}}^{\varepsilon}}-{{\mathit{\Phi}}^{[2\varepsilon]}}){\mathring{\mathit{\Phi}}^{[2\varepsilon]}}d\Gamma_{\bm{y}}\\
&=-\sum\limits _{\mathbf{z}\in I_\xi}\int_{\partial E_\mathbf{z}}\sigma_{\mathit{\Phi}\mathbf{Y}}({{\mathit{\Phi}}^{[2\varepsilon]}}){\mathring{\mathit{\Phi}}^{[2\varepsilon]}}d\Gamma_{\bm{y}}=0.
\end{aligned}
\end{equation}
Furthermore, employing assumption (A$_1$) and Poincar$\mathrm{\acute{e}}$-Friedrichs inequality to the left side of (3.6), it follows that
\begin{equation}
\Big|\int_{\Omega}{\lambda_{ij}^{\varepsilon}({\bm{x}},{T^\varepsilon })\frac{{\partial {\mathring{\mathit{\Phi}}^{[2\varepsilon]}}}}{\partial {x_j}}}\frac{{\partial {\mathring{\mathit{\Phi}}^{[2\varepsilon]}}}}{\partial {x_i}}d\Omega\Big| \geq C\left \|{\mathring{\mathit{\Phi}}^{[2\varepsilon]}}\right\|_{H_0^1(\Omega)}^2.
\end{equation}
After that, exploiting the Schwarz's inequality and equality (3.7), we derive the following inequality by transforming the right side of (3.6)
\begin{equation}
\begin{aligned}
&\Big|\int_{\Omega}\varepsilon E_2(\bm{x},\bm{y},t){\mathring{\mathit{\Phi}}^{[2\varepsilon]}}d\Omega+\sum\limits _{\mathbf{z}\in I_\xi}\int_{\partial E_\mathbf{z}}\sigma_{\mathit{\Phi}\mathbf{Y}}({\mathring{\mathit{\Phi}}^{[2\varepsilon]}}){\mathring{\mathit{\Phi}}^{[2\varepsilon]}}d\Gamma_{\bm{y}}\Big|\\
&\leq \left\|\varepsilon E_2(\bm{x},\bm{y},t)\right\|_{L^2(\Omega)} \left\|\mathring{\mathit{\Phi}}^{[2\varepsilon]}\right\|_{L^2(\Omega)}\\
&\leq C\varepsilon\left\| E_2(\bm{x},\bm{y},t)\right\|_{L^2(\Omega)} \left\|\mathring{\mathit{\Phi}}^{[2\varepsilon]}\right\|_{H_0^1(\Omega)}.
\end{aligned}
\end{equation}
A combination of (3.8) and (3.9) leads to the following inequality
\begin{equation}
\big\|{\mathit{\Phi}^{\varepsilon}}-{\mathit{\Phi}^{[2\varepsilon]}}\big\|_{H_0^1(\Omega)}=\left\|\mathring{\mathit{\Phi}}^{[2\varepsilon]}\right\|_{H_0^1(\Omega)}\leq C\varepsilon\left\| E_2(\bm{x},\bm{y},t)\right\|_{L^2(\Omega)}.
\end{equation}
Then, leveraging the arbitrariness of temporal variable $t$ in (3.10), we readily arrive at the explicit convergence estimate $\big\|{\mathit{\Phi}^{\varepsilon}}-{\mathit{\Phi}^{[2\varepsilon]}}\big\|_{L^\infty(0,\mathcal T;H_0^1(\Omega))}\leq C\varepsilon$ for the electric potential field.

In what follows, we proceed to derive the error estimates for the temperature and displacement fields. Since $\displaystyle\frac{\partial \mathring{U}_{i}^{[2\varepsilon]}}{\partial t}\in L^\infty(0,\mathcal T;L^2(\Omega))$, it cannot be directly utilized as a test function for hyperbolic dynamics part of the residual equations (3.3). To tackle this challenge, the density argument approach in \cite{R14} is employed. For the sake of brevity, the details of this argument are omitted. Next, multiplying the thermal residual equation in (3.3) by $\mathring{T}^{[2\varepsilon]}$ and the mechanical equations in (3.3) by $\displaystyle\frac{\partial \mathring{U}_{i}^{[2\varepsilon]}}{\partial t}$ respectively, and integrating over $\Omega$, the following equalities are derived by further applying Green's formula and substituting associated boundary conditions in (3.4).
\begin{equation}
\left\{
\begin{aligned}
&\int_{\Omega}\rho^{\varepsilon}c^{\varepsilon}\frac{{\partial {{\mathring{T}^{[2\varepsilon]}}}}}{{\partial t}}\mathring{T}^{[2\varepsilon]}d\Omega+\int_{\Omega}{{k_{ij}^{\varepsilon}}\frac{{\partial {{\mathring{T}^{[2\varepsilon]}}}}}{{\partial {x_j}}}}\frac{\partial \mathring{T}^{[2\varepsilon]}}{{\partial {x_i}}}d\Omega\\
&-\int_{\Omega}{{\lambda_{ij}^{\varepsilon}}\frac{\partial {\mathring{\mathit{\Phi}}^{[2\varepsilon]}}}{\partial {x_i}}}\frac{\partial {\mathring{\mathit{\Phi}}^{[2\varepsilon]}}}{\partial {x_j}}\mathring{T}^{[2\varepsilon]}d\Omega+\int_{\Omega}\widetilde T\beta_{ij}^{\varepsilon}\frac{\partial}{\partial t}\Big(\frac{\partial \mathring{U}_{i}^{[2\varepsilon]}}{\partial x_j}\Big)\mathring{T}^{[2\varepsilon]}d\Omega\\
&=\int_{\Omega}\varepsilon H_2(\bm{x},\bm{y},t)\mathring{T}^{[2\varepsilon]}d\Omega+\sum\limits _{\mathbf{z}\in I_\xi}\int_{\partial E_\mathbf{z}}\sigma_{T\mathbf{Y}}({\mathring{T}^{[2\varepsilon]}}){\mathring{T}^{[2\varepsilon]}}d\Gamma_{\bm{y}},\\
&\int_{\Omega}\rho^{\varepsilon} \frac{{\partial^2 \mathring{U}_{i}^{[2\varepsilon]}}}{{\partial t^2}}\frac{\partial \mathring{U}_{i}^{[2\varepsilon]}}{\partial t}d\Omega+\int_{\Omega}\Big( {{c_{ijkl}^{\varepsilon}}\frac{{\partial \mathring{U}_{k}^{[2\varepsilon]}}}{{\partial {x_l}}}}-\beta_{ij}^{\varepsilon}\mathring{T}^{[2\varepsilon]}\Big)\frac{\partial }{{\partial {x_j}}}\Big(\frac{\partial \mathring{U}_{i}^{[2\varepsilon]}}{\partial t}\Big)d\Omega\\
&=\int_{\Omega}\varepsilon F_{i2}(\bm{x},\bm{y},t)\frac{\partial \mathring{U}_{i}^{[2\varepsilon]}}{\partial t}d\Omega+\sum\limits _{\mathbf{z}\in I_\xi}\int_{\partial E_\mathbf{z}}\sigma_{i\mathbf{Y}}({\mathring{\bm{U}}^{[2\varepsilon]}})\frac{\partial \mathring{U}_{i}^{[2\varepsilon]}}{\partial t}d\Gamma_{\bm{y}}.
\end{aligned}\right.
\end{equation}
Invoking lemma 3.1, the following results are established for the boundary integral terms on $\partial E_\mathbf{z}$ in (3.11).
\begin{equation}
\left\{
\begin{aligned}
&\sum\limits _{\mathbf{z}\in I_\xi}\int_{\partial E_\mathbf{z}}\sigma_{T\mathbf{Y}}({\mathring{T}^{[2\varepsilon]}}){\mathring{T}^{[2\varepsilon]}}d\Gamma_{\bm{y}}=\sum\limits _{\mathbf{z}\in I_\xi}\int_{\partial E_\mathbf{z}}\sigma_{T\mathbf{Y}}({{T}^{\varepsilon}}-{{T}^{[2\varepsilon]}}){\mathring{T}^{[2\varepsilon]}}d\Gamma_{\bm{y}}\\
&=-\sum\limits _{\mathbf{z}\in I_\xi}\int_{\partial E_\mathbf{z}}\sigma_{T\mathbf{Y}}({{T}^{[2\varepsilon]}}){\mathring{T}^{[2\varepsilon]}}d\Gamma_{\bm{y}}=0,\\
&\sum\limits _{\mathbf{z}\in I_\xi}\int_{\partial E_\mathbf{z}}\sigma_{i\mathbf{Y}}({\mathring{\bm{U}}^{[2\varepsilon]}})\frac{\partial \mathring{U}_{i}^{[2\varepsilon]}}{\partial t}d\Gamma_{\bm{y}}=\sum\limits _{\mathbf{z}\in I_\xi}\int_{\partial E_\mathbf{z}}\sigma_{i\mathbf{Y}}({{\bm{U}}^{\varepsilon}}-{{\bm{U}}^{[2\varepsilon]}})\frac{\partial \mathring{U}_{i}^{[2\varepsilon]}}{\partial t}d\Gamma_{\bm{y}}\\
&=-\sum\limits _{\mathbf{z}\in I_\xi}\int_{\partial E_\mathbf{z}}\sigma_{i\mathbf{Y}}({{\bm{U}}^{[2\varepsilon]}})\frac{\partial \mathring{U}_{i}^{[2\varepsilon]}}{\partial t}d\Gamma_{\bm{y}}=0.
\end{aligned}\right.
\end{equation}
After that, the following two equalities are derived separately for multi-scale thermal and mechanical equations with temperature-dependent material parameters in (3.11).
\begin{equation}
\begin{aligned}
&\frac{1}{2}\frac{\partial}{\partial t}\Big[\int_{\Omega}\rho^{\varepsilon}({\bm{x}},{T^\varepsilon })c^{\varepsilon}({\bm{x}},{T^\varepsilon })\big(\mathring{T}^{[2\varepsilon]}\big)^2d\Omega\Big]+\int_{\Omega}{{k_{ij}^{\varepsilon}}({\bm{x}},{T^\varepsilon })\frac{{\partial {{\mathring{T}^{[2\varepsilon]}}}}}{{\partial {x_j}}}}\frac{\partial \mathring{T}^{[2\varepsilon]}}{{\partial {x_i}}}d\Omega\\
&-\frac{1}{2}\int_{\Omega}\frac{\partial\rho^{\varepsilon}({\bm{x}},\!{T^\varepsilon })}{\partial t}c^{\varepsilon}({\bm{x}},\!{T^\varepsilon })\big(\mathring{T}^{[2\varepsilon]}\big)^2d\Omega-\frac{1}{2}\int_{\Omega}\rho^{\varepsilon}({\bm{x}},\!{T^\varepsilon })\frac{\partial c^{\varepsilon}({\bm{x}},\!{T^\varepsilon })}{\partial t}\big(\mathring{T}^{[2\varepsilon]}\big)^2d\Omega\\
&-\int_{\Omega}{{\lambda_{ij}^{\varepsilon}}({\bm{x}},{T^\varepsilon })\frac{\partial {\mathring{\mathit{\Phi}}^{[2\varepsilon]}}}{\partial {x_i}}}\frac{\partial {\mathring{\mathit{\Phi}}^{[2\varepsilon]}}}{\partial {x_j}}\mathring{T}^{[2\varepsilon]}d\Omega+\int_{\Omega}\widetilde T\beta_{ij}^{\varepsilon}({\bm{x}},{T^\varepsilon })\frac{\partial}{\partial t}\Big(\frac{\partial \mathring{U}_{i}^{[2\varepsilon]}}{\partial x_j}\Big)\mathring{T}^{[2\varepsilon]}d\Omega\\
&=\int_{\Omega}\varepsilon H_2(\bm{x},\bm{y},t)\mathring{T}^{[2\varepsilon]}d\Omega.
\end{aligned}
\end{equation}
\begin{equation}
\begin{aligned}
&\frac{1}{2}\frac{\partial}{\partial t}\Big[\int_{\Omega}\rho^{\varepsilon}({\bm{x}},{T^\varepsilon })\big(\frac{\partial \mathring{U}_{i}^{[2\varepsilon]}}{\partial t}\big)^2d\Omega+\int_{\Omega} {{c_{ijkl}^{\varepsilon}}({\bm{x}},{T^\varepsilon })\frac{{\partial \mathring{U}_{k}^{[2\varepsilon]}}}{{\partial {x_l}}}}\frac{\partial \mathring{U}_{i}^{[2\varepsilon]}}{{\partial {x_j}}}d\Omega\Big]\\
&-\frac{1}{2}\int_{\Omega}\frac{\partial\rho^{\varepsilon}({\bm{x}},{T^\varepsilon })}{\partial t}\big(\frac{\partial \mathring{U}_{i}^{[2\varepsilon]}}{\partial t}\big)^2d\Omega-\frac{1}{2}\int_{\Omega}{\frac{\partial c_{ijkl}^{\varepsilon}({\bm{x}},{T^\varepsilon })}{\partial t}\frac{{\partial \mathring{U}_{k}^{[2\varepsilon]}}}{{\partial {x_l}}}}\frac{\partial \mathring{U}_{i}^{[2\varepsilon]}}{{\partial {x_j}}}d\Omega\\
&-\int_{\Omega}\beta_{ij}^{\varepsilon}({\bm{x}},{T^\varepsilon })\mathring{T}^{[2\varepsilon]}\frac{\partial }{{\partial {x_j}}}\Big(\frac{\partial \mathring{U}_{i}^{[2\varepsilon]}}{\partial t}\Big)d\Omega=\int_{\Omega}\varepsilon F_{i2}(\bm{x},\bm{y},t)\frac{\partial \mathring{U}_{i}^{[2\varepsilon]}}{\partial t}d\Omega.
\end{aligned}
\end{equation}
Then, implementing the calculation (3.13)+$\widetilde T\times$(3.14) yields the following equality.
\begin{equation}
\begin{aligned}
&\frac{1}{2}\frac{\partial}{\partial t}\Big[\int_{\Omega}\rho^{\varepsilon}({\bm{x}},{T^\varepsilon })c^{\varepsilon}({\bm{x}},{T^\varepsilon })\big(\mathring{T}^{[2\varepsilon]}\big)^2d\Omega\Big]+\int_{\Omega}{{k_{ij}^{\varepsilon}}({\bm{x}},{T^\varepsilon })\frac{{\partial {{\mathring{T}^{[2\varepsilon]}}}}}{{\partial {x_j}}}}\frac{\partial \mathring{T}^{[2\varepsilon]}}{{\partial {x_i}}}d\Omega\\
&-\frac{1}{2}\int_{\Omega}\frac{\partial\rho^{\varepsilon}({\bm{x}},\!{T^\varepsilon })}{\partial t}c^{\varepsilon}({\bm{x}},\!{T^\varepsilon })\big(\mathring{T}^{[2\varepsilon]}\big)^2d\Omega-\frac{1}{2}\int_{\Omega}\rho^{\varepsilon}({\bm{x}},\!{T^\varepsilon })\frac{\partial c^{\varepsilon}({\bm{x}},\!{T^\varepsilon })}{\partial t}\big(\mathring{T}^{[2\varepsilon]}\big)^2d\Omega\\
&-\int_{\Omega}{{\lambda_{ij}^{\varepsilon}}({\bm{x}},{T^\varepsilon })\frac{\partial {\mathring{\mathit{\Phi}}^{[2\varepsilon]}}}{\partial {x_i}}}\frac{\partial {\mathring{\mathit{\Phi}}^{[2\varepsilon]}}}{\partial {x_j}}\mathring{T}^{[2\varepsilon]}d\Omega+\frac{1}{2}\frac{\partial}{\partial t}\Big[\int_{\Omega}\rho^{\varepsilon}({\bm{x}},{T^\varepsilon })\big(\frac{\partial \mathring{U}_{i}^{[2\varepsilon]}}{\partial t}\big)^2\widetilde Td\Omega\Big]\\
&+\frac{1}{2}\frac{\partial}{\partial t}\Big[\int_{\Omega} {{c_{ijkl}^{\varepsilon}}({\bm{x}},{T^\varepsilon })\frac{{\partial \mathring{U}_{k}^{[2\varepsilon]}}}{{\partial {x_l}}}}\frac{\partial \mathring{U}_{i}^{[2\varepsilon]}}{{\partial {x_j}}}\widetilde Td\Omega\Big]-\frac{1}{2}\int_{\Omega}\frac{\partial\rho^{\varepsilon}({\bm{x}},{T^\varepsilon })}{\partial t}\big(\frac{\partial \mathring{U}_{i}^{[2\varepsilon]}}{\partial t}\big)^2\widetilde Td\Omega\\
&-\frac{1}{2}\int_{\Omega}{\frac{\partial c_{ijkl}^{\varepsilon}({\bm{x}},{T^\varepsilon })}{\partial t}\frac{{\partial \mathring{U}_{k}^{[2\varepsilon]}}}{{\partial {x_l}}}}\frac{\partial \mathring{U}_{i}^{[2\varepsilon]}}{{\partial {x_j}}}\widetilde Td\Omega=\int_{\Omega}\varepsilon H_2(\bm{x},\bm{y},t)\mathring{T}^{[2\varepsilon]}d\Omega\\
&+\int_{\Omega}\varepsilon F_{i2}(\bm{x},\bm{y},t)\frac{\partial \mathring{U}_{i}^{[2\varepsilon]}}{\partial t}\widetilde Td\Omega.
\end{aligned}
\end{equation}
Afterwards, integrating both sides of (3.15) over $(0,t]$ with $0<t\leq \mathcal T$ and imposing the initial conditions of thermal and mechanical residual equations given in (3.4), it follows that
\begin{equation}
\begin{aligned}
&\int_{\Omega}\rho^{\varepsilon}c^{\varepsilon}\big(\mathring{T}^{[2\varepsilon]}(\bm{x},t)\big)^2d\Omega+\int_0^t\int_{\Omega}2{{k_{ij}^{\varepsilon}}\frac{{\partial {{\mathring{T}^{[2\varepsilon]}}(\bm{x},\tau)}}}{{\partial {x_j}}}}\frac{\partial \mathring{T}^{[2\varepsilon]}(\bm{x},\tau)}{{\partial {x_i}}}d\Omega d\tau\\
&+\int_{\Omega}\rho^{\varepsilon}\big(\frac{\partial \mathring{U}_{i}^{[2\varepsilon]}(\bm{x},t)}{\partial t}\big)^2\widetilde Td\Omega+\int_{\Omega} {{c_{ijkl}^{\varepsilon}}\frac{{\partial \mathring{U}_{k}^{[2\varepsilon]}(\bm{x},t)}}{{\partial {x_l}}}}\frac{\partial \mathring{U}_{i}^{[2\varepsilon]}(\bm{x},t)}{{\partial {x_j}}}\widetilde Td\Omega\\
&=\int_{\Omega}\rho^{\varepsilon}c^{\varepsilon}\big(\varepsilon \widetilde\chi_2(\bm{x})\big)^2d\Omega+\int_0^t\int_{\Omega}2{{\lambda_{ij}^{\varepsilon}}\frac{\partial {\mathring{\mathit{\Phi}}^{[2\varepsilon]}}(\bm{x},\tau)}{\partial {x_i}}}\frac{\partial {\mathring{\mathit{\Phi}}^{[2\varepsilon]}}(\bm{x},\tau)}{\partial {x_j}}\mathring{T}^{[2\varepsilon]}(\bm{x},\tau)d\Omega d\tau\\
&+\int_0^t\int_{\Omega}\frac{\partial\rho^{\varepsilon}}{\partial t}c^{\varepsilon}\big(\mathring{T}^{[2\varepsilon]}(\bm{x},\tau)\big)^2d\Omega d\tau+\int_0^t\int_{\Omega}\rho^{\varepsilon}\frac{\partial c^{\varepsilon}}{\partial t}\big(\mathring{T}^{[2\varepsilon]}(\bm{x},\tau)\big)^2d\Omega d\tau\\
&+\int_{\Omega}\rho^{\varepsilon}\big(\varepsilon \widetilde\zeta_{2i}(\bm{x})\big)^2\widetilde Td\Omega+\int_0^t\int_{\Omega}{\frac{\partial c_{ijkl}^{\varepsilon}}{\partial t}\frac{{\partial \mathring{U}_{k}^{[2\varepsilon]}(\bm{x},\tau)}}{{\partial {x_l}}}}\frac{\partial \mathring{U}_{i}^{[2\varepsilon]}(\bm{x},\tau)}{{\partial {x_j}}}\widetilde Td\Omega d\tau\\
&+\int_0^t\int_{\Omega}\frac{\partial\rho^{\varepsilon}}{\partial t}\big(\frac{\partial \mathring{U}_{i}^{[2\varepsilon]}(\bm{x},\tau)}{\partial t}\big)^2\widetilde Td\Omega d\tau+\int_{\Omega} {{c_{ijkl}^{\varepsilon}}\frac{{\partial [\varepsilon \widetilde\eta_{2k}(\bm{x})]}}{{\partial {x_l}}}}\frac{\partial [\varepsilon \widetilde\eta_{2i}(\bm{x})]}{{\partial {x_j}}}\widetilde Td\Omega\\
&+\int_0^t\int_{\Omega}2\varepsilon H_2\mathring{T}^{[2\varepsilon]}(\bm{x},\tau)d\Omega d\tau+\int_0^t\int_{\Omega}2\varepsilon F_{i2}\frac{\partial \mathring{U}_{i}^{[2\varepsilon]}(\bm{x},\tau)}{\partial t}\widetilde Td\Omega d\tau.
\end{aligned}
\end{equation}
Owing to the assumptions (A$_1$) and (A$_2$), and employing the Poincar$\rm{\acute{e}}$-Friedrichs inequality and Korn's inequality, the following inequality can be obtained from left side of equality (3.16).
\begin{equation}
\begin{aligned}
&\int_{\Omega}\rho^{\varepsilon}c^{\varepsilon}\big(\mathring{T}^{[2\varepsilon]}(\bm{x},t)\big)^2d\Omega+\int_0^t\int_{\Omega}2{{k_{ij}^{\varepsilon}}\frac{{\partial {{\mathring{T}^{[2\varepsilon]}}(\bm{x},\tau)}}}{{\partial {x_j}}}}\frac{\partial \mathring{T}^{[2\varepsilon]}(\bm{x},\tau)}{{\partial {x_i}}}d\Omega d\tau\\
&+\int_{\Omega}\rho^{\varepsilon}\big(\frac{\partial \mathring{U}_{i}^{[2\varepsilon]}(\bm{x},t)}{\partial t}\big)^2\widetilde Td\Omega+\int_{\Omega} {{c_{ijkl}^{\varepsilon}}\frac{{\partial \mathring{U}_{k}^{[2\varepsilon]}(\bm{x},t)}}{{\partial {x_l}}}}\frac{\partial \mathring{U}_{i}^{[2\varepsilon]}(\bm{x},t)}{{\partial {x_j}}}\widetilde Td\Omega\\
&\geq\rho_*c_*\left \|\mathring{T}^{[2\varepsilon]}(\bm{x},t)\right\|_{L^2(\Omega)}^2
+C_1\int_0^t\left \|\mathring{T}^{[2\varepsilon]}(\bm{x},\tau)\right\|_{H_0^1(\Omega)}^2d\tau\\
&+\rho_*\widetilde T\Big\|\frac{\partial\mathring{\bm{U}}^{[2\varepsilon]}(\bm{x},t)}{\partial t}\Big\|_{(L^2(\Omega))^n}+C_2\widetilde T\big\|\mathring{\bm{U}}^{[2\varepsilon]}(\bm{x},t)\big\|_{(H^1_0(\Omega))^n}\\
&\geq \underline{C}\Big(\left \|\mathring{T}^{[2\varepsilon]}(\bm{x},t)\right\|_{L^2(\Omega)}^2
+\int_0^t\left \|\mathring{T}^{[2\varepsilon]}(\bm{x},\tau)\right\|_{H_0^1(\Omega)}^2d\tau\\
&+\Big\|\frac{\partial\mathring{\bm{U}}^{[2\varepsilon]}(\bm{x},t)}{\partial t}\Big\|_{(L^2(\Omega))^n}^2+\big\|\mathring{\bm{U}}^{[2\varepsilon]}(\bm{x},t)\big\|_{(H^1_0(\Omega))^n}^2\Big).
\end{aligned}
\end{equation}
Here $C_1$ and $C_2$ arise from the Poincar$\rm{\acute{e}}$-Friedrichs inequality and Korn's inequality separately. Moreover, we denote $\underline{C}=min(\rho_{*}c_{*},C_1,\rho_{*}\widetilde T,C_2\widetilde T)$. Subsequently, utilizing the assumptions (A$_1$), (A$_2$), (B$_2$), Schwarz's inequality, Young's inequality and error estimate for electric potential field, the following inequality is derived from the right side of equality (3.16).
\begin{equation}
\begin{aligned}
&\int_{\Omega}\rho^{\varepsilon}c^{\varepsilon}\big(\varepsilon \widetilde\chi_2(\bm{x})\big)^2d\Omega+\int_0^t\int_{\Omega}2{{\lambda_{ij}^{\varepsilon}}\frac{\partial {\mathring{\mathit{\Phi}}^{[2\varepsilon]}}(\bm{x},\tau)}{\partial {x_i}}}\frac{\partial {\mathring{\mathit{\Phi}}^{[2\varepsilon]}}(\bm{x},\tau)}{\partial {x_j}}\mathring{T}^{[2\varepsilon]}(\bm{x},\tau)d\Omega d\tau\\
&+\int_0^t\int_{\Omega}\frac{\partial\rho^{\varepsilon}}{\partial t}c^{\varepsilon}\big(\mathring{T}^{[2\varepsilon]}(\bm{x},\tau)\big)^2d\Omega d\tau+\int_0^t\int_{\Omega}\rho^{\varepsilon}\frac{\partial c^{\varepsilon}}{\partial t}\big(\mathring{T}^{[2\varepsilon]}(\bm{x},\tau)\big)^2d\Omega d\tau\\
&+\int_{\Omega}\rho^{\varepsilon}\big(\varepsilon \widetilde\zeta_{2i}(\bm{x})\big)^2\widetilde Td\Omega+\int_0^t\int_{\Omega}{\frac{\partial c_{ijkl}^{\varepsilon}}{\partial t}\frac{{\partial \mathring{U}_{k}^{[2\varepsilon]}(\bm{x},\tau)}}{{\partial {x_l}}}}\frac{\partial \mathring{U}_{i}^{[2\varepsilon]}(\bm{x},\tau)}{{\partial {x_j}}}\widetilde Td\Omega d\tau\\
&+\int_0^t\int_{\Omega}\frac{\partial\rho^{\varepsilon}}{\partial t}\big(\frac{\partial \mathring{U}_{i}^{[2\varepsilon]}(\bm{x},\tau)}{\partial t}\big)^2\widetilde Td\Omega d\tau+\int_{\Omega} {{c_{ijkl}^{\varepsilon}}\frac{{\partial [\varepsilon \widetilde\eta_{2k}(\bm{x})]}}{{\partial {x_l}}}}\frac{\partial [\varepsilon \widetilde\eta_{2i}(\bm{x})]}{{\partial {x_j}}}\widetilde Td\Omega\\
&+\int_0^t\int_{\Omega}2\varepsilon H_2\mathring{T}^{[2\varepsilon]}(\bm{x},\tau)d\Omega d\tau+\int_0^t\int_{\Omega}2\varepsilon F_{i2}\frac{\partial \mathring{U}_{i}^{[2\varepsilon]}(\bm{x},\tau)}{\partial t}\widetilde Td\Omega d\tau\\
&\leq C_3\varepsilon^2+C_4\int_0^t\left \|\mathring{T}^{[2\varepsilon]}(\bm{x},\tau)\right\|_{L^2(\Omega)}^2d\tau+C_5\int_0^t\left \|\mathring{\bm{U}}^{[2\varepsilon]}(\bm{x},\tau)\right\|_{(H_0^1(\Omega))^n}^2d\tau\\
&+C_6\int_0^t\Big\|\frac{\partial \mathring{\bm{U}}^{[2\varepsilon]}(\bm{x},\tau)}{\partial t}\Big\|_{(L^2(\Omega))^n}^2d\tau+C_7\int_0^t\left \|\mathring{\mathit{\Phi}}^{[2\varepsilon]}(\bm{x},\tau)\right\|_{H_0^1(\Omega)}^2d\tau\\
&\leq \overline{C}\varepsilon^2+\overline{C}\Bigg[\int_0^t\left \|\mathring{T}^{[2\varepsilon]}(\bm{x},\tau)\right\|_{L^2(\Omega)}^2d\tau+\int_0^t\int_0^{\tau}\left \|\mathring{T}^{[2\varepsilon]}(\bm{x},s)\right\|_{H_0^1(\Omega)}^2dsd\tau\\
&+\int_0^t\Big\|\frac{\partial \mathring{\bm{U}}^{[2\varepsilon]}(\bm{x},\tau)}{\partial t}\Big\|_{(L^2(\Omega))^n}^2d\tau+\int_0^t\left \|\mathring{\bm{U}}^{[2\varepsilon]}(\bm{x},\tau)\right\|_{(H_0^1(\Omega))^n}^2d\tau\Bigg].
\end{aligned}
\end{equation}
Here, we denote $\overline{C}=max(C_3,C_4,C_5,C_6,C_7)$. Furthermore, combining (3.17) and (3.18), it is apparent that
\begin{equation}
\begin{aligned}
&\underline{C}\Big(\left \|\mathring{T}^{[2\varepsilon]}(\bm{x},t)\right\|_{L^2(\Omega)}^2
+\int_0^t\left \|\mathring{T}^{[2\varepsilon]}(\bm{x},\tau)\right\|_{H_0^1(\Omega)}^2d\tau\\
&+\Big\|\frac{\partial\mathring{\bm{U}}^{[2\varepsilon]}(\bm{x},t)}{\partial t}\Big\|_{(L^2(\Omega))^n}+\big\|\mathring{\bm{U}}^{[2\varepsilon]}(\bm{x},t)\big\|_{(H^1_0(\Omega))^n}\Big)\\
&\leq \overline{C}\varepsilon^2+\overline{C}\Bigg[\int_0^t\left \|\mathring{T}^{[2\varepsilon]}(\bm{x},\tau)\right\|_{L^2(\Omega)}^2d\tau+\int_0^t\int_0^{\tau}\left \|\mathring{T}^{[2\varepsilon]}(\bm{x},s)\right\|_{H_0^1(\Omega)}^2dsd\tau\\
&+\int_0^t\Big\|\frac{\partial \mathring{\bm{U}}^{[2\varepsilon]}(\bm{x},\tau)}{\partial t}\Big\|_{(L^2(\Omega))^n}^2d\tau+\int_0^t\left \|\mathring{\bm{U}}^{[2\varepsilon]}(\bm{x},\tau)\right\|_{(H_0^1(\Omega))^n}^2d\tau\Bigg].
\end{aligned}
\end{equation}
Without loss of generality, let us introduce a constant $C=\overline{C}/\underline{C}$ and define a function $\displaystyle\Upsilon(t)=\left \|\mathring{T}^{[2\varepsilon]}\right\|_{L^2(\Omega)}^2
+\int_0^t\left \|\mathring{T}^{[2\varepsilon]}\right\|_{H_0^1(\Omega)}^2d\tau+\Big\|\frac{\partial\mathring{\bm{U}}^{[2\varepsilon]}}{\partial t}\Big\|_{(L^2(\Omega))^n}^2+\big\|\mathring{\bm{U}}^{[2\varepsilon]}\big\|_{(H^1_0(\Omega))^n}^2$. Subsequently, we have $\Upsilon(t)\leq C(\Omega)\varepsilon^2+C(\Omega)\mathlarger{\int}_0^t\Upsilon(\tau)d\tau$ according to (3.19). Next, it follows from the Gronwall's inequality in chapter 12 of reference \cite{R14} that $\Upsilon(t)\leq C(\Omega,\mathcal T)\varepsilon^2$, which yields the following inequality
\begin{equation}
\begin{aligned}
&\left\|\mathring{T}^{[2\varepsilon]}\right\|_{L^2(\Omega)}^2
+\int_0^t\left \|\mathring{T}^{[2\varepsilon]}\right\|_{H_0^1(\Omega)}^2d\tau\\
&+\Big\|\frac{\partial\mathring{\bm{U}}^{[2\varepsilon]}}{\partial t}\Big\|_{(L^2(\Omega))^n}^2+\big\|\mathring{\bm{U}}^{[2\varepsilon]}\big\|_{(H^1_0(\Omega))^n}^2\leq C(\Omega,\mathcal T)\varepsilon^2.
\end{aligned}
\end{equation}
Then using the arithmetic and geometric means inequality to the left side of the inequality (3.20) and squaring root on both sides of the inequality (3.20), the following inequality can be readily derived.
\begin{equation}
\begin{aligned}
&\left\|\mathring{T}^{[2\varepsilon]}\right\|_{L^2(\Omega)}
+\big\|\mathring{T}^{[2\varepsilon]}\big\|_{L^2(0,t;H^1_0(\Omega))}\\
&+\Big\|\frac{\partial\mathring{\bm{U}}^{[2\varepsilon]}}{\partial t}\Big\|_{(L^2(\Omega))^n}+\big\|\mathring{\bm{U}}^{[2\varepsilon]}\big\|_{(H^1_0(\Omega))^n}\leq C(\Omega,\mathcal T)\varepsilon.
\end{aligned}
\end{equation}
To account for the arbitrariness of time variable $t$, we obtain the explicit convergence estimate (3.5) from (3.10) and (3.21).
\section{Two-stage numerical algorithm}
For the investigated multi-scale nonlinear multi-physics problem (1.1), we propose a high-accuracy multi-scale computational framework that yields a closed solving system, consisting of microscopic cell models, a macroscopic homogenized model, and high-order multi-scale solutions. Leveraging the inherent characteristics of this multi-scale computational model, a two-stage numerical algorithm, comprising off-line and on-line stages, is developed to enable efficient multi-scale simulation of the time-dependent nonlinear thermo-electro-mechanical coupling problem (1.1) of composite structures, details of which are presented as follows.
\subsection{Off-line stage: computation for microscopic cell problems}
\begin{enumerate}
\item[(1)]
Identify the geometric configuration of PUC $\mathbf{Y}=[0,1]^n$ and generate a family of triangular or tetrahedral finite element meshes $\mathcal{T}_{h_1}=\{K\}$ for PUC $\mathbf{Y}$, where $h_1=$max$_K\{h_K\}$ denotes the mesh size. Whereupon denoting the linear conforming finite element space $\mathcal{V}_{h_1}(\mathbf{Y})=\{\nu\in C^0(\bar{\mathbf{Y}}):\nu\mid_{\partial \mathbf{Y}}=0,\nu\mid_{K}\in P_1(K)\}\subset H^1(\mathbf{Y})$ for auxiliary cell problems.
\item[(2)]
Define the computational temperature range $[T_0^{min},T_0^{max}]$ and select a set of representative macroscopic temperature $\bar T_0^{s}$ within the concerned temperature range. Next, employ FEM to solve the first-order cell functions defined by (2.20)-(2.23) on $\mathcal{V}_{h_1}(\mathbf{Y})$ corresponding to each distinct representative macroscopic temperature $\bar T_0^{s}$. Note that classical periodic boundary condition of auxiliary cell problems is replaced by homogeneous Dirichlet boundary condition for practical numerical implementation \cite{R5,R24,R37}. Below, the specific finite element scheme for first-order unit cell problem (2.20) is established.
\begin{equation}
\begin{aligned}
&\int_{\mathbf{Y}}{ k_{ij}^{(0)}(\bm{y},\bar T_0^{s}){\frac{\partial \mathcal{M}_{\alpha_1}(\bm{y},\bar T_0^{s})}{\partial y_j}}}\frac{\partial \upsilon^{h_1}}{\partial y_i}d\mathbf{Y}\\
&=-\int_{\mathbf{Y}}k_{i{\alpha_1}}^{(0)}(\bm{y},\bar T_0^{s})\frac{\partial \upsilon^{h_1}}{\partial y_i}d\mathbf{Y},\;\forall\upsilon^{h_1}\in \mathcal{V}_{h_1}(\mathbf{Y}).
\end{aligned}
\end{equation}
\item[(3)]
The macroscopic material parameters $\widehat{S}(T_0)$, $\widehat{k}_{ij}(T_0)$, $\widehat{\lambda}_{ij}(T_0)$, $\widehat{\rho}(T_0)$, $\widehat{c}_{ijkl}(T_0)$ and $\widehat\beta_{ij}(T_0)$ are evaluated via formula (2.25) for each distinct macroscopic temperature $\bar T_0^{s}$.
\item[(4)]
Utilizing the identical mesh as that employed for the first-order cell functions, the second-order auxiliary cell functions governed by (2.26)-(2.41), which correspond to distinct representative macroscopic temperature $\bar T_0^{s}$, are evaluated on $\mathcal{V}_{h_1}(\mathbf{Y})$ via FEM respectively.
\end{enumerate}
\subsection{On-line stage: computation for macroscopic homogenized problem and high-order multi-scale solutions}
\begin{enumerate}
\item[(1)]
Let $\mathcal{T}_{h_0}=\{e\}$ be a family of triangular or tetrahedral finite element meshes for the macroscopic region $\Omega$, where $h_0=$max$_e\{h_e\}$ denotes the mesh size. Then define the linear conforming finite element space $\mathcal{V}_{h_0}(\Omega)=\{\nu\in C^0(\bar{\Omega}):\nu\mid_{\partial\Omega}=0,\nu\mid_{e}\in P_1(e)\}\subset H^1(\Omega)$ for the macroscopic homogenized equations (2.24). Moreover, the homogenized material parameters can be determined through the interpolation approach at each node $\bm{x}$ of $\mathcal{V}_{h_0}(\Omega)$.
\item[(2)]
Solve the macroscopic homogenized equations (2.24) with homogeneous coefficients on the computational domain $\Omega\times(0,\mathcal T)$ using a mixed FDM-FEM inspired by references \cite{R11,R35,R36,R38}, which is implemented on a coarser mesh with a larger time step. This approach employs FEM for spatial discretization and and FDM for temporal discretization. Using the equidistant time step $\displaystyle\Delta t={\mathcal T}/{N}$ to divide time-domain $(0,\mathcal T)$ into $0=t_0<t_1<\cdots<t_N=\mathcal T$ and $t_m=m\Delta t(m=0,\cdots,N)$, then we denote $T_0^{m}=T_0(\bm{x},t_m)$, $\mathit{\Phi}_0^{m}=\mathit{\Phi}_0(\bm{x},t_m)$ and $U_{i0}^{m}=U_{i0}(\bm{x},t_m)$. Additionally, define $\widehat T_{0}^{m+1/2}=(3T_{0}^{m}-T_{0}^{m-1})/2$ and $\bar T_{0}^{m+1/2}=(T_{0}^{m}+T_{0}^{m+1})/2$ for $m=0,\cdots,N-1$. Then, the computational scheme inspired by reference \cite{R38} is employed to precompute $\mathit{\Phi}_0^0$ and $\widehat T_{0}^{1/2}$, details of which are presented below.
\begin{equation}
\left\{\begin{aligned}
&\int_{\Omega}{{\widehat \lambda _{ij}}(T_{0}^0)\frac{{\partial \mathit{\Phi}_0^0}}{{\partial {x_j}}}} \frac{\partial \varphi}{{\partial {x_i}}}d\Omega=\int_{\Omega}f_{\mathit{\Phi}}(\bm{x},t_0)\varphi d\Omega,\;\;\forall\varphi\in \mathcal{V}_{h_0}(\Omega),\\
&\mathit{\Phi}_0^0=\widehat{\mathit{\Phi}}(\bm{x},t_0),\;\;\text{on}\;\;\partial\Omega.
\end{aligned}\right.
\end{equation}
\begin{equation}
\left\{ \begin{aligned}
&\int_{\Omega}\widehat S(T_{0}^0)\frac{{\widehat T_{0}^{1/2}-T_{0}^{0}}}{{\Delta t/2}}\upsilon d\Omega+ \int_{\Omega}{{\widehat k_{ij}}(T_{0}^0)\frac{{\partial {\widehat T_{0}^{1/2}}}}{{\partial {x_j}}}}\frac{\partial \upsilon}{{\partial {x_i}}}d\Omega\\
&=\int_{\Omega} {\widehat \lambda _{ij}^*}(T_{0}^0)\frac{\partial \mathit{\Phi}_{0}^0}{\partial {x_i}}\frac{\partial \mathit{\Phi}_{0}^0}{\partial {x_j}}\upsilon d\Omega-\int_{\Omega} T_0^0\widehat \beta_{ij}^*(T_{0}^0)\frac{\partial \widetilde{{U}}_i^1}{\partial x_j}\upsilon d\Omega\\
&+\int_{\Omega} f_T(\bm{x},t_{1/2})\upsilon d\Omega,\;\;\forall\upsilon\in \mathcal{V}_{h_0}(\Omega),\\
&\widehat T_0^{1/2} = \widehat T(\bm{x},t_{1/2}),\;\;\text{on}\;\;\partial {\Omega}.
\end{aligned} \right.
\end{equation}
Next, the following computational scheme inspired by references \cite{R11,R35,R36,R38} is adopted to successively compute $\mathit{\Phi}_0^{m+1/2}$, $T_{0}^{m+1}$ and $U_{i0}^{m+1}$ for $m=0,\cdots,N-1$.
\begin{equation}
\left\{\begin{aligned}
&\!\!\!\!\int_{\Omega}\!{{\widehat \lambda _{ij}}(\widehat T_{0}^{m+1/2})\frac{{\partial \mathit{\Phi}_0^{m+1/2}}}{{\partial {x_j}}}}\!\frac{\partial \varphi}{{\partial {x_i}}}d\Omega\!=\!\!\!\int_{\Omega}\!f_{\mathit{\Phi}}(\bm{x},t_{m+1/2})\varphi d\Omega,\;\;\forall\varphi\!\in\! \mathcal{V}_{h_0}\!(\Omega),\\
&\!\!\!\mathit{\Phi}_0^{m+1/2}\!=\!\widehat{\mathit{\Phi}}(\bm{x},t_{m+1/2}),\;\;\text{on}\;\;\partial\Omega.
\end{aligned}\right.
\end{equation}
\begin{equation}
\left\{ \begin{aligned}
&\int_{\Omega}\widehat S(\widehat T_{0}^{m+1/2})\frac{{T_{0}^{m+1}-T_{0}^{m}}}{{\Delta t}}\upsilon d\Omega+\int_{\Omega}{{\widehat k_{ij}}(\widehat T_{0}^{m+1/2})\frac{{\partial {\bar T_{0}^{m+1/2}}}}{{\partial {x_j}}}}\frac{\partial \upsilon}{{\partial {x_i}}}d\Omega\\
&= \int_{\Omega}{\widehat \lambda _{ij}^*}(\widehat T_{0}^{m+1/2})\frac{\partial \mathit{\Phi}_{0}^{m+1/2}}{\partial {x_i}}\frac{\partial \mathit{\Phi}_{0}^{m+1/2}}{\partial {x_j}}\upsilon d\Omega\\
&-\int_{\Omega} \widehat T_{0}^{m+1/2}\widehat \beta_{ij}^*(\widehat T_{0}^{m+1/2})\frac{1}{{\Delta t}}\Big(\frac{\partial {{U}}_{i0}^{m}}{\partial x_j}-\frac{\partial {{U}}_{i0}^{m-1}}{\partial x_j}\Big)\upsilon d\Omega\\
&+\int_{\Omega}f_T(\bm{x},t_{m+1/2})\upsilon d\Omega,\;\;\forall\upsilon\in \mathcal{V}_{h_0}(\Omega),\\
&T_0^{m+1} = \widehat T(\bm{x},t_{m+1}),\;\;\text{on}\;\;\partial {\Omega}.
\end{aligned} \right.
\end{equation}
\begin{equation}
\left\{ \begin{aligned}
&\int_{\Omega}\widehat \rho(T_{0}^{m+1})\frac{{U_{i0}^{m+1}-2U_{i0}^{m}+U_{i0}^{m-1}}}{{(\Delta t)^2}}\nu_i d\Omega\\
&+\int_{\Omega}{{\widehat c_{ijkl}}(T_{0}^{m+1})\frac{{\partial {U_{k0}^{m+1}}}}{{\partial {x_l}}}}\frac{\partial \nu_i}{{\partial {x_j}}}d\Omega\\
&-\int_{\Omega}\widehat \beta_{ij}(T_{0}^{m+1})(T_0^{m+1}-\widetilde T)\frac{\partial \nu_i}{{\partial {x_j}}}d\Omega\\
&+\int_{\Omega}f_i(\bm{x},t_{m+1})\upsilon d\Omega,\;\;\forall\bm{\nu}\in \big(\mathcal{V}_{h_0}(\Omega)\big)^n,\\
&\bm{U}_0^{m+1} = \widehat {\bm{U}}(\bm{x},t_{m+1}),\;\;\text{on}\;\;\partial {\Omega}.
\end{aligned} \right.
\end{equation}
Consequently, we can solve for the macroscopic electric potential, temperature and displacement fields at each temporal step by means of three alternating sub-problems (4.5)-(4.6).
\item[(3)]
At any given point $(\bm{x},t)\in \Omega\times(0,\mathcal T)$, we utilize the interpolation method to determine the corresponding values of the first-order cell functions, second-order cell functions, and homogenized solutions.
\item[(4)]
The average approach on relative elements \cite{R37} is utilized to compute the spatial derivatives, and the difference scheme is applied to approximate the temporal derivatives at each time steps involving in formulas (2.17)-(2.19).
\item[(5)]
Ultimately, the temperature field $T^{[2\varepsilon]}(\bm{x},t)$, electric potential field $\mathit{\Phi}^{[2\varepsilon]}(\bm{x},t)$ and displacement field $\bm{U}^{[2\varepsilon]}(\bm{x},t)$ are computed  separately via the formulas (2.17)-(2.19). Furthermore, we can employ the high-order interpolation and post-processing techniques to gain high-accuracy HOMS solutions \cite{R39}.
\end{enumerate}
\section{The error analyses for two-stage numerical algorithm}
\begin{lemma}
Let $M_{\alpha_1}^{h_1}$ and $N_{\alpha_1}^{h_1}$ denote the corresponding finite element solutions for the microscopic cell functions, respectively. If all microscopic cell functions belong to $H^2(\mathcal{Y})$ for any fixed $T_0$, then the following inequalities hold
\begin{equation}
\begin{aligned}
&\left\|\mathcal{M}_{\alpha_1}^{h_1}(\bm{y},T_0)-\mathcal{M}_{\alpha_1}(\bm{y},T_0)\right\|_{H^m(\mathbf{Y})}\leq Ch_1^{2-m}\left\| \mathcal{M}_{\alpha_1}(\bm{y},T_0)\right\|_{H^2(\mathbf{Y})},\\
&\left\|\mathcal{H}_{\alpha_1}^{h_1}(\bm{y},T_0)-\mathcal{H}_{\alpha_1}(\bm{y},T_0)\right\|_{H^m(\mathbf{Y})}\leq Ch_1^{2-m}\left\|\mathcal{H}_{\alpha_1}(\bm{y},T_0)\right\|_{H^2(\mathbf{Y})},\\
&\left\|\mathcal{N}_{km}^{\alpha_1,h_1}(\bm{y},T_0)-\mathcal{N}_{km}^{\alpha_1}(\bm{y},T_0)\right\|_{H^m(\mathbf{Y})}\leq Ch_1^{2-m}\left\| \mathcal{N}_{km}^{\alpha_1}(\bm{y},T_0)\right\|_{H^2(\mathbf{Y})},\\
&\left\|\mathcal{P}_{k}^{h_1}(\bm{y},T_0)-\mathcal{P}_{k}(\bm{y},T_0)\right\|_{H^m(\mathbf{Y})}\leq Ch_1^{2-m}\left\|\mathcal{P}_{k}(\bm{y},T_0)\right\|_{H^2(\mathbf{Y})},
\end{aligned}
\end{equation}
where $m=0,1$ and $C$ is the finite element estimate constant independent of $h_1$ and dependent on $\mathbf{Y}$.
\end{lemma}
$\mathbf{Proof:}$ These inequalities can be easily derived by employing classical finite element theory.
\begin{lemma}
Denote $\widehat S^{h_1}(T_{0})$, $\widehat k_{ij}^{h_1}(T_{0})$, $\widehat \lambda_{ij}^{h_1}(T_{0})$, $\widehat c_{ijkl}^{h_1}(T_{0})$ and $\widehat \beta_{ij}^{h_1}(T_{0})$ be the finite element approximations of the corresponding homogenized parameters, the following results hold
\begin{equation}
\begin{aligned}
&\left|\widehat S^{h_1}(T_{0})-\widehat S(T_{0})\right|\leq Ch_1^2\left\| \bm{\mathcal{P}}(\bm{y},T_0)\right\|_{(H^2(\mathbf{Y}))^n}^2,\\
&\widetilde S_0\leq \widehat S^{h_1}(T_{0})\leq \widetilde S_1,
\end{aligned}
\end{equation}
\begin{equation}
\begin{aligned}
&\left|\widehat k_{ij}^{h_1}(T_{0})-\widehat k_{ij}(T_{0})\right|\leq Ch_1^2\left\| \mathcal{M}_i(\bm{y},T_0)\right\|_{H^2(\mathbf{Y})}\left\| \mathcal{M}_j(\bm{y},T_0)\right\|_{H^2(\mathbf{Y})},\\
&\widetilde\kappa_0|\bm{\xi}|^2\leq \widehat k_{ij}^{h_1}(T_{0})\xi_i\xi_j\leq\widetilde\kappa_1|\bm{\xi}|^2,
\end{aligned}
\end{equation}
\begin{equation}
\begin{aligned}
&\left|\widehat \lambda_{ij}^{h_1}(T_{0})-\widehat \lambda_{ij}(T_{0})\right|\leq Ch_1^2\left\| \mathcal{H}_i(\bm{y},T_0)\right\|_{H^2(\mathbf{Y})}\left\| \mathcal{H}_j(\bm{y},T_0)\right\|_{H^2(\mathbf{Y})},\\
&\widetilde\kappa_0|\bm{\xi}|^2\leq \widehat \lambda_{ij}^{h_1}(T_{0})\xi_i\xi_j\leq\widetilde\kappa_1|\bm{\xi}|^2,
\end{aligned}
\end{equation}
\begin{equation}
\begin{aligned}
&\left|\widehat c_{ijkl}^{h_1}(T_{0})-\widehat c_{ijkl}(T_{0})\right|\leq Ch_1^2\left\| \bm{\mathcal{N}}_{i}^j(\bm{y},T_0)\right\|_{(H^2(\mathbf{Y}))^n}\left\| \bm{\mathcal{N}}_{k}^l(\bm{y},T_0)\right\|_{(H^2(\mathbf{Y}))^n},\\
&\widetilde\kappa_0|\bm{\zeta}|^2\leq \widehat c_{ijkl}^{h_1}(T_{0})\zeta_{ij}\zeta_{kl}\leq\widetilde\kappa_1|\bm{\zeta}|^2,
\end{aligned}
\end{equation}
\begin{equation}
\begin{aligned}
&\left|\widehat \beta_{ij}^{h_1}(T_{0})-\widehat \beta_{ij}(T_{0})\right|\leq Ch_1^2\left\| \bm{\mathcal{N}}_{i}^j(\bm{y},T_0)\right\|_{(H^2(\mathbf{Y}))^n}\left\| \bm{\mathcal{P}}(\bm{y},T_0)\right\|_{(H^2(\mathbf{Y}))^n}^2,\\
&\widetilde\kappa_0|\bm{\xi}|^2\leq \widehat \beta_{ij}^{h_1}(T_{0})\xi_i\xi_j\leq\widetilde\kappa_1|\bm{\xi}|^2,
\end{aligned}
\end{equation}
where $C$ is a constant independent of $h_1$.
\end{lemma}
$\mathbf{Proof:}$ Inspired by reference \cite{R40}, by employing the definitions of macroscopic homogenized material parameters given in (2.25), along with assumption $A_1$, remark 2 and lemma 5.1, it follows that
\begin{equation}
\begin{aligned}
&\left|\widehat c_{ijkl}^{h_1}(T_{0})-\widehat c_{ijkl}(T_{0})\right|\\
&=\left|\frac{1}{|\mathbf{Y}|}{\int_{\mathbf{Y}}}\big({c_{ijkl}^{(0)} + c_{ij\alpha_1\alpha_2}^{(0)}{\frac{\partial \mathcal{N}_{\alpha_1k}^{l,h_1}}{\partial y_{\alpha_2}}}}\big)d\mathbf{Y}-\frac{1}{|\mathbf{Y}|}{\int_{\mathbf{Y}}}\big({c_{ijkl}^{(0)} + c_{ij\alpha_1\alpha_2}^{(0)}{\frac{\partial \mathcal{N}_{\alpha_1k}^{l}}{\partial y_{\alpha_2}}}}\big)d\mathbf{Y}\right|\\
&=\left|\frac{1}{|\mathbf{Y}|}{\int_{\mathbf{Y}}}{c_{ij\alpha_1\alpha_2}^{(0)}{\frac{\partial \big(\mathcal{N}_{\alpha_1k}^{l,h_1}-\mathcal{N}_{\alpha_1k}^{l}\big)}{\partial y_{\alpha_2}}}}d\mathbf{Y}\right|\\
&=\frac{1}{|\mathbf{Y}|}\left|{\int_{\mathbf{Y}}}{\frac{\partial \mathcal{N}_{\alpha_3i}^{j}}{\partial y_{\alpha_4}}}{c_{\alpha_3\alpha_4\alpha_1\alpha_2}^{(0)}{\frac{\partial \big(\mathcal{N}_{\alpha_1k}^{l,h_1}-\mathcal{N}_{\alpha_1k}^{l}\big)}{\partial y_{\alpha_2}}}}d\mathbf{Y}\right|\\
&=\frac{1}{|\mathbf{Y}|}\left|{\int_{\mathbf{Y}}}{\frac{\partial \big(\mathcal{N}_{\alpha_3i}^{j,h_1}-\mathcal{N}_{\alpha_3i}^{j}\big)}{\partial y_{\alpha_4}}}{c_{\alpha_3\alpha_4\alpha_1\alpha_2}^{(0)}{\frac{\partial \big(\mathcal{N}_{\alpha_1k}^{l,h_1}-\mathcal{N}_{\alpha_1k}^{l}\big)}{\partial y_{\alpha_2}}}}d\mathbf{Y}\right|\\
&\leq C\left\|\bm{\mathcal{N}}_{i}^{j,h_1}(\bm{y},T_0)\!-\!\bm{\mathcal{N}}_{i}^j(\bm{y},T_0)\right\|_{(H^1(\mathbf{Y}))^n}\left\| \bm{\mathcal{N}}_{k}^{l,h_1}(\bm{y},T_0)\!-\!\bm{\mathcal{N}}_{k}^l(\bm{y},T_0)\right\|_{(H^1(\mathbf{Y}))^n}\\
&\leq Ch_1^2\left\| \bm{\mathcal{N}}_{i}^j(\bm{y},T_0)\right\|_{(H^2(\mathbf{Y}))^n}\left\| \bm{\mathcal{N}}_{k}^l(\bm{y},T_0)\right\|_{(H^2(\mathbf{Y}))^n}.
\end{aligned}
\end{equation}
Choosing a sufficiently small $h_1>0$ ensures that $Ch_1^2\left\| \bm{\mathcal{N}}_{i}^j\right\|_{(H^2(\mathbf{Y}))^n}\left\| \bm{\mathcal{N}}_{k}^l\right\|_{(H^2(\mathbf{Y}))^n}\leq\bar\kappa_0/2$, from which we can verify that the lower bound stated in (5.5) holds as below.
\begin{equation}
\begin{aligned}
&\widehat c_{ijkl}^{h_1}(T_{0})\zeta_{ij}\zeta_{kl}=\widehat c_{ijkl}(T_{0})\zeta_{ij}\zeta_{kl}+\big[\widehat c_{ijkl}^{h_1}(T_{0})-\widehat c_{ijkl}(T_{0})\big]\zeta_{ij}\zeta_{kl}\\
&\geq(\bar\kappa_0-\bar\kappa_0/2)|\bm{\zeta}|^2=\widetilde\kappa_0|\bm{\zeta}|^2,
\end{aligned}
\end{equation}
where $\widetilde\kappa_0=\bar\kappa_0/2$ is a constant independent of $h_1$. Moreover, the upper bound in (5.5) is easily derived when setting $\widetilde\kappa_1=\bar\kappa_1+\bar\kappa_0/2$. Finally, following the similar way, we can obtain the results (5.2)-(5.4) and (5.6).

According to lemmas 5.1 and 5.2, the values of macroscopic homogenized material parameters depend on the finite element computation of the first-order auxiliary cell functions $\mathcal{M}_{\alpha_1}(\bm{y},T_0)$, $\mathcal{H}_{\alpha_1}(\bm{y},T_0)$, $\mathcal{N}_{km}^{\alpha_1}(\bm{y},T_0)$ and $\mathcal{P}_{k}(\bm{y},T_0)$. Thus, the following modified homogenized equations are numerically solved in practice.
\begin{equation}
\left\{ \begin{aligned}
&\widehat S^{h_1}(T_{0}^{h_1})\frac{{\partial {T_{0}^{h_1}}}}{{\partial t}}- \frac{\partial }{{\partial {x_i}}}\Big( {{\widehat k_{ij}^{h_1}}(T_{0}^{h_1})\frac{{\partial {T_{0}^{h_1}}}}{{\partial {x_j}}}}\Big)\\
&\quad={\widehat \lambda_{ij}^{h_1}}(T_{0}^{h_1})\frac{\partial \mathit{\Phi}_{0}^{h_1}}{\partial {x_i}}\frac{\partial \mathit{\Phi}_{0}^{h_1}}{\partial {x_j}}-T_0^{h_1}\widehat \beta_{ij}^{h_1}(T_{0}^{h_1})\frac{\partial}{\partial t}\Big(\frac{\partial U_{i0}^{h_1}}{\partial x_j}\Big)+f_T,\;\text{in}\;\Omega\times(0,\mathcal T),\\
&-\frac{\partial }{{\partial {x_i}}}\Big( {{\widehat \lambda _{ij}^{h_1}}(T_{0}^{h_1})\frac{{\partial \mathit{\Phi}_0^{h_1}}}{{\partial {x_j}}}}\Big)=f_\mathit{\Phi},\;\text{in}\;\Omega\times(0,\mathcal T),\\
&\widehat \rho(T_0^{h_1})\frac{{\partial^2 {U_{i0}^{h_1}}}}{{\partial t^2}}-\frac{\partial }{{\partial {x_j}}}\Big( {{\widehat c_{ijkl}^{h_1}}(T_0^{h_1})\frac{{\partial {U_{k0}^{h_1}}}}{{\partial {x_l}}}}\Big)\\
&\quad=-\frac{\partial }{{\partial {x_j}}}\Big(\widehat {\beta}_{ij}^{h_1}(T_0^{h_1})(T_0^{h_1}-\widetilde T)\Big)+f_i,\;\text{in}\;\Omega\times(0,\mathcal T),\\
&T_0^{h_1}(\bm{x},\!t) = \widehat T(\bm{x},\!t),\!\;\mathit{\Phi}_0^{h_1}(\bm{x},\!t)=\widehat{\mathit{\Phi}}(\bm{x},\!t),\!\;{\bm{U}}_0^{h_1}(\bm{x},\!t)=\widehat {\bm{U}}(\bm{x},\!t),\;\text{on}\;\partial\Omega\!\times\!(0,\!\mathcal T),\\
&T_0^{h_1}({\bm{x}},0)=\widetilde T,\;{\bm{U}}_0^{h_1}({\bm{x}},0)=\widetilde{\bm{U}}^0({\bm{x}}),\;\frac{\partial {\bm{U}}_0^{h_1}({\bm{x}},0)}{\partial t}=\widetilde{\bm{U}}^1({\bm{x}}),\;\text{in}\;\Omega.
\end{aligned} \right.
\end{equation}
Based on lemma 5.2, the solutions of the modified homogenized equations (5.9) possess the same regularity as those established in references \cite{R8,R9,R10,R11,R12}.
\begin{lemma}
Denote $T_{0}^{h_1}$, $\mathit{\Phi}_{0}^{h_1}$ and $\bm{U}_{0}^{h_1}$ be the exact solution of the modified homogenized equations (5.9), and assume that the two-way coupled system is reduced to a one-way coupled system in which the displacement field does not affect the temperature field, then the following estimate holds
\begin{equation}
{\begin{aligned}
&\big\|{\mathit{\Phi}_0^{h_1}}\!-\!{\mathit{\Phi}_0}\big\|_{L^\infty(0,T;H_0^1(\Omega))}+\big\|{T_0^{h_1}}\!-\!{T_0}\big\|_{L^\infty(0,\mathcal T;L^2(\Omega))}+\big\|{T_0^{h_1}}\!-\!{T_0}\big\|_{L^2(0,\mathcal T;H^1_0(\Omega))}\\
&+\Big\|\frac{\partial{\bm{U}_0^{h_1}}}{\partial t}-\frac{\partial{\bm{U}_0}}{\partial t}\Big\|_{L^\infty(0,\mathcal T;(L^2(\Omega))^n)}+\big\|{\bm{U}_0^{h_1}}-{\bm{U}_0}\big\|_{L^\infty(0,\mathcal T;(H^1_0(\Omega))^n)}\leq Ch_1^2,
\end{aligned}}
\end{equation}
where $C$ is a constant independent of $h_1$.
\end{lemma}
$\mathbf{Proof:}$ Through subtracting the electric equation in (2.24) from corresponding electric equation in (5.9), one can directly verify that
\begin{equation}
\begin{aligned}
&- \frac{\partial }{{\partial {x_i}}}\Big[ {{\widehat \lambda_{ij}^{h_1}}(T_{0}^{h_1})\frac{{\partial \big(\mathit{\Phi}_0^{h_1}-\mathit{\Phi}_0\big)}}{{\partial {x_j}}}}\Big]=-\frac{\partial }{{\partial {x_i}}}\Big[ {{\big(\widehat \lambda _{ij}(T_{0})-\widehat \lambda _{ij}^{h_1}(T_{0}^{h_1})\big)}\frac{{\partial \mathit{\Phi}_0}}{{\partial {x_j}}}}\Big]\\
&=-\frac{\partial }{{\partial {x_i}}}\Big[ {{\big(\widehat \lambda_{ij}(T_{0})-\widehat \lambda _{ij}^{h_1}(T_{0})\big)}\frac{{\partial \mathit{\Phi}_0}}{{\partial {x_j}}}}\Big]- \frac{\partial }{{\partial {x_i}}}\Big[ {{\big(\widehat \lambda_{ij}^{h_1}(T_{0})-\widehat \lambda_{ij}^{h_1}(T_{0}^{h_1})\big)}\frac{{\partial \mathit{\Phi}_0}}{{\partial {x_j}}}}\Big].
\end{aligned}
\end{equation}
Next, multiplying on both sides of equality (5.11) by $\mathit{\Phi}_0^{h_1}-\mathit{\Phi}_0$ and integrating over $\Omega$, we obtain the following equality.
\begin{equation}
\begin{aligned}
&\int_{\Omega}{{\widehat \lambda _{ij}^{h_1}}(T_{0}^{h_1})\frac{{\partial \big(\mathit{\Phi}_0^{h_1}-\mathit{\Phi}_0\big)}}{{\partial {x_j}}}}\frac{\partial \big(\mathit{\Phi}_0^{h_1}-\mathit{\Phi}_0\big)}{{\partial {x_i}}}d\Omega\\
&=\int_{\Omega}{{\big(\widehat \lambda _{ij}(T_{0})-\widehat \lambda _{ij}^{h_1}(T_{0})\big)}\frac{{\partial \mathit{\Phi}_0}}{{\partial {x_j}}}}\frac{\partial \big(\mathit{\Phi}_0^{h_1}-\mathit{\Phi}_0\big)}{{\partial {x_i}}}d\Omega\\
&+\int_{\Omega} {{\big(\widehat \lambda _{ij}^{h_1}(T_{0})-\widehat \lambda _{ij}^{h_1}(T_{0}^{h_1})\big)}\frac{{\partial\mathit{\Phi}_0}}{{\partial {x_j}}}}\frac{\partial \big(\mathit{\Phi}_0^{h_1}-\mathit{\Phi}_0\big)}{{\partial {x_i}}}d\Omega.
\end{aligned}
\end{equation}
Recalling the inequality in (5.4) and employing the Cauchy-Schwarz inequality, we naturally derive the following inequality from equality (5.14) if $|\widehat \lambda _{ij}^{h_1}(T_{0})-\widehat \lambda _{ij}^{h_1}(T_{0}^{h_1})|\leq C|T_{0}-T_{0}^{h_1}|$.
\begin{equation}
{\begin{aligned}
\big\|{\mathit{\Phi}_0^{h_1}}-{\mathit{\Phi}_0}\big\|_{H_0^1(\Omega)}\leq Ch_1^2+C\big\|{T_0^{h_1}}-{T_0}\big\|_{L^2(\Omega)}.
\end{aligned}}
\end{equation}

Afterwards, by subtracting the thermal equation in (2.24) from the corresponding thermal equation in (5.9) respectively, we shall obtain
\begin{equation}
\begin{aligned}
&\widehat S^{h_1}(T_{0}^{h_1})\frac{{\partial \big(T_0^{h_1}-T_0\big)}}{{\partial t}}- \frac{\partial }{{\partial {x_i}}}\Big[ {{\widehat k_{ij}^{h_1}}(T_{0}^{h_1})\frac{{\partial \big(T_0^{h_1}-T_0\big)}}{{\partial {x_j}}}}\Big]\\
&=\big[\widehat S(T_{0})-\widehat S^{h_1}(T_{0})\big]\frac{{\partial T_0}}{{\partial t}}+\big[\widehat S^{h_1}(T_{0})-\widehat S^{h_1}(T_{0}^{h_1})\big]\frac{{\partial T_0}}{{\partial t}}\\
&-\frac{\partial }{{\partial {x_i}}}\Big[ {{\big(\widehat k_{ij}(T_{0})-\widehat k_{ij}^{h_1}(T_{0})\big)}\frac{{\partial T_0}}{{\partial {x_j}}}}\Big]- \frac{\partial }{{\partial {x_i}}}\Big[ {{\big(\widehat k_{ij}^{h_1}(T_{0})-\widehat k_{ij}^{h_1}(T_{0}^{h_1})\big)}\frac{{\partial T_0}}{{\partial {x_j}}}}\Big]\\
&+ {\widehat \lambda _{ij}^{h_1}}(T_{0}^{h_1})\frac{\partial \big(\mathit{\Phi}_0^{h_1}-\mathit{\Phi}_0\big)}{\partial {x_i}}\frac{\partial \mathit{\Phi}_{0}^{h_1}}{\partial {x_j}}+ \big[{\widehat \lambda _{ij}^{h_1}}(T_{0}^{h_1})-{\widehat \lambda _{ij}^{h_1}}(T_{0})\big]\frac{\partial \mathit{\Phi}_0}{\partial {x_i}}\frac{\partial \mathit{\Phi}_{0}^{h_1}}{\partial {x_j}}\\
&+ \big[{\widehat \lambda _{ij}^{h_1}}(T_{0})-{\widehat \lambda _{ij}}(T_{0})\big]\frac{\partial \mathit{\Phi}_0}{\partial {x_i}}\frac{\partial \mathit{\Phi}_{0}^{h_1}}{\partial {x_j}}+ {\widehat \lambda _{ij}}(T_{0})\frac{\partial \mathit{\Phi}_0}{\partial {x_i}}\frac{\partial \big(\mathit{\Phi}_0^{h_1}-\mathit{\Phi}_0\big)}{\partial {x_j}}.
\end{aligned}
\end{equation}
Furthermore, multiplying both sides of equality (5.14) by $(T_0^{h_1}-T_0)$ and integrating over $\Omega$, it can derive that
\begin{equation}
\begin{aligned}
&\frac{1}{2}\frac{\partial}{\partial t}\Big[\int_{\Omega}\widehat S^{h_1}(T_{0}^{h_1})\big(T_0^{h_1}-T_0\big)^2d\Omega\Big]\!\!+\!\!\int_{\Omega}{{\widehat k_{ij}^{h_1}}(T_{0}^{h_1})\frac{{\partial \big(T_0^{h_1}-T_0\big)}}{{\partial {x_j}}}}\frac{\partial \big(T_0^{h_1}-T_0\big) }{{\partial {x_i}}}d\Omega\\
&=\frac{1}{2}\int_{\Omega}\frac{\partial \widehat S(T_{0}^{h_1})}{\partial t}\big(T_0^{h_1}-T_0\big)^2d\Omega+\int_{\Omega}\big[\widehat S(T_{0})-\widehat S^{h_1}(T_{0})\big]\frac{{\partial T_0}}{{\partial t}}\big(T_0^{h_1}-T_0\big)d\Omega\\
&+\int_{\Omega}\big[\widehat S^{h_1}(T_{0})-\widehat S^{h_1}(T_{0}^{h_1})\big]\frac{{\partial T_0}}{{\partial t}}\big(T_0^{h_1}-T_0\big)d\Omega\\
&+\int_{\Omega}{{\big(\widehat k_{ij}^{h_1}(T_{0}^{h_1})-\widehat k_{ij}^{h_1}(T_{0})\big)}\frac{{\partial T_0}}{{\partial {x_j}}}}\frac{\partial \big(T_0^{h_1}-T_0\big)}{{\partial {x_i}}}d\Omega\\
&+\int_{\Omega}{{\big(\widehat k_{ij}^{h_1}(T_{0})-\widehat k_{ij}(T_{0})\big)}\frac{{\partial T_0}}{{\partial {x_j}}}}\frac{\partial \big(T_0^{h_1}-T_0\big)}{{\partial {x_i}}}d\Omega\\
&+\int_{\Omega} {\widehat \lambda_{ij}^{h_1}}(T_{0}^{h_1})\frac{\partial \big(\mathit{\Phi}_0^{h_1}-\mathit{\Phi}_0\big)}{\partial {x_i}}\frac{\partial \mathit{\Phi}_{0}^{h_1}}{\partial {x_j}}\big(T_0^{h_1}-T_0\big)d\Omega\\
&+ \int_{\Omega}\big[{\widehat \lambda_{ij}^{h_1}}(T_{0}^{h_1})-{\widehat \lambda _{ij}^{h_1}}(T_{0})\big]\frac{\partial \mathit{\Phi}_0}{\partial {x_i}}\frac{\partial \mathit{\Phi}_{0}^{h_1}}{\partial {x_j}}\big(T_0^{h_1}-T_0\big)d\Omega\\
&+\int_{\Omega}\big[{\widehat \lambda_{ij}^{h_1}}(T_{0})-{\widehat \lambda_{ij}}(T_{0})\big]\frac{\partial \mathit{\Phi}_0}{\partial {x_i}}\frac{\partial \mathit{\Phi}_{0}^{h_1}}{\partial {x_j}}\big(T_0^{h_1}-T_0\big)d\Omega\\
&+ \int_{\Omega}{\widehat \lambda_{ij}}(T_{0})\frac{\partial \mathit{\Phi}_0}{\partial {x_i}}\frac{\partial \big(\mathit{\Phi}_0^{h_1}-\mathit{\Phi}_0\big)}{\partial {x_j}}\big(T_0^{h_1}-T_0\big)d\Omega.
\end{aligned}
\end{equation}
Subsequently, integrating both sides of (5.15) from $0$ to $t$ $(0<t\leq \mathcal )$, recalling inequalities (5.2)-(5.6) and (5.13), and employing the Cauchy-Schwarz inequality and Young inequality, we can obtain the following inequality from equality (5.15) if $|\widehat S^{h_1}(T_{0})-\widehat S^{h_1}(T_{0}^{h_1})|\leq C|T_{0}-T_{0}^{h_1}|$ and $|\widehat k_{ij}^{h_1}(T_{0})-\widehat k_{ij}^{h_1}(T_{0}^{h_1})|\leq C|T_{0}-T_{0}^{h_1}|$.
\begin{equation}
\begin{aligned}
&C\big\|T_0^{h_1}-T_0\big\|_{L^2(\Omega)}^2+C\mathlarger{\int}_0^t\big\|T_0^{h_1}-T_0\big\|_{H_0^1(\Omega)}^2d\tau\\
&\leq \mathlarger{\int}_0^t C\big\|T_0^{h_1}-T_0\big\|_{L^2(\Omega)}^2d\tau+\mathlarger{\int}_0^t\frac{C}{2}h_1^4d\tau+\mathlarger{\int}_0^t\frac{C}{2}\big\|T_0^{h_1}-T_0\big\|_{L^2(\Omega)}^2d\tau\\
&+\mathlarger{\int}_0^tC\big\|T_0^{h_1}-T_0\big\|_{L^2(\Omega)}^2d\tau+\mathlarger{\int}_0^t\frac{C}{2\lambda}\big\|T_0^{h_1}-T_0\big\|_{L^2(\Omega)}^2d\tau\\
&+\mathlarger{\int}_0^t\frac{\lambda}{2}\big\|T_0^{h_1}-T_0\big\|_{H_0^1(\Omega)}^2d\tau+\mathlarger{\int}_0^t\frac{C}{2\lambda}h_1^4d\tau+\mathlarger{\int}_0^t\frac{\lambda}{2}\big\|T_0^{h_1}-T_0\big\|_{H_0^1(\Omega)}^2d\tau\\
&+\mathlarger{\int}_0^tC\big\|\mathit{\Phi}_0^{h_1}-\mathit{\Phi}_0\big\|_{H_0^1(\Omega)}^2d\tau+\mathlarger{\int}_0^tC\big\|T_0^{h_1}-T_0\big\|_{L^2(\Omega)}^2d\tau\\
&+\mathlarger{\int}_0^tC\big\|T_0^{h_1}-T_0\big\|_{L^2(\Omega)}^2d\tau+\mathlarger{\int}_0^tC\big\|T_0^{h_1}-T_0\big\|_{L^2(\Omega)}^2d\tau\\
&+\mathlarger{\int}_0^t\frac{C}{2}h_1^4d\tau+\mathlarger{\int}_0^t\frac{C}{2}\big\|T_0^{h_1}-T_0\big\|_{L^2(\Omega)}^2d\tau+\mathlarger{\int}_0^tC\big\|\mathit{\Phi}_0^{h_1}-\mathit{\Phi}_0\big\|_{H_0^1(\Omega)}^2d\tau\\
&+\mathlarger{\int}_0^tC\big\|T_0^{h_1}-T_0\big\|_{L^2(\Omega)}^2d\tau\\
&\leq \mathlarger{\int}_0^tCh_1^4d\tau+\mathlarger{\int}_0^tC\big\|T_0^{h_1}-T_0\big\|_{L^2(\Omega)}^2d\tau\\
&+\mathlarger{\int}_0^t{\lambda}\big\|T_0^{h_1}-T_0\big\|_{H_0^1(\Omega)}^2d\tau+\mathlarger{\int}_0^tC\big\|\mathit{\Phi}_0^{h_1}-\mathit{\Phi}_0\big\|_{H_0^1(\Omega)}^2d\tau\\
&\leq \mathlarger{\int}_0^tCh_1^4d\tau+\mathlarger{\int}_0^tC\big\|T_0^{h_1}-T_0\big\|_{L^2(\Omega)}^2d\tau+\mathlarger{\int}_0^t{\lambda}\big\|T_0^{h_1}-T_0\big\|_{H_0^1(\Omega)}^2d\tau.
\end{aligned}
\end{equation}
By choosing a sufficiently small $\lambda$ and defining $\Upsilon_1(t)=\big\|T_0^{h_1}-T_0\big\|_{L^2(\Omega)}^2+\mathlarger{\int}_0^t\big\|T_0^{h_1}-T_0\big\|_{H_0^1(\Omega)}^2d\tau$, then we can derive $\Upsilon_1(t)\leq Ch_1^4+C\mathlarger{\int}_0^t\Upsilon_1(\tau)d\tau$ from (5.16). Consequently, by applying Gronwall inequality and leveraging the arbitrariness of time variable $t$, there holds the inequality $\big\|{T_0^{h_1}}-{T_0}\big\|_{L^\infty(0,\mathcal T;L^2(\Omega))}+
\big\|{T_0^{h_1}}-{T_0}\big\|_{L^2(0,\mathcal T;H^1_0(\Omega))}\leq Ch_1^2$. Additionally, the inequality, $\big\|{\mathit{\Phi}_0^{h_1}}-{\mathit{\Phi}_0}\big\|_{L^\infty(0,T;H_0^1(\Omega))}\leq Ch_1^2$, holds true.

Afterwards, subtracting the mechanical equations in (2.24) from the corresponding mechanical equations in (5.10) respectively, we can derive
\begin{equation}
\begin{aligned}
&\widehat \rho(T_0^{h_1})\frac{{\partial^2 \big(U_{i0}^{h_1}}-U_{i0}\big)}{{\partial t^2}}-\frac{\partial }{{\partial {x_j}}}\Big[{{\widehat c_{ijkl}^{h_1}}(T_0^{h_1})\frac{{\partial \big(U_{k0}^{h_1}}-U_{k0}\big)}{{\partial {x_l}}}}\Big]\\
&=\big[\widehat \rho(T_0)-\widehat \rho(T_0^{h_1})\big]\frac{\partial^2 U_{i0}}{{\partial t^2}}+\frac{\partial }{{\partial {x_j}}}\Big[ {{\big(\widehat c_{ijkl}^{h_1}(T_{0}^{h_1})-\widehat c_{ijkl}^{h_1}(T_{0})\big)}\frac{{\partial U_{k0}}}{{\partial {x_l}}}}\Big]\\
&+\frac{\partial }{{\partial {x_j}}}\Big[ {{\big(\widehat c_{ijkl}^{h_1}(T_{0})-\widehat c_{ijkl}(T_{0})\big)}\frac{{\partial U_{k0}}}{{\partial {x_l}}}}\Big]+\frac{\partial }{{\partial {x_j}}}\Big[\big(\widehat {\beta}_{ij}(T_0)-\widehat {\beta}_{ij}^{h_1}(T_0)\big)(T_0-\widetilde T)\Big]\\
&+\frac{\partial }{{\partial {x_j}}}\Big[\big(\widehat {\beta}_{ij}^{h_1}(T_0)-\widehat {\beta}_{ij}^{h_1}(T_0^{h_1})\big)(T_0-\widetilde T)\Big]+\frac{\partial }{{\partial {x_j}}}\Big[\widehat {\beta}_{ij}^{h_1}(T_0^{h_1})(T_0-T_0^{h_1})\Big].
\end{aligned}
\end{equation}
Furthermore, implementing the calculation $\displaystyle\int_{\Omega}(5.17)\times \frac{{\partial \big(U_{i0}^{h_1}}-U_{i0}\big)}{{\partial t}}d\Omega$, suffice it to have the following
equality
\begin{equation}
\begin{aligned}
&\frac{1}{2}\frac{\partial}{\partial t}\Big[\int_{\Omega}\widehat \rho(T_0^{h_1})\big[\frac{\partial \big(U_{i0}^{h_1}-U_{i0}\big)}{\partial t}\big]^2d\Omega\\
&+\int_{\Omega}\widehat c_{ijkl}^{h_1}(T_0^{h_1})\frac{\partial \big(U_{k0}^{h_1}-U_{k0}\big)}{\partial {x_l}}\frac{\partial \big(U_{i0}^{h_1}-U_{i0}\big)}{\partial {x_j}}d\Omega\Big]\\
&=\frac{1}{2}\int_{\Omega}\frac{\partial\widehat \rho(T_0^{h_1})}{\partial t}\big[\frac{\partial \big(U_{i0}^{h_1}-U_{i0}\big)}{\partial t}\big]^2d\Omega\\
&+\frac{1}{2}\int_{\Omega}\frac{\partial\widehat c_{ijkl}^{h_1}(T_0^{h_1})}{\partial t}\frac{\partial \big(U_{k0}^{h_1}-U_{k0}\big)}{\partial {x_l}}\frac{\partial \big(U_{i0}^{h_1}-U_{i0}\big)}{\partial {x_j}}d\Omega\\
&+\int_{\Omega}\big[\widehat \rho(T_0)-\widehat \rho(T_0^{h_1})\big]\frac{\partial^2 U_{i0}}{{\partial t^2}}\frac{\partial \big(U_{i0}^{h_1}-U_{i0}\big)}{\partial {t}}d\Omega\\
&-\int_{\Omega} {{\big(\widehat c_{ijkl}^{h_1}(T_{0}^{h_1})-\widehat c_{ijkl}^{h_1}(T_{0})\big)}\frac{{\partial U_{k0}}}{{\partial {x_l}}}}\frac{\partial }{{\partial {x_j}}}\Big(\frac{\partial \big(U_{i0}^{h_1}-U_{i0}\big)}{\partial {t}}\Big)d\Omega\\
&-\int_{\Omega}{{\big(\widehat c_{ijkl}^{h_1}(T_{0})-\widehat c_{ijkl}(T_{0})\big)}\frac{{\partial U_{k0}}}{{\partial {x_l}}}}\frac{\partial }{{\partial {x_j}}}\Big(\frac{\partial \big(U_{i0}^{h_1}-U_{i0}\big)}{\partial {t}}\Big)d\Omega\\
&+\int_{\Omega}\big(\widehat {\beta}_{ij}(T_0)-\widehat {\beta}_{ij}^{h_1}(T_0)\big)(T_0-\widetilde T)\frac{\partial }{{\partial {x_j}}}\Big(\frac{\partial \big(U_{i0}^{h_1}-U_{i0}\big)}{\partial {t}}\Big)d\Omega\\
&+\int_{\Omega}\big(\widehat {\beta}_{ij}^{h_1}(T_0)-\widehat {\beta}_{ij}^{h_1}(T_0^{h_1})\big)(T_0-\widetilde T)\frac{\partial }{{\partial {x_j}}}\Big(\frac{\partial \big(U_{i0}^{h_1}-U_{i0}\big)}{\partial {t}}\Big)d\Omega\\
&+\int_{\Omega}\widehat {\beta}_{ij}^{h_1}(T_0^{h_1})(T_0-T_0^{h_1})\frac{\partial }{{\partial {x_j}}}\Big(\frac{\partial \big(U_{i0}^{h_1}-U_{i0}\big)}{\partial {t}}\Big)d\Omega.
\end{aligned}
\end{equation}
Subsequently, integrating both sides of (5.18) from $0$ to $t$ $(0<t\leq \mathcal T)$, recalling inequalities (5.2)-(5.6), and employing the Cauchy-Schwarz inequality and Young inequality, we can obtain the following inequality from equality (5.18) if $|\widehat \rho(T_{0})-\widehat \rho(T_{0}^{h_1})|\leq C|T_{0}-T_{0}^{h_1}|$ , $|\widehat \beta_{ij}^{h_1}(T_{0})-\widehat \beta_{ij}^{h_1}(T_{0}^{h_1})|\leq C|T_{0}-T_{0}^{h_1}|$, $|\widehat c_{ijkl}^{h_1}(T_{0})-\widehat c_{ijkl}^{h_1}(T_{0}^{h_1})|\leq C|T_{0}-T_{0}^{h_1}|$, $\displaystyle|\frac{\partial\widehat \beta_{ij}^{h_1}(T_{0})}{\partial t}-\frac{\partial\widehat \beta_{ij}^{h_1}(T_{0}^{h_1})}{\partial t}|\leq C|T_{0}-T_{0}^{h_1}|$, $\displaystyle|\frac{\partial\widehat c_{ijkl}^{h_1}(T_{0})}{\partial t}-\frac{\partial\widehat c_{ijkl}^{h_1}(T_{0}^{h_1})}{\partial t}|\leq C|T_{0}-T_{0}^{h_1}|$ and $\displaystyle|\frac{\partial T_{0}}{\partial t}-\frac{\partial T_{0}^{h_1}}{\partial t}|\leq C|T_{0}-T_{0}^{h_1}|$.
\begin{equation}
\begin{aligned}
&C\Big\|\frac{\partial{\bm{U}_0^{h_1}}}{\partial t}-\frac{\partial{\bm{U}_0}}{\partial t}\Big\|_{(L^2(\Omega))^n}^2+
C\big\|{\bm{U}_0^{h_1}}-{\bm{U}_0}\big\|_{(H^1_0(\Omega))^n}^2\\
&\leq \mathlarger{\int}_0^tC\Big\|\frac{\partial{\bm{U}_0^{h_1}}}{\partial \tau}-\frac{\partial{\bm{U}_0}}{\partial \tau}\Big\|_{(L^2(\Omega))^n}^2d\tau+\mathlarger{\int}_0^tC\big\|{\bm{U}_0^{h_1}}-{\bm{U}_0}\big\|_{(H^1_0(\Omega))^n}^2d\tau\\
&+\mathlarger{\int}_0^tC\big\|T_0^{h_1}-T_0\big\|_{L^2(\Omega)}^2d\tau+\mathlarger{\int}_0^tC\Big\|\frac{\partial{\bm{U}_0^{h_1}}}{\partial \tau}-\frac{\partial{\bm{U}_0}}{\partial \tau}\Big\|_{(L^2(\Omega))^n}^2d\tau\\
&+\frac{C}{2\lambda}\big\|T_0^{h_1}-T_0\big\|_{L^2(\Omega)}^2+\frac{\lambda}{2}\big\|{\bm{U}_0^{h_1}}-{\bm{U}_0}\big\|_{(H^1_0(\Omega))^n}^2\\
&+\mathlarger{\int}_0^tC\big\|T_0^{h_1}-T_0\big\|_{L^2(\Omega)}^2d\tau+\mathlarger{\int}_0^tC\big\|{\bm{U}_0^{h_1}}-{\bm{U}_0}\big\|_{(H^1_0(\Omega))^n}^2d\tau\\
&+\mathlarger{\int}_0^tC\big\|T_0^{h_1}-T_0\big\|_{L^2(\Omega)}^2d\tau+\mathlarger{\int}_0^tC\big\|{\bm{U}_0^{h_1}}-{\bm{U}_0}\big\|_{(H^1_0(\Omega))^n}^2d\tau\\
&+\frac{C}{2\lambda}h_1^4+\frac{\lambda}{2}\big\|{\bm{U}_0^{h_1}}-{\bm{U}_0}\big\|_{(H^1_0(\Omega))^n}^2\\
&+\mathlarger{\int}_0^tCh_1^4d\tau+\mathlarger{\int}_0^tC\big\|{\bm{U}_0^{h_1}}-{\bm{U}_0}\big\|_{(H^1_0(\Omega))^n}^2d\tau\\
&+\mathlarger{\int}_0^tCh_1^4d\tau+\mathlarger{\int}_0^tC\big\|{\bm{U}_0^{h_1}}-{\bm{U}_0}\big\|_{(H^1_0(\Omega))^n}^2d\tau\\
&+\frac{C}{2\lambda}\big\|T_0^{h_1}-T_0\big\|_{L^2(\Omega)}^2+\frac{\lambda}{2}\big\|{\bm{U}_0^{h_1}}-{\bm{U}_0}\big\|_{(H^1_0(\Omega))^n}^2\\
&+\mathlarger{\int}_0^tC\big\|T_0^{h_1}-T_0\big\|_{L^2(\Omega)}^2d\tau+\mathlarger{\int}_0^tC\big\|{\bm{U}_0^{h_1}}-{\bm{U}_0}\big\|_{(H^1_0(\Omega))^n}^2d\tau\\
&+\mathlarger{\int}_0^tC\big\|T_0^{h_1}-T_0\big\|_{L^2(\Omega)}^2d\tau+\mathlarger{\int}_0^tC\big\|{\bm{U}_0^{h_1}}-{\bm{U}_0}\big\|_{(H^1_0(\Omega))^n}^2d\tau\\
&+\frac{C}{2\lambda}h_1^4+\frac{\lambda}{2}\big\|{\bm{U}_0^{h_1}}-{\bm{U}_0}\big\|_{(H^1_0(\Omega))^n}^2\\
&+\mathlarger{\int}_0^tCh_1^4d\tau+\mathlarger{\int}_0^tC\big\|{\bm{U}_0^{h_1}}-{\bm{U}_0}\big\|_{(H^1_0(\Omega))^n}^2d\tau\\
&+\mathlarger{\int}_0^tCh_1^4d\tau+\mathlarger{\int}_0^tC\big\|{\bm{U}_0^{h_1}}-{\bm{U}_0}\big\|_{(H^1_0(\Omega))^n}^2d\tau\\
&+\frac{C}{2\lambda}\big\|T_0^{h_1}-T_0\big\|_{L^2(\Omega)}^2+\frac{\lambda}{2}\big\|{\bm{U}_0^{h_1}}-{\bm{U}_0}\big\|_{(H^1_0(\Omega))^n}^2\\
&+\mathlarger{\int}_0^tC\big\|T_0^{h_1}-T_0\big\|_{L^2(\Omega)}^2d\tau+\mathlarger{\int}_0^tC\big\|{\bm{U}_0^{h_1}}-{\bm{U}_0}\big\|_{(H^1_0(\Omega))^n}^2d\tau\\
&+\mathlarger{\int}_0^tC\big\|T_0^{h_1}-T_0\big\|_{L^2(\Omega)}^2d\tau+\mathlarger{\int}_0^tC\big\|{\bm{U}_0^{h_1}}-{\bm{U}_0}\big\|_{(H^1_0(\Omega))^n}^2d\tau\\
&\leq \mathlarger{\int}_0^tCh_1^4d\tau\!\!+\!\!\mathlarger{\int}_0^tC\Big\|\frac{\partial{\bm{U}_0^{h_1}}}{\partial \tau}-\frac{\partial{\bm{U}_0}}{\partial \tau}\Big\|_{(L^2(\Omega))^n}^2d\tau\!\!+\!\!\mathlarger{\int}_0^tC\big\|{\bm{U}_0^{h_1}}-{\bm{U}_0}\big\|_{(H^1_0(\Omega))^n}^2d\tau.
\end{aligned}
\end{equation}
When choosing a sufficiently small $\lambda$ and setting $\displaystyle\Upsilon_2(t)=\Big\|\frac{\partial{\bm{U}_0^{h_1}}}{\partial t}-\frac{\partial{\bm{U}_0}}{\partial t}\Big\|_{(L^2(\Omega))^n}^2+
\big\|{\bm{U}_0^{h_1}}-{\bm{U}_0}\big\|_{(H^1_0(\Omega))^n}^2$, then we can derive $\Upsilon_2(t)\leq Ch_1^4+C\mathlarger{\int}_0^t\Upsilon_2(\tau)d\tau$ from (5.19). Consequently, by taking advantage of Gronwall inequality and the arbitrariness of time variable $t$, there holds the inequality $\displaystyle\Big\|\frac{\partial{\bm{U}_0^{h_1}}}{\partial t}-\frac{\partial{\bm{U}_0}}{\partial t}\Big\|_{L^\infty(0,\mathcal T;(L^2(\Omega))^n)}+\big\|{\bm{U}_0^{h_1}}-{\bm{U}_0}\big\|_{L^\infty(0,\mathcal T;(H^1_0(\Omega))^n)}\leq Ch_1^2$. In summary, combining the inequalities for electric potential, temperature and displacement fields, we can easily obtain inequality (5.10).
\begin{lemma}
Denote $\mathit{\Phi}_{0}^{h_1,h_0}$, $T_{0}^{h_1,h_0}$ and $\bm{U}_{0}^{h_1,h_0}$ be the finite element solutions of the modified homogenized equations (5.9) by mixed FDM-FEM scheme proposed in reference \cite{R11}, the following estimate holds
\begin{equation}
{\begin{aligned}
&\mathop{\max}_{1\leq n\leq N}\big\|{\mathit{\Phi}_0^{h_1,h_0}}(\bm{x},t_n)-{\mathit{\Phi}_0^{h_1}(\bm{x},t_n)}\big\|_{L^2(\Omega)}\\
&+\mathop{\max}_{1\leq n\leq N}\big\|{T_0^{h_1,h_0}}(\bm{x},t_{n})-{T_0^{h_1}(\bm{x},t_{n})}\big\|_{L^{2}(\Omega)}\\
&+\mathop{\max}_{1\leq n\leq N}\big\|\frac{\partial{\bm{U}_0^{h_1,h_0}(\bm{x},t_{n})}}{\partial t}-\frac{\partial{\bm{U}_0^{h_1}}(\bm{x},t_{n})}{\partial t}\big\|_{(L^{2}(\Omega))^n}\leq C\Delta t+Ch_0^2,
\end{aligned}}
\end{equation}
where $C$ is a constant irrespective of $h_1$ and $\Delta t$.
\end{lemma}
$\mathbf{Proof:}$ As shown in the modified homogenized equations (5.9), it satisfy the conditions of error estimates in reference \cite{R11}. Hence, employing the same proof technique as reference \cite{R11}, the estimate (5.20) can be derived.
\begin{theorem}
Let $\mathit{\Phi}_{0}^{h_1,h_0}$, $T_{0}^{h_1,h_0}$ and $\bm{U}_{0}^{h_1,h_0}$ be the finite element solutions of the modified homogenized equations (5.9), and $\mathit{\Phi}_{0}$, $T_{0}$ and $\bm{U}_{0}$ be the exact solutions of the homogenized equations (2.24), then the following estimate holds
\begin{equation}
{\begin{aligned}
&\mathop{\max}_{1\leq n\leq N}\big\|{\mathit{\Phi}_0^{h_1,h_0}}(\bm{x},t_n)-{\mathit{\Phi}_0(\bm{x},t_n)}\big\|_{L^2(\Omega)}\\
&+\mathop{\max}_{1\leq n\leq N}\big\|{T_0^{h_1,h_0}}(\bm{x},t_{n})-{T_0(\bm{x},t_{n})}\big\|_{L^{2}(\Omega)}\\
&+\mathop{\max}_{1\leq n\leq N}\big\|\frac{\partial{\bm{U}_0^{h_1,h_0}(\bm{x},t_{n})}}{\partial t}-\frac{\partial{\bm{U}_0}(\bm{x},t_{n})}{\partial t}\big\|_{(L^{2}(\Omega))^n}\leq C\Delta t+Ch_0^2+Ch_1^2,
\end{aligned}}
\end{equation}
where $C$ is a positive constant independent of $h_1$, $h_0$ and $\Delta t$.
\end{theorem}
$\mathbf{Proof:}$ Firstly, employing the triangle inequality, there exists a inequality such that
\begin{equation}
{\begin{aligned}
&\mathop{\max}_{1\leq n\leq N}\big\|{\mathit{\Phi}_0^{h_1,h_0}}(\bm{x},t_n)-{\mathit{\Phi}_0(\bm{x},t_n)}\big\|_{L^2(\Omega)}\\
&+\mathop{\max}_{1\leq n\leq N}\big\|{T_0^{h_1,h_0}}(\bm{x},t_{n})-{T_0(\bm{x},t_{n})}\big\|_{L^{2}(\Omega)}\\
&+\mathop{\max}_{1\leq n\leq N}\big\|\frac{\partial{\bm{U}_0^{h_1,h_0}(\bm{x},t_{n})}}{\partial t}-\frac{\partial{\bm{U}_0}(\bm{x},t_{n})}{\partial t}\big\|_{(L^{2}(\Omega))^n}\\
&\leq \mathop{\max}_{1\leq n\leq N}\big\|{\mathit{\Phi}_0^{h_1,h_0}}(\bm{x},t_n)-{\mathit{\Phi}_0^{h_1}(\bm{x},t_n)}\big\|_{L^2(\Omega)}\\
&+\mathop{\max}_{1\leq n\leq N}\big\|{T_0^{h_1,h_0}}(\bm{x},t_{n})-{T_0^{h_1}(\bm{x},t_{n})}\big\|_{L^{2}(\Omega)}\\
&+\mathop{\max}_{1\leq n\leq N}\big\|\frac{\partial{\bm{U}_0^{h_1,h_0}(\bm{x},t_{n})}}{\partial t}-\frac{\partial{\bm{U}_0^{h_1}}(\bm{x},t_{n})}{\partial t}\big\|_{(L^{2}(\Omega))^n}\\
&+\mathop{\max}_{1\leq n\leq N}\big\|{\mathit{\Phi}_0^{h_1}}(\bm{x},t_n)-{\mathit{\Phi}_0(\bm{x},t_n)}\big\|_{L^2(\Omega)}\\
&+\mathop{\max}_{1\leq n\leq N}\big\|{T_0^{h_1}}(\bm{x},t_{n})-{T_0(\bm{x},t_{n})}\big\|_{L^{2}(\Omega)}\\
&+\mathop{\max}_{1\leq n\leq N}\big\|\frac{\partial{\bm{U}_0^{h_1}(\bm{x},t_{n})}}{\partial t}-\frac{\partial{\bm{U}_0}(\bm{x},t_{n})}{\partial t}\big\|_{(L^{2}(\Omega))^n}.
\end{aligned}}
\end{equation}
Moreover, with the help of the embedding theorem, suffice it to have the following inequality
\begin{equation}
\begin{aligned}
&\mathop{\max}_{1\leq n\leq N}\big\|{\mathit{\Phi}_0^{h_1,h_0}}(\bm{x},t_n)-{\mathit{\Phi}_0^{h_1}(\bm{x},t_n)}\big\|_{L^2(\Omega)}\\
&\leq
\big\|{\mathit{\Phi}_0^{h_1}}-{\mathit{\Phi}_0}\big\|_{L^\infty(0,\mathcal T;H_0^1(\Omega))}\leq Ch_1^2.
\end{aligned}
\end{equation}
Finally, substituting (5.10), (5.20) and (5.23) into (5.22), one can directly obtain the error estimate (5.21).

In summary, we remark that the above-mentioned theoretical analysis rigorously ensures the convergence of the proposed two-stage numerical algorithm in both microscopic and macroscopic computations.
\section{Numerical examples and results}
In this section, numerical experiments are conducted on a HPS desktop workstation equipped with an Intel(R) Xeon(R) Gold 6146 CPU (3.20 GHz) and internal memory (96.0 GB), and implemented based on Freefem++ software. In addition, since exact solutions for nonlinear multi-scale problems are unattainable, we use direct numerical simulation (DNS) solutions $T_{\text{DNS}}^\varepsilon(\bm{x},t)$, $\mathit{\Phi}_{{\text{DNS}}}^\varepsilon(\bm{x},t)$ and ${\bm{U}}_{{\text{DNS}}}^{{\varepsilon}}(\bm{x},t)$ to replace $T^{\varepsilon}(\bm{x},t)$, $\mathit{\Phi}^{{\varepsilon}}(\bm{x},t)$ and ${\bm{U}}^{{\varepsilon}}(\bm{x},t)$ for evaluating the proposed HOMS method. Furthermore, we define the following notations: $\text{T}_{L^2}^0(t)$, $\text{T}_{L^2}^1(t)$ and $\text{T}_{L^2}^2(t)$ denote the relative errors for homogenized solutions, LOMS solutions and HOMS solutions of temperature field in the $L^2$ norm, while $\text{T}_{H^1}^0(t)$, $\text{T}_{H^1}^1(t)$ and $\text{T}_{H^1}^2(t)$ denote the relative errors for those solutions in the $H^1$ semi-norm. Analogously, $\text{P}_{L^2}^0(t)$, $\text{P}_{L^2}^1(t)$ and $\text{P}_{L^2}^2(t)$ represent the relative errors for electric potential field in the $L^2$ norm, and $\text{P}_{H^1}^0(t)$, $\text{P}_{H^1}^1(t)$ and $\text{P}_{H^1}^2(t)$ represent those in the $H^1$ semi-norm. $\text{D}_{L^2}^0(t)$, $\text{D}_{L^2}^1(t)$ and $\text{D}_{L^2}^2(t)$ represent the relative errors for displacement field in the $L^2$ norm, and $\text{D}_{H^1}^0(t)$, $\text{D}_{H^1}^1(t)$ and $\text{D}_{H^1}^2(t)$ represent those in the $H^1$ semi-norm.
\subsection{Example 1: nonlinear thermo-electro-mechanical coupling simulation of 2D composite structure}
This example investigates the nonlinear thermo-electro-mechanical coupling behaviors of a 2D composite structure, which is composed of repeating unit cells with the characteristic periodic length $\varepsilon=1/10$. The whole domain $\displaystyle\Omega=(x_1,x_2)=[0,1]\times[0,1]$ and PUC $\mathbf{Y}$ are displayed in Fig.\hspace{1mm}1.
\begin{figure}[!htb]
\centering
\begin{minipage}[c]{0.32\textwidth}
  \centering
  \includegraphics[width=0.60\linewidth,totalheight=1.0in]{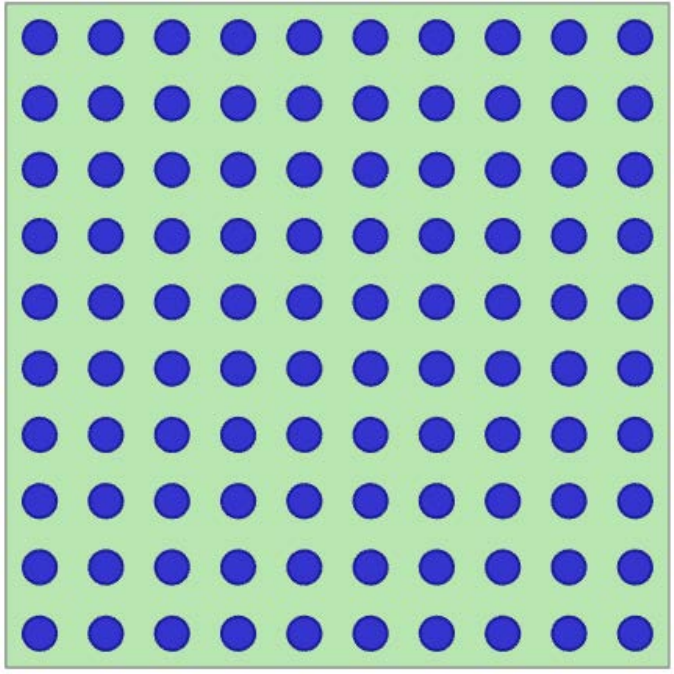} \\
  (a)
\end{minipage}
\begin{minipage}[c]{0.32\textwidth}
  \centering
  \includegraphics[width=0.98\linewidth,totalheight=1.0in]{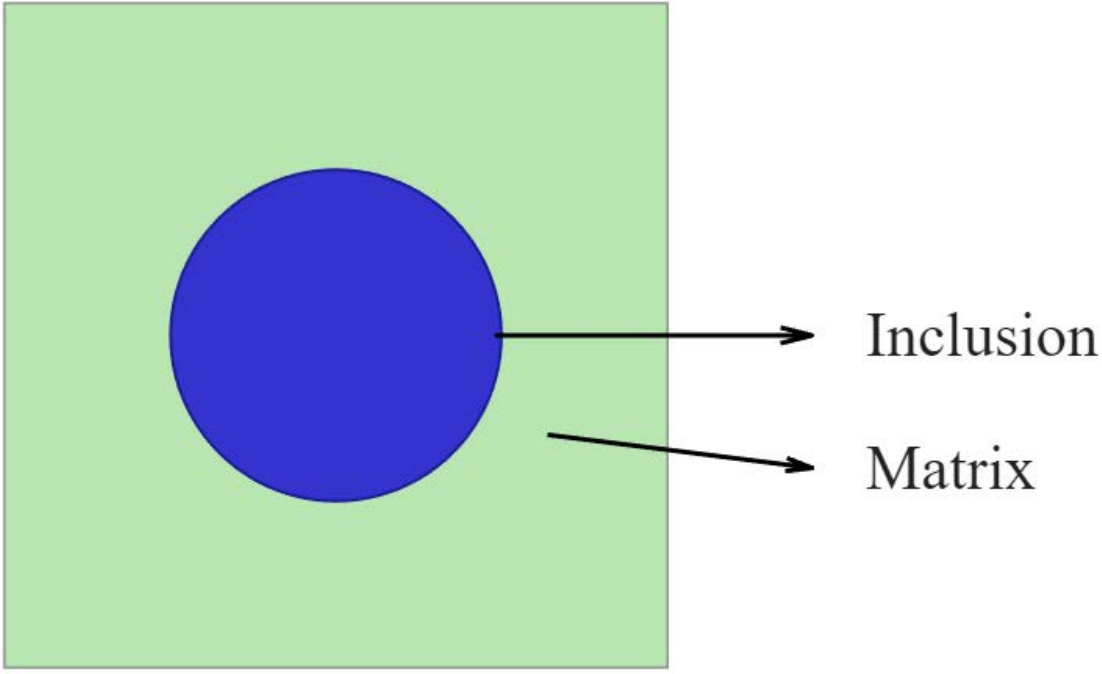} \\
  (b)
\end{minipage}
\begin{minipage}[c]{0.32\textwidth}
  \centering
  \includegraphics[width=0.60\linewidth,totalheight=1.0in]{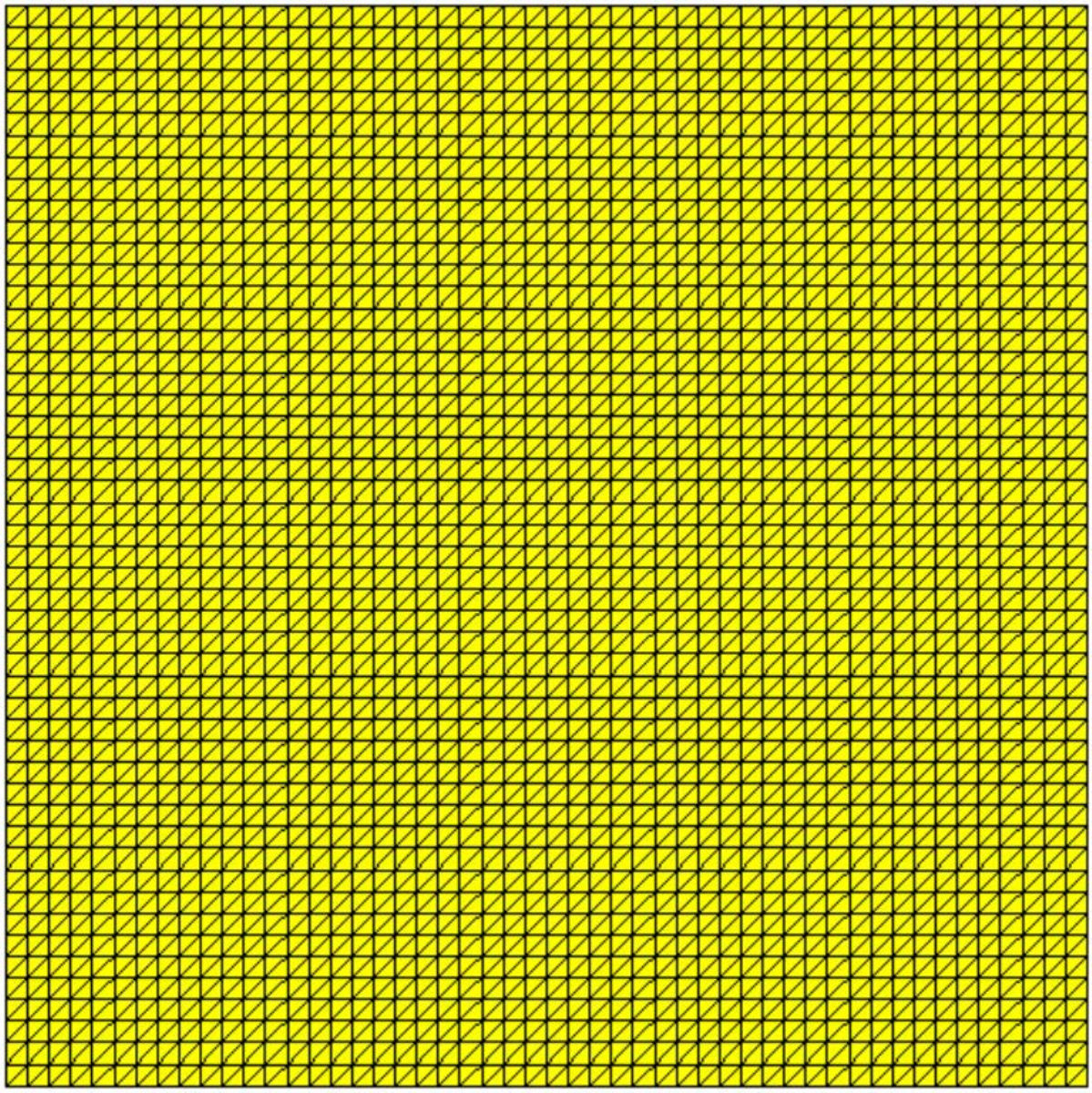} \\
  (c)
\end{minipage}
\caption{(a) The 2D composite structure $\Omega$; (b) PUC $\mathbf{Y}$; (c) homogenized structure $\Omega$.}
\end{figure}

Additionally, the material parameters of matrix and inclusion phases of this composite structure are detailed in Table 1.
\begin{table}[h]{\caption{Material parameters of composite structure ($T$ stands for the temperature value).}}
\centering
\begin{tabular}{cc}
\hline
Material parameters & Matrix/Inclusion \\
\hline
$\rho^\varepsilon$ & 0.008/0.002  \\
$c^\varepsilon$ & 562.5/750.0  \\
$k_{ij}^\varepsilon$ & 4.0+0.0004$T$/0.04+0.000004$T$  \\
$\lambda_{ij}^{\varepsilon}$ & 300.0-0.015$T$/0.075-0.00000325$T$ \\
$\beta_{ij}^{\varepsilon}$ & 3.0-0.0003$T$/7.5-0.00075$T$ \\
$E^{\varepsilon}$ & 3.5$\times 10^6$-3.5$\times 10^3T$/2.2$\times 10^6$-2.2$\times 10^3T$ \\
$\nu^{\varepsilon}$ & 0.25/0.20 \\
\hline
\end{tabular}
\end{table}

Moreover, the source items, boundary conditions and initial conditions in the multi-scale nonlinear problem (1.1) of this example are presented as below.
\begin{equation}
\begin{aligned}
&f_T(\bm{x},t)=20000,\;f_\mathit{\Phi}(\bm{x},t)=200,\;f_1(\bm{x},t)=5000,\;f_2(\bm{x},t)=5000,\;\;\text{in}\;\;\Omega,\\
&\widehat{T}(\bm{x},t)=300.0,\;\widehat{\mathit{\Phi}}(\bm{x},t)=0.0,\;\widehat{\bm{U}}(\bm{x},t)=(0.0,0.0,0.0),\;\;\text{in}\;\;\partial\Omega,\\
&\widetilde T=300.0,\;\widetilde{\bm{U}}^0({\bm{x}})=(0.0,0.0,0.0),\;\widetilde{\bm{U}}^1({\bm{x}})=(0.0,0.0,0.0),\;\;\text{in}\;\;\Omega.
\end{aligned}
\end{equation}

In this example, the 2D heterogeneous structure contains a multitude of microscopic unit cells. Resolving these microscale features via a direct finite element method (FEM) demands excessively fine meshes, leading to severe computational overhead. Consequently, separate finite element meshes are generated for original multi-scale equations, auxiliary cell equations and corresponding homogenized equations. Detailed mesh statistics and CPU times are summarized in Table 2.
\begin{table}[!htb]{\caption{Summary of computational cost in 2D composite structure.}}
\centering
\begin{tabular}{cccc}
\hline
 & Multi-scale eqs. & Cell eqs. & Homogenized eqs. \\
\hline
FEM elements & 70800 & 856 & 5000\\
FEM nodes    & 35761 & 469 & 2601\\
\hline
& DNS & off-line stage & on-line stage \\
\hline
Computational time & 13086.440s  & 41.129s & 8294.787s\\
\hline
\end{tabular}
\end{table}

For this example, 10 equidistant macroscopic interpolation temperatures are prescribed within a single unit cell. It is noteworthy that auxiliary cell problems are off-line computation before on-line multi-scale computation, which remain valid for different heterogeneous structures. The dynamic nonlinear thermo-electro-mechanical coupling responses of the 2D heterogeneous structure are simulated over the time interval $t\in[0,1]$. Setting the temporal step as $\Delta t = 0.001$, the macroscopic homogenized equations (2.24) and multi-scale nonlinear equations (1.1) are on-line simulated respectively. Following this, the ultimate simulation outcomes of temperature, electric potential and displacement fields at $t=1.0$ are depicted in Figs.\hspace{1mm}2-4, respectively. Besides, Fig.\hspace{1mm}5 is plotted to visualize the evolutionary relative errors of three physical fields tracked in the $L^2$ norm and $H^1$ semi-norm senses.
\begin{figure}[!htb]
\centering
\begin{minipage}[c]{0.24\textwidth}
  \centering
  \includegraphics[width=31mm]{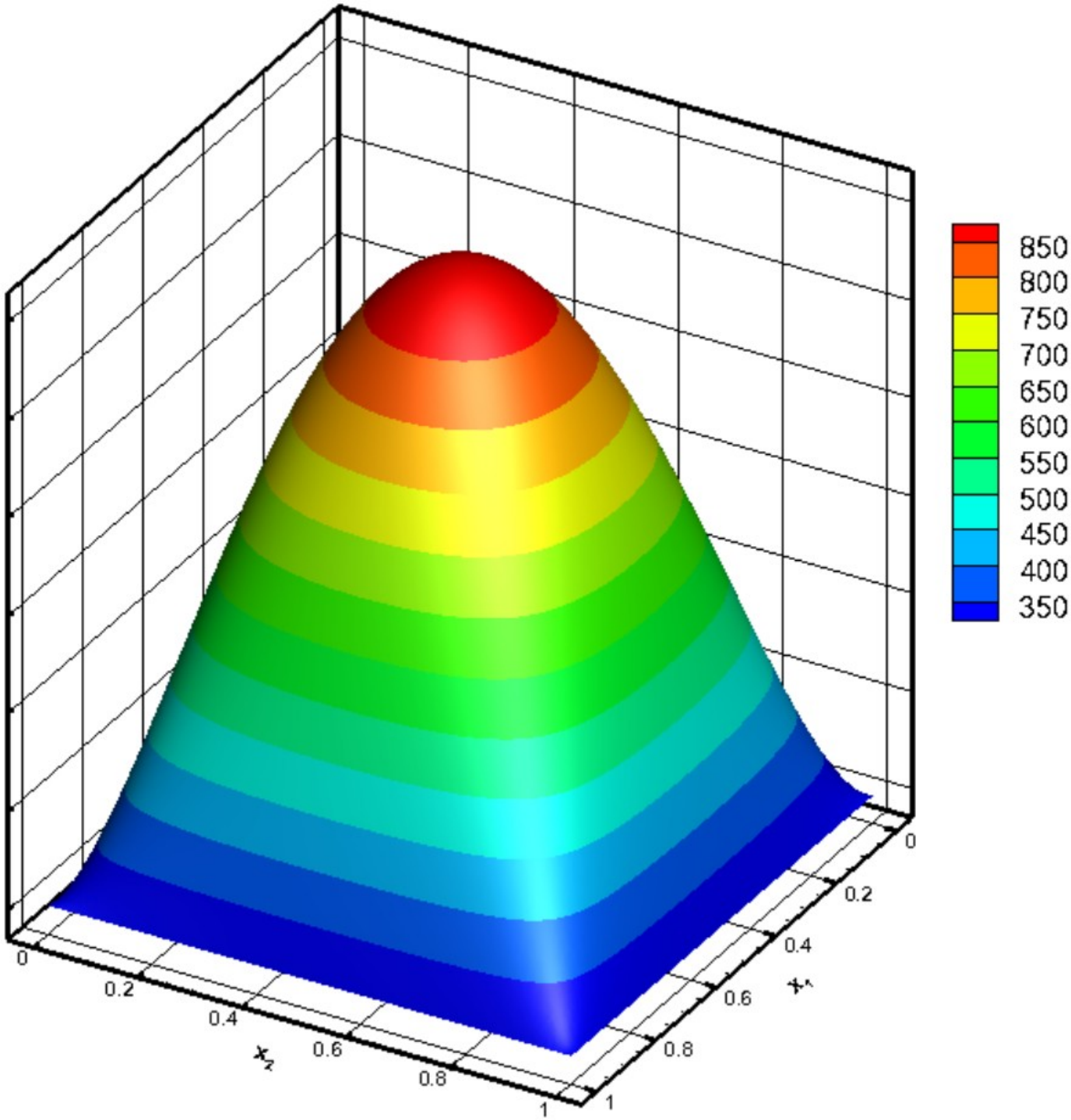}\\
  (a)
\end{minipage}
\begin{minipage}[c]{0.24\textwidth}
  \centering
  \includegraphics[width=31mm]{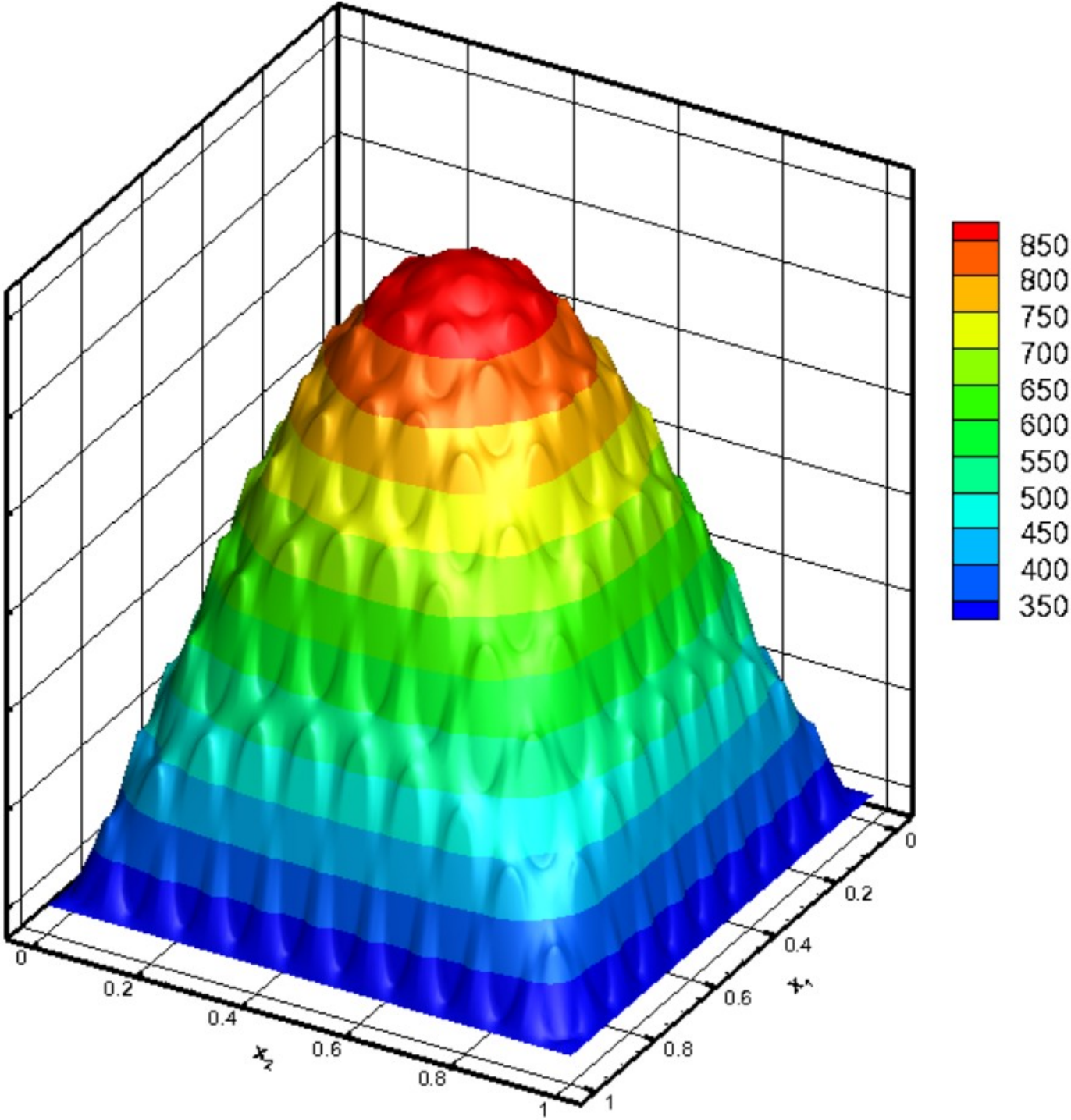}\\
  (b)
\end{minipage}
\begin{minipage}[c]{0.24\textwidth}
  \centering
  \includegraphics[width=31mm]{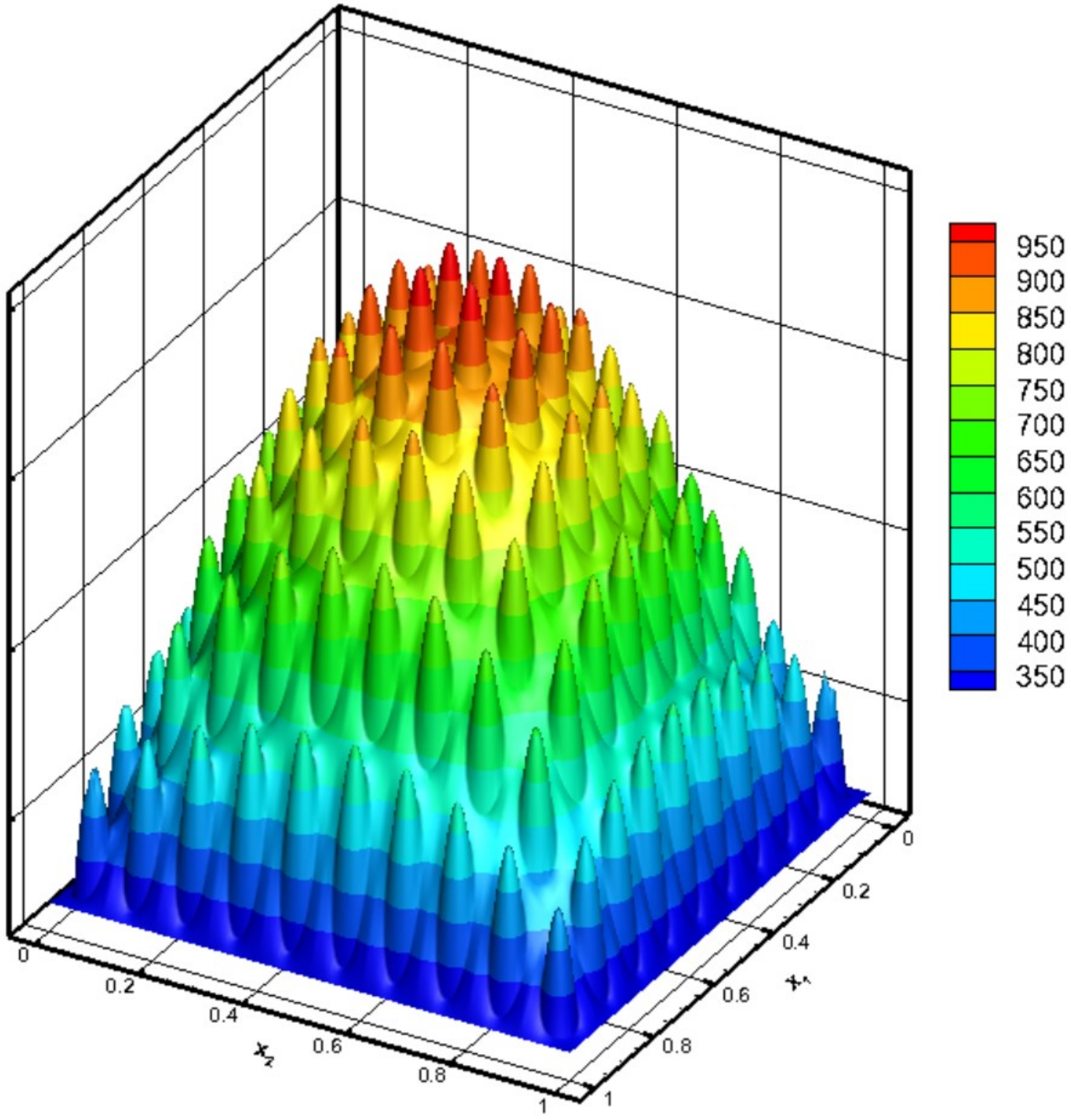}\\
  (c)
\end{minipage}
\begin{minipage}[c]{0.24\textwidth}
  \centering
  \includegraphics[width=31mm]{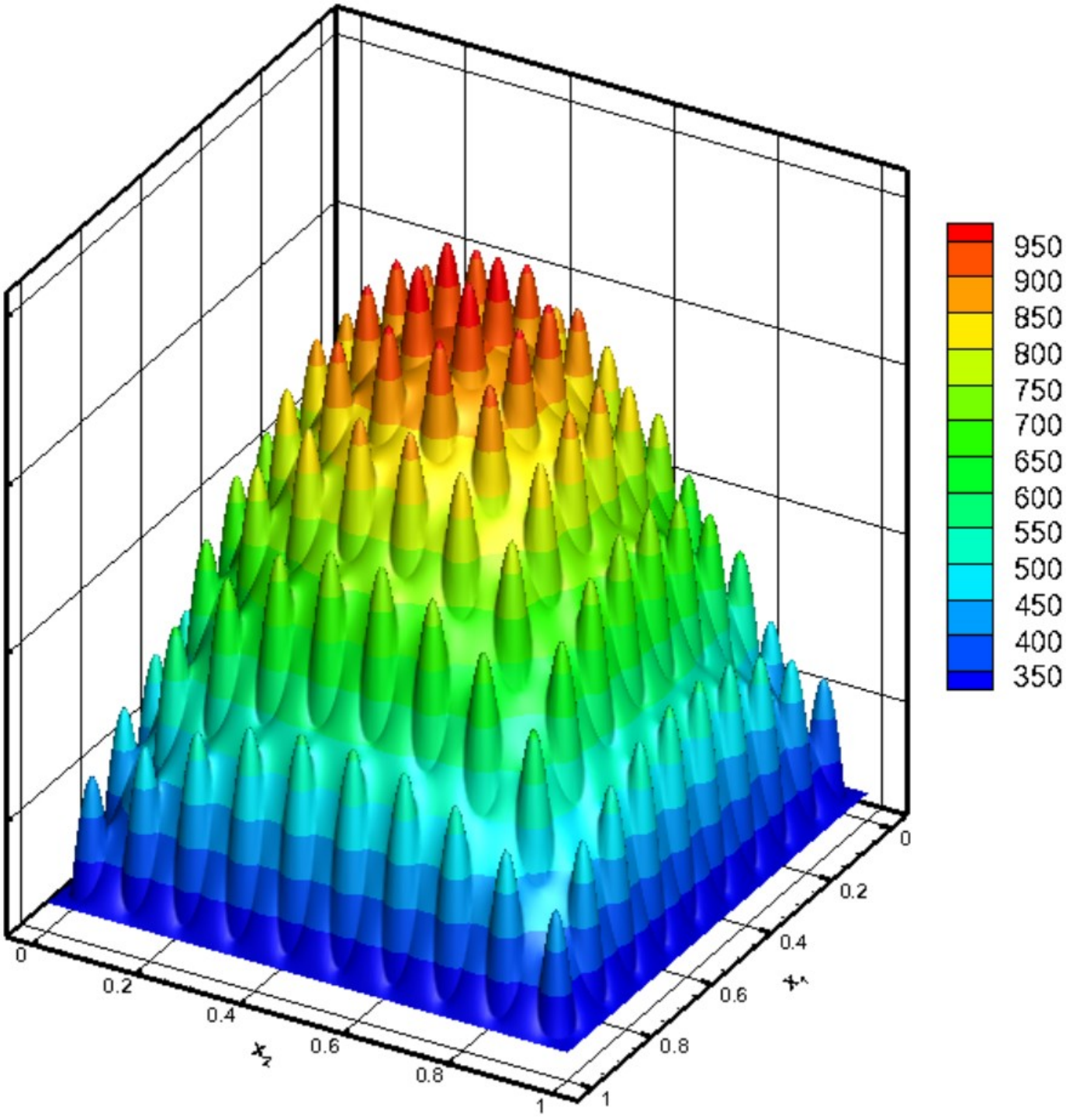}\\
  (d)
\end{minipage}
\caption{The temperature field at $t=1.0$: (a) $T_{0}$; (b) $T^{[1\varepsilon]}$; (c) $T^{[2\varepsilon]}$; (d) $T_{\rm{DNS}}^\varepsilon$.}
\end{figure}
\begin{figure}[!htb]
\centering
\begin{minipage}[c]{0.24\textwidth}
  \centering
  \includegraphics[width=30mm]{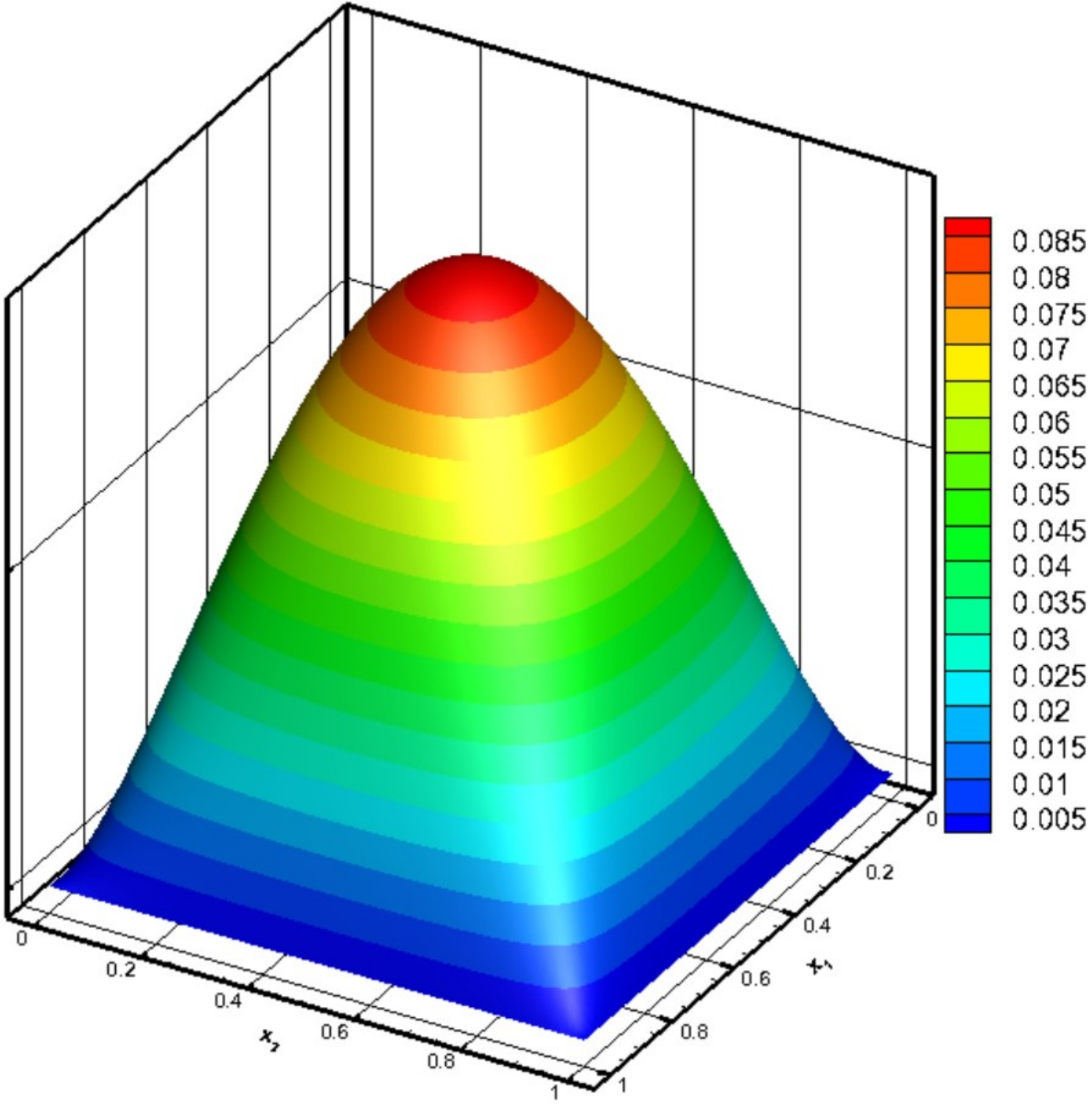}\\
  (a)
\end{minipage}
\begin{minipage}[c]{0.24\textwidth}
  \centering
  \includegraphics[width=30mm]{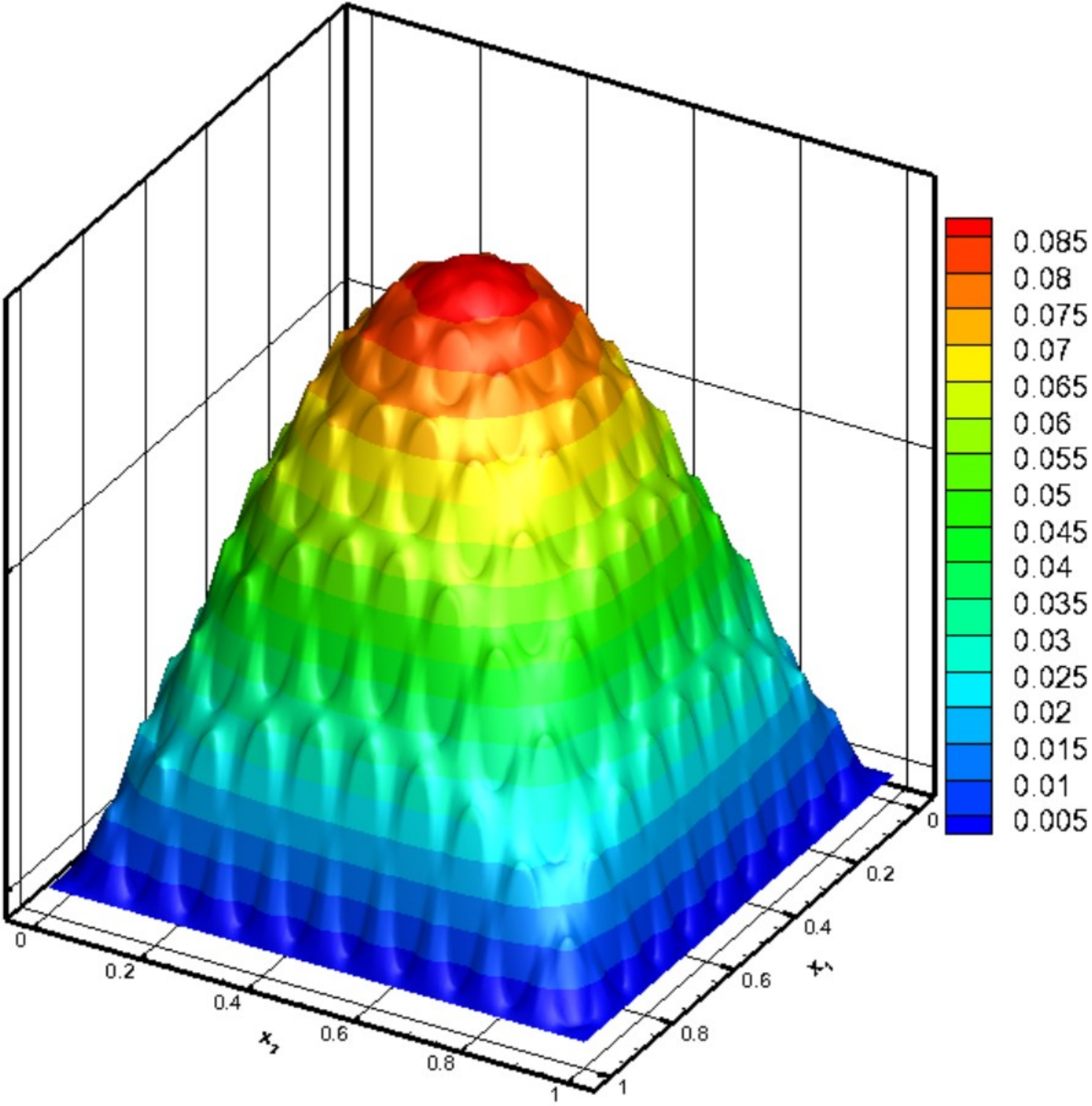}\\
  (b)
\end{minipage}
\begin{minipage}[c]{0.24\textwidth}
  \centering
  \includegraphics[width=30mm]{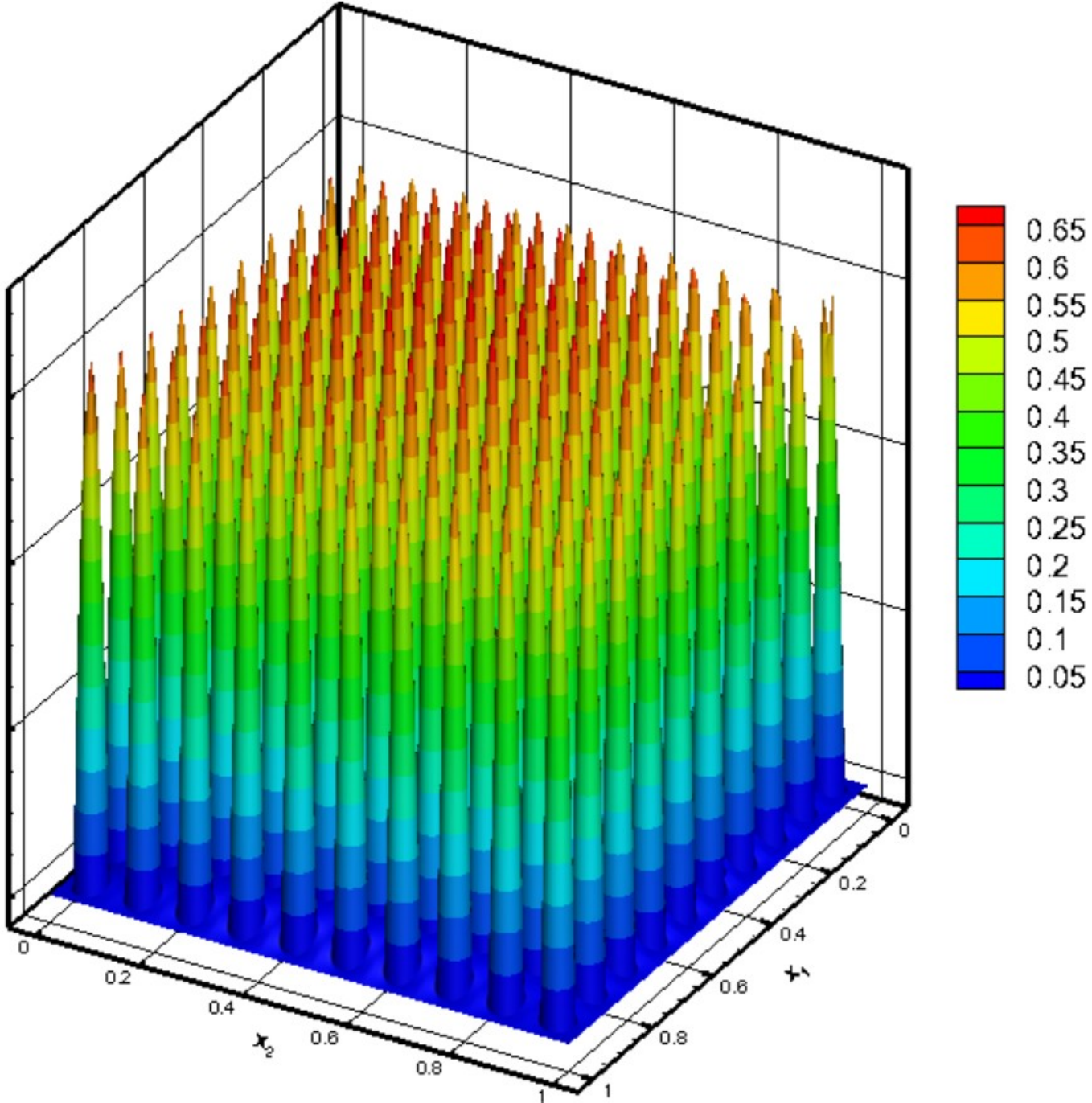}\\
  (c)
\end{minipage}
\begin{minipage}[c]{0.24\textwidth}
  \centering
  \includegraphics[width=30mm]{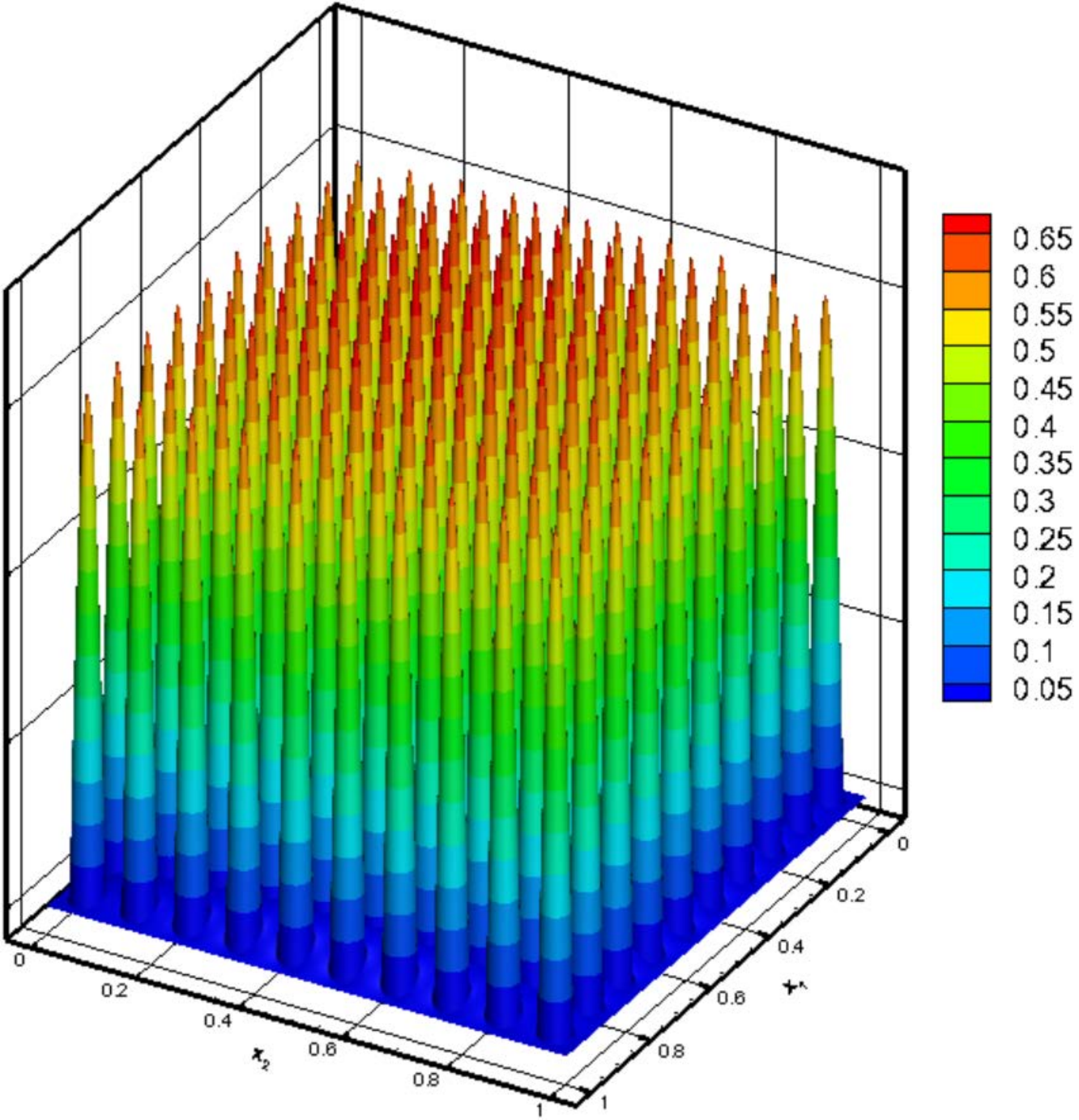}\\
  (d)
\end{minipage}
\caption{The electric potential field at $t=1.0$: (a) $\mathit{\Phi}_{0}$; (b) $\mathit{\Phi}^{[1\varepsilon]}$; (c) $\mathit{\Phi}^{[2\varepsilon]}$; (d) $\mathit{\Phi}_{{\rm{DNS}}}^\varepsilon$.}
\end{figure}
\begin{figure}[!htb]
\centering
\begin{minipage}[c]{0.24\textwidth}
  \centering
  \includegraphics[width=30mm]{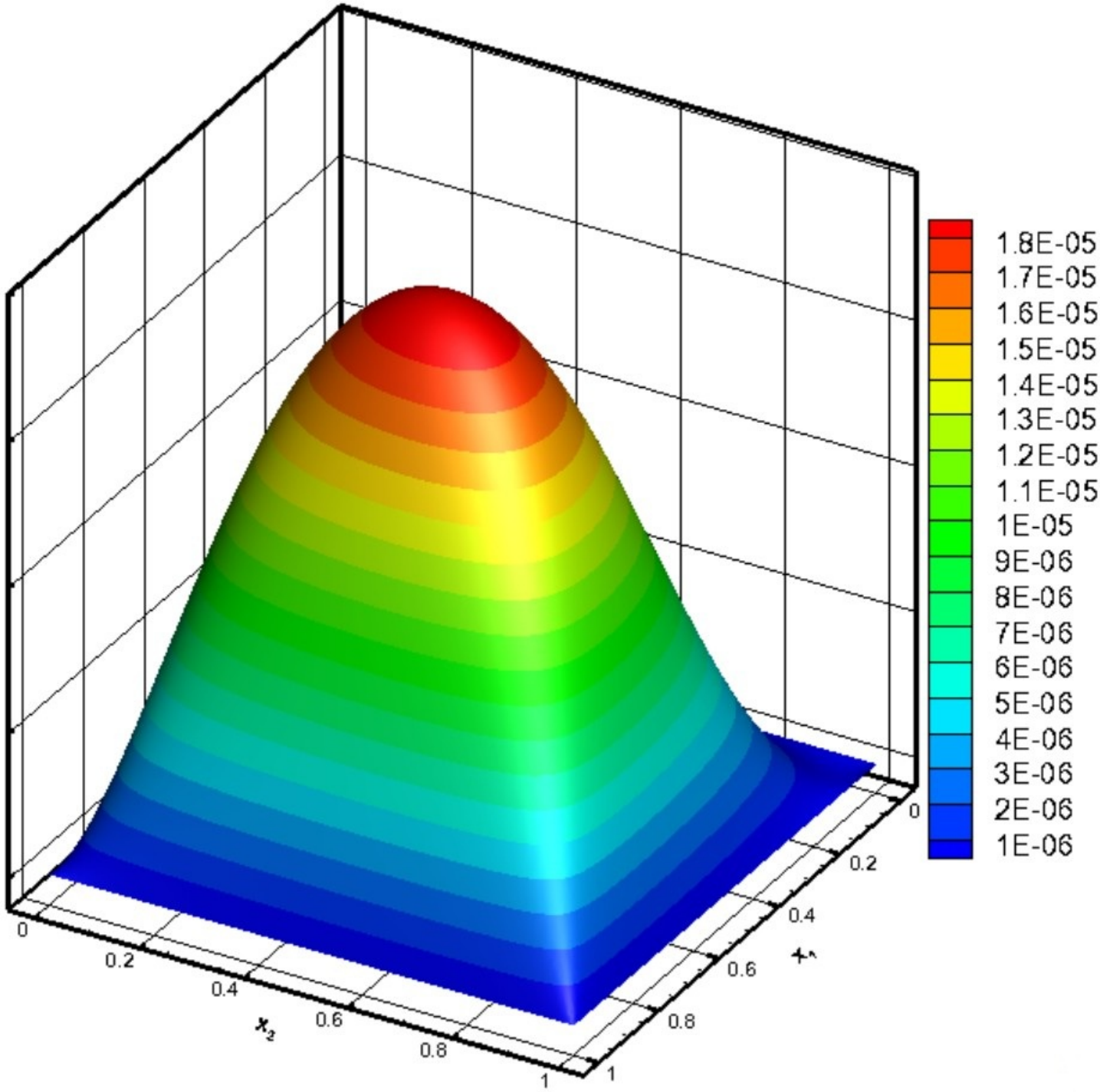}\\
  (a)
\end{minipage}
\begin{minipage}[c]{0.24\textwidth}
  \centering
  \includegraphics[width=30mm]{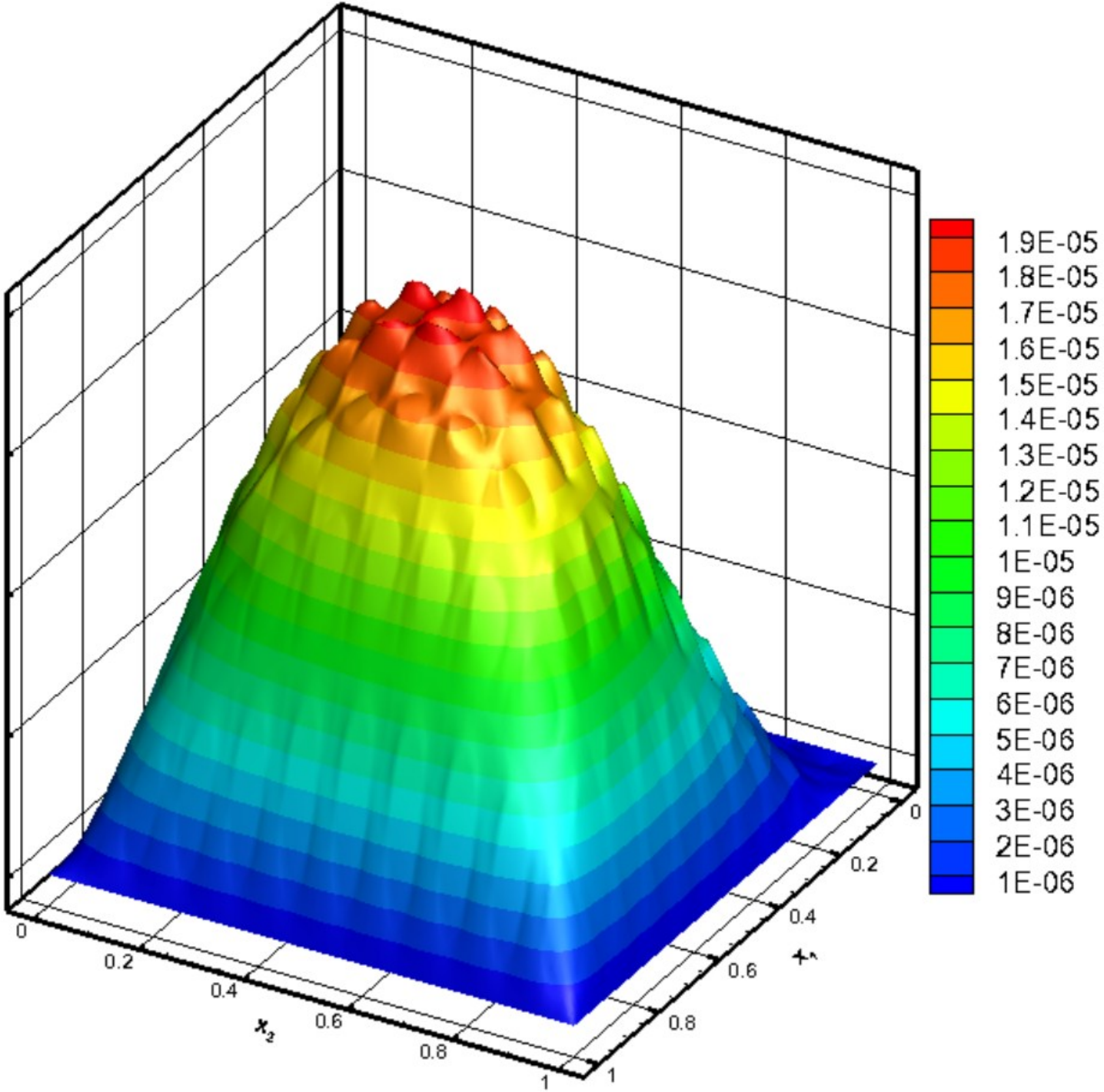}\\
  (b)
\end{minipage}
\begin{minipage}[c]{0.24\textwidth}
  \centering
  \includegraphics[width=30mm]{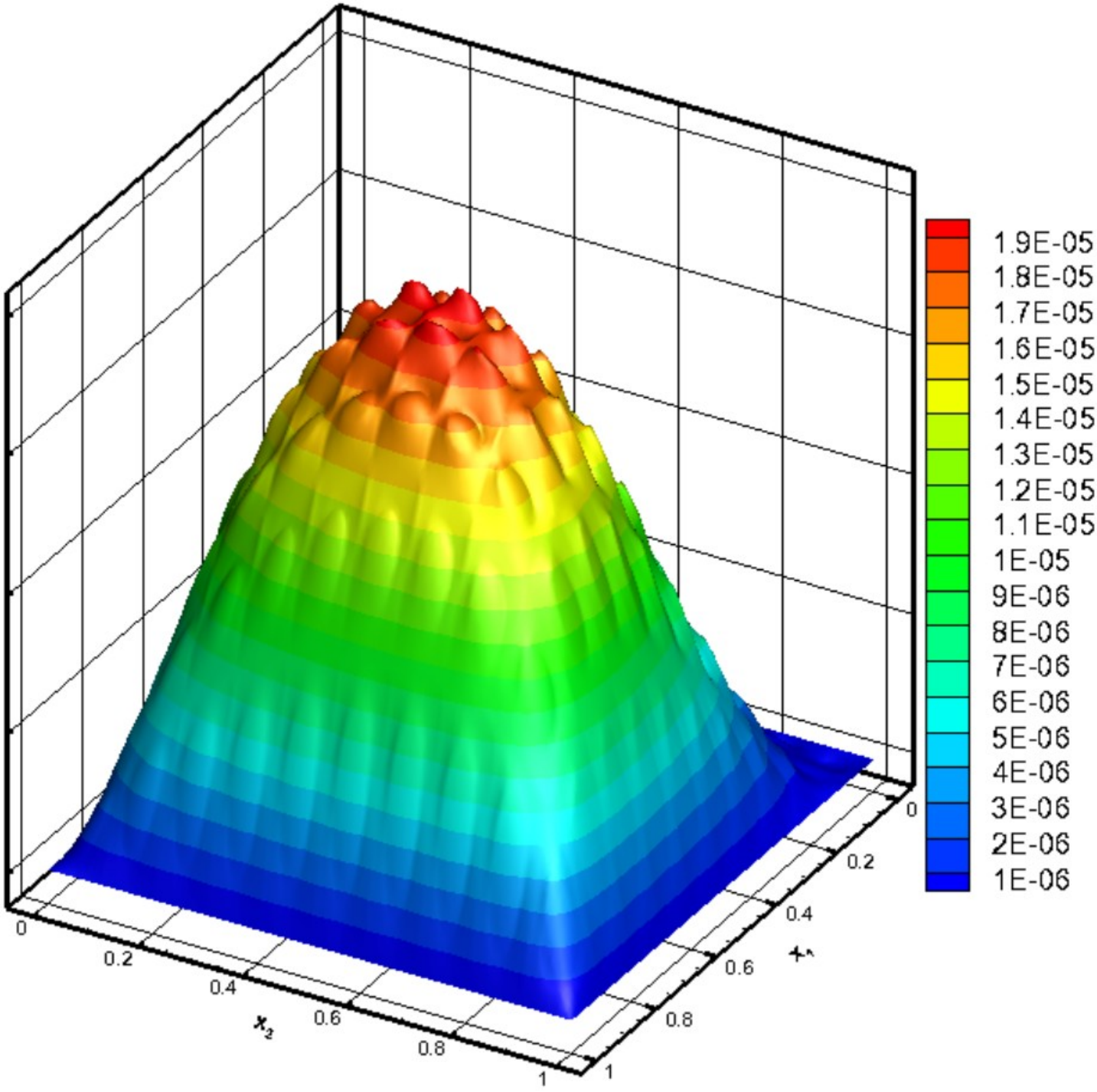}\\
  (c)
\end{minipage}
\begin{minipage}[c]{0.24\textwidth}
  \centering
  \includegraphics[width=30mm]{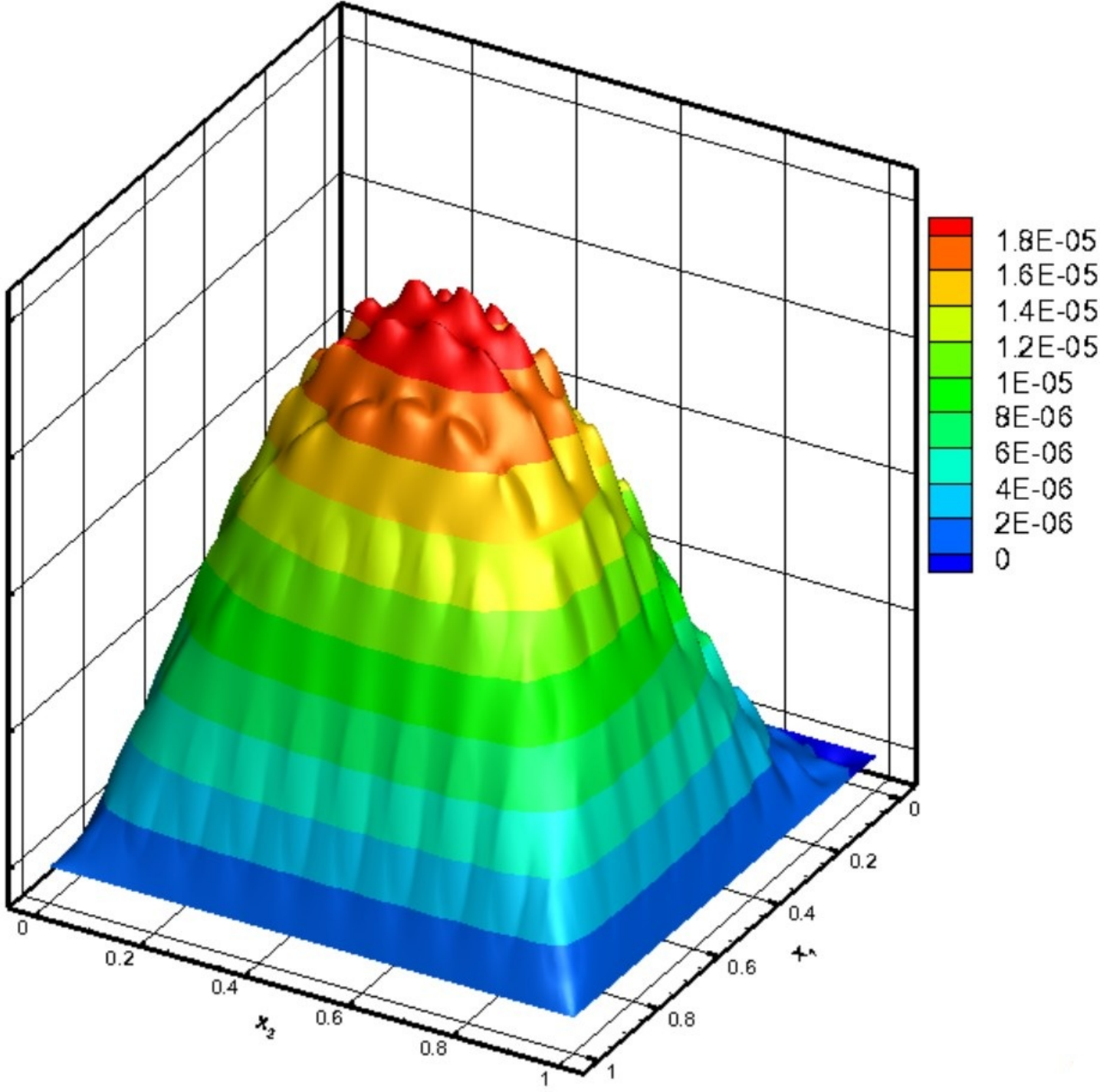}\\
  (d)
\end{minipage}
\caption{The displacement field at $t=1.0$: (a) $U_{10}$; (b) $U_1^{[1\varepsilon]}$; (c) $U_1^{[2\varepsilon]}$; (d) $U_{{1\rm{DNS}}}^\varepsilon$.}
\end{figure}
\begin{figure}[!htb]
\centering
\begin{minipage}[c]{0.42\textwidth}
  \centering
  \includegraphics[width=48mm]{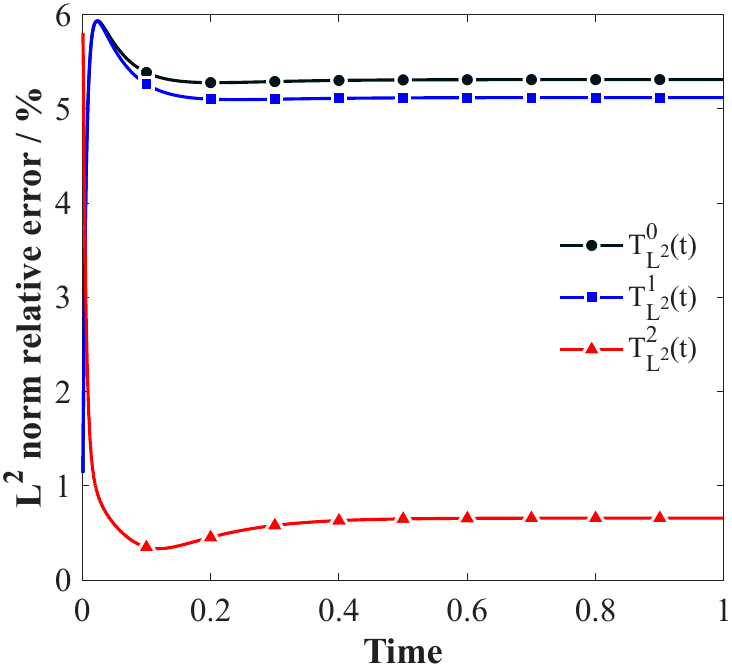}\\
  (a)
\end{minipage}
\begin{minipage}[c]{0.42\textwidth}
  \centering
  \includegraphics[width=48mm]{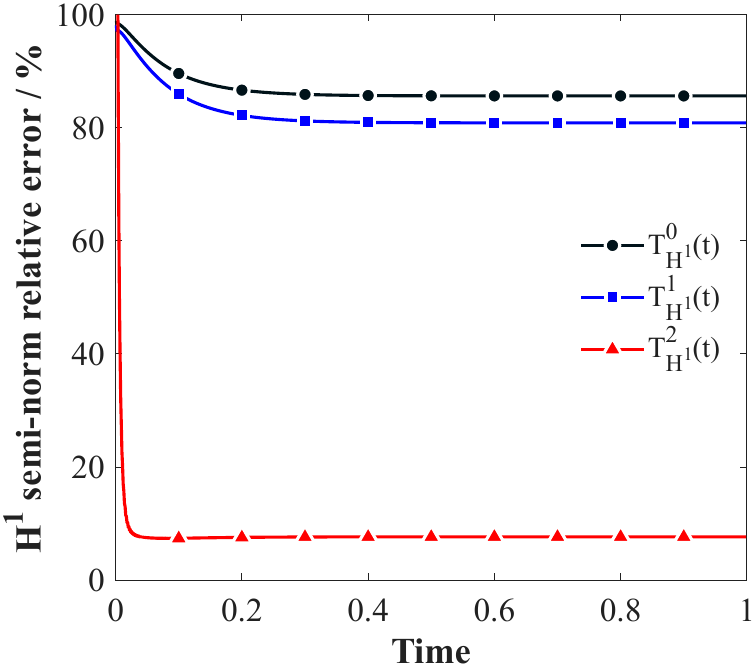}\\
  (b)
\end{minipage}
\begin{minipage}[c]{0.42\textwidth}
  \centering
  \includegraphics[width=48mm]{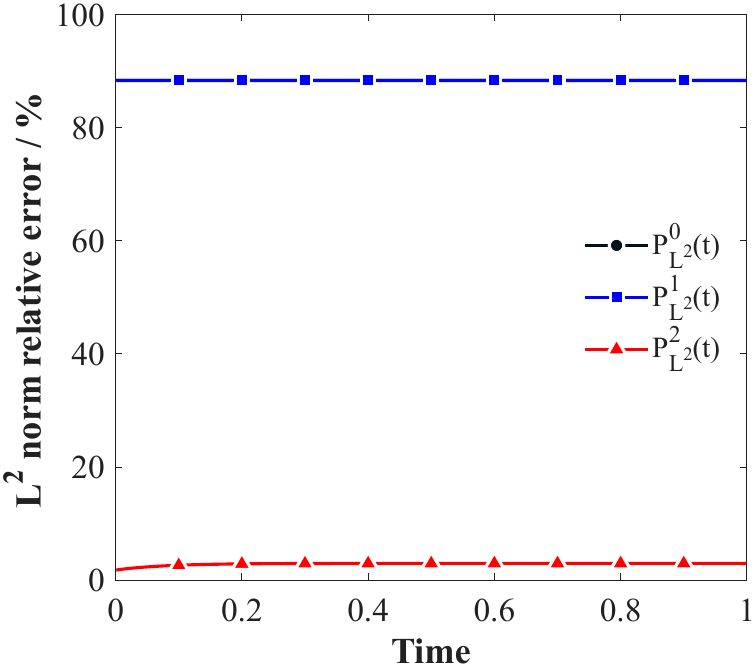}\\
  (c)
\end{minipage}
\begin{minipage}[c]{0.42\textwidth}
  \centering
  \includegraphics[width=48mm]{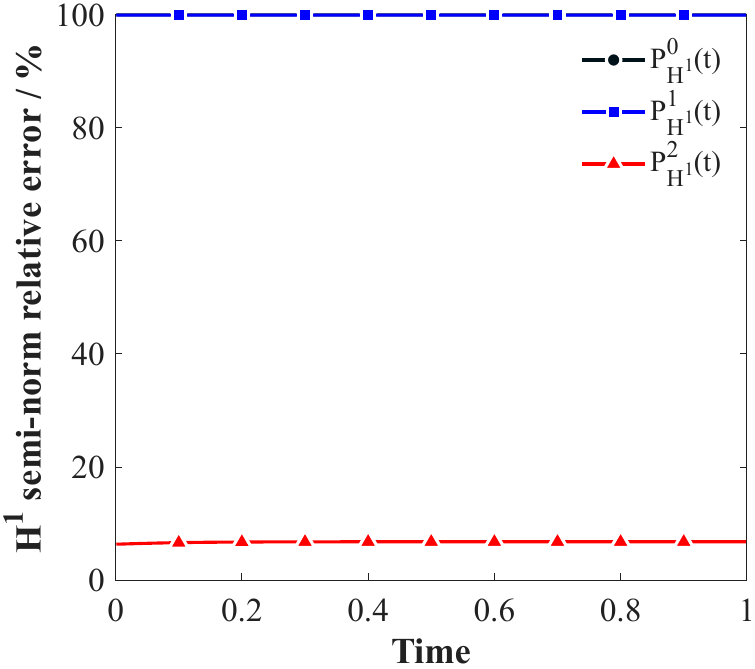}\\
  (d)
\end{minipage}
\begin{minipage}[c]{0.42\textwidth}
  \centering
  \includegraphics[width=48mm]{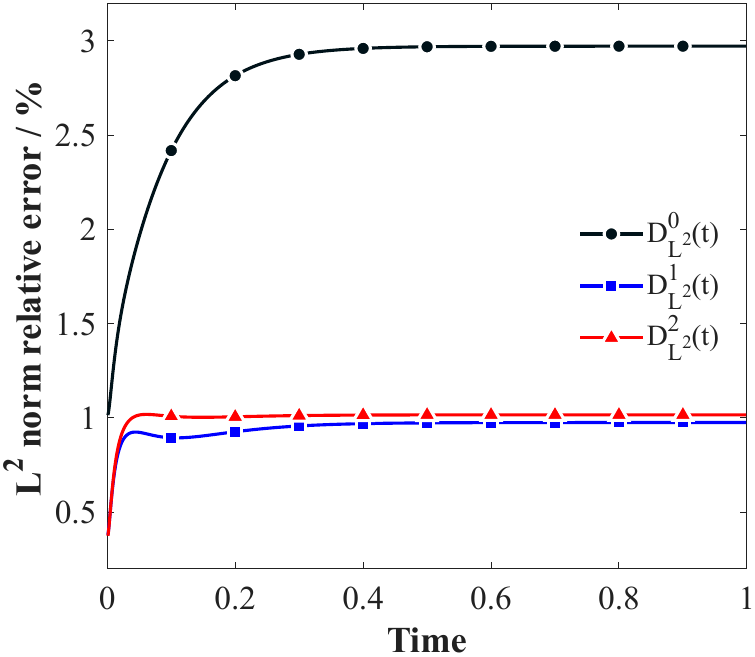}\\
  (e)
\end{minipage}
\begin{minipage}[c]{0.42\textwidth}
  \centering
  \includegraphics[width=48mm]{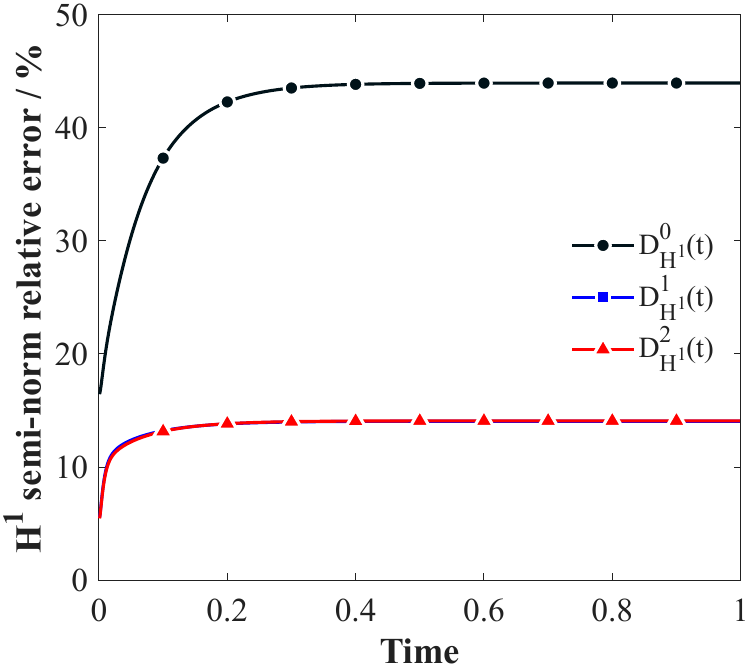}\\
  (f)
\end{minipage}
\caption{The evolutive relative errors of physical fields: (a) $\text{T}_{L^2}(t)$; (b) $\text{T}_{H^1}(t)$; (c) $\text{P}_{L^2}(t)$; (d) $\text{P}_{H^1}(t)$; (e) $\text{D}_{L^2}(t)$; (f) $\text{D}_{H^1}(t)$.}
\end{figure}

As reported in Table 2, the proposed HOMS approach yields a substantial reduction in computational resource utilization, particularly in terms of CPU memory and time, when contrasted with high-resolution DNS. Specifically,
the proposed HOMS approach achieves an approximate saving of 36.30$\%$ in computational time. Furthermore, the computational savings afforded by HOMS become more pronounced as the duration of the simulation increases. From the computational results in Figs.\hspace{1mm}2-4, it can be clearly found that only HOMS solutions have the capacity to accurately simulate the nonlinear thermo-electro-mechanical coupling behaviors of 2D heterogeneous structure, and the computational accuracy of homogenized and LOMS solutions is significantly inferior. Obviously, the homogenized method can only catch its macroscopic behaviors and the LOMS method can only capture its inadequate microscopic responses. Notably, as evidenced in Fig.\hspace{1mm}5, the HOMS approach maintains long-term numerical stability, exhibiting no blow-up phenomena up to the terminal time $t=1.0$. In real-world engineering applications, since direct FEM simulations are computationally prohibitive for large-scale composite structures involving a multitude of unit cells, the HOMS method proposed herein provides a efficient and robust alternative for multi-physics multi-scale simulations with minimal resource consumption.
\subsection{Example 2: nonlinear thermo-electro-mechanical coupling simulation of 3D composite structure}
This example evaluates the nonlinear thermo-electro-mechanical coupling behaviors of a 3D composite structure characterized by a periodic array of unit cells with the characteristic length parameter $\varepsilon=1/5$. The entire domain $\displaystyle\Omega=(x_1,x_2,x_3)=[0,1]\times[0,1]\times[0,1]$ and PUC $\mathbf{Y}$ are detailedly illustrated in Fig.\hspace{1mm}6.
\begin{figure}[!htb]
\centering
\begin{minipage}[c]{0.32\textwidth}
  \centering
  \includegraphics[width=0.60\linewidth,totalheight=1.1in]{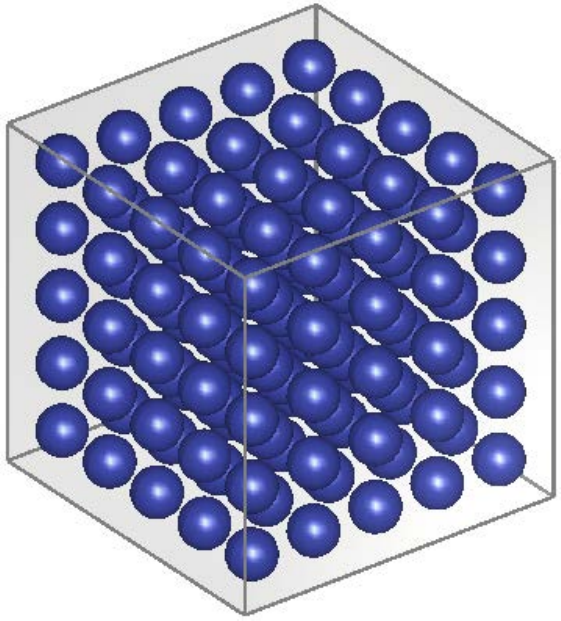} \\
  (a)
\end{minipage}
\begin{minipage}[c]{0.32\textwidth}
  \centering
  \includegraphics[width=0.98\linewidth,totalheight=1.1in]{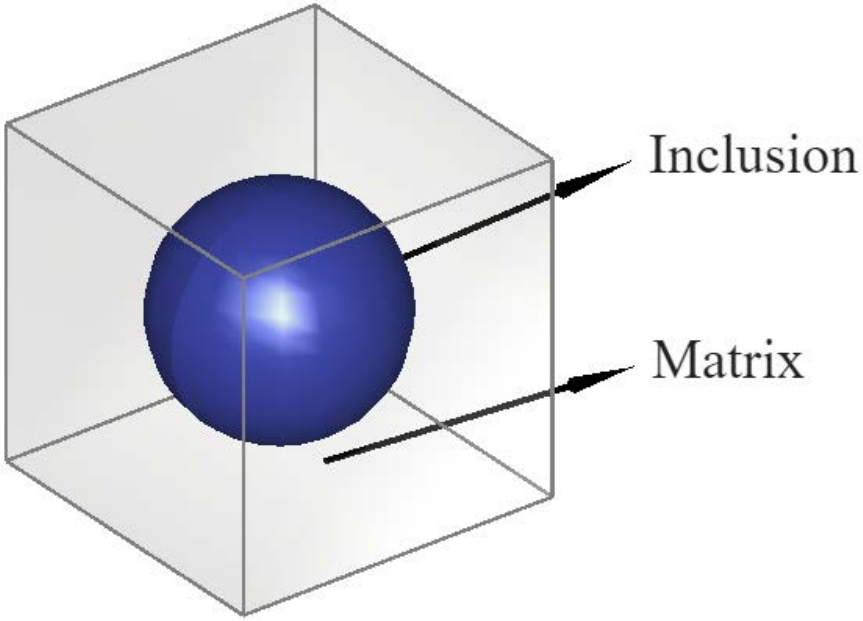} \\
  (b)
\end{minipage}
\begin{minipage}[c]{0.32\textwidth}
  \centering
  \includegraphics[width=0.60\linewidth,totalheight=1.1in]{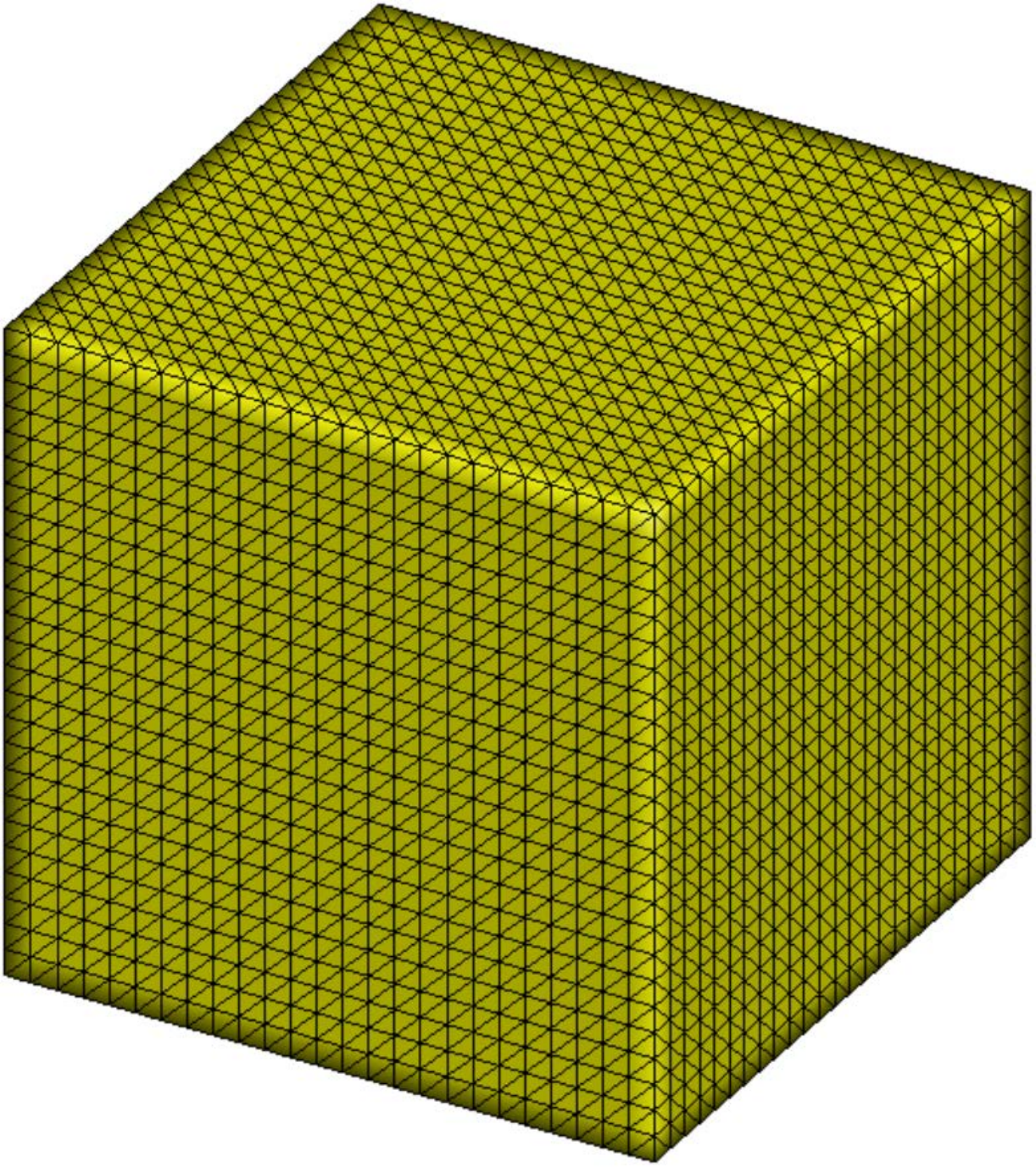} \\
  (c)
\end{minipage}
\caption{(a) The 3D composite structure $\Omega$; (b) PUC $\mathbf{Y}$; (c) homogenized structure $\Omega$.}
\end{figure}

The material parameters of this example are defined as those of Example 1. Furthermore, the internal heat source, electric charge density and body forces are prescribed as $f_T(\bm{x},t)=80000x_1(1-x_1)x_2(1-x_2)$, $f_\mathit{\Phi}(\bm{x},t)=800x_1(1-x_1)x_2(1-x_2)$ and $(f_1(\bm{x},t),f_2(\bm{x},t),f_3(\bm{x},t))=(0,0,20000)$. The boundary conditions of temperature and electric potential fields are defined as $\widehat{T}(\bm{x},t)=300.0$ and $\widehat{\mathit{\Phi}}(\bm{x},t)=0.0$. The boundary condition of displacement field on the four lateral surfaces of the investigated 3D composite structure, that are perpendicular to $x_1$-axis and $x_2$-axis, is enforced as $\widehat{\bm{U}}(\bm{x},t)=(0.0,0.0,0.0)$. The initial conditions of temperature and displacement fields are given by $\widetilde T=300.0$, $\widetilde{\bm{U}}^0({\bm{x}})=(0.0,0.0,0.0)$, and $\widetilde{\bm{U}}^1({\bm{x}})=(0.0,0.0,0.0)$.

Direct numerical simulation of the 3D heterogeneous structure considered herein, which contains a vast number of unit cells, necessitates extremely fine meshes and incurs unaffordable computational costs in turn. For the purpose of comparison, finite element discretizations are constructed for original multi-scale system, auxiliary cell problems, and corresponding homogenized system. A detailed breakdown of the mesh resolutions and execution times is provided in Table 3.
\begin{table}[!htb]{\caption{Summary of computational cost in 3D composite structure.}}
\centering
\begin{tabular}{cccc}
\hline
 & Multi-scale eqs. & Cell eqs. & Homogenized eqs. \\
\hline
FEM elements & 1275910 & 139775 & 162000\\
FEM nodes    & 206563 & 23807 & 29791\\
\hline
& DNS & off-line stage & on-line stage \\
\hline
Computational time & 405368.064s  & 53542.102s & 231230.402s\\
\hline
\end{tabular}
\end{table}

For this example, we allocate 10 equidistant macroscopic interpolation temperature within a single unit cell. It should note that although solving the auxiliary cell problems constitutes a computationally expensive off-line stage, these obtained data can be reused for various heterogeneous structures during the on-line phase. The dynamic nonlinear thermo-electro-mechanical coupling responses of the 3D heterogeneous structure are simulated over the time interval $t\in[0,1]$. Setting the temporal step as $\Delta t = 0.001$, the macroscopic homogenized equations (2.24) and multi-scale nonlinear equations (1.1) are on-line simulated respectively. Following this, the final distributions of temperature, electric potential and displacement fields at $t=1.0$ are illustrated in Figs.\hspace{1mm}7-10, respectively. Besides, the temporal evolution of the relative errors of three physical fields in the $L^2$ norm and $H^1$ semi-norm senses is visualized in Fig.\hspace{1mm}11.
\begin{figure}[!htb]
\centering
\begin{minipage}[c]{0.24\textwidth}
  \centering
  \includegraphics[width=31mm]{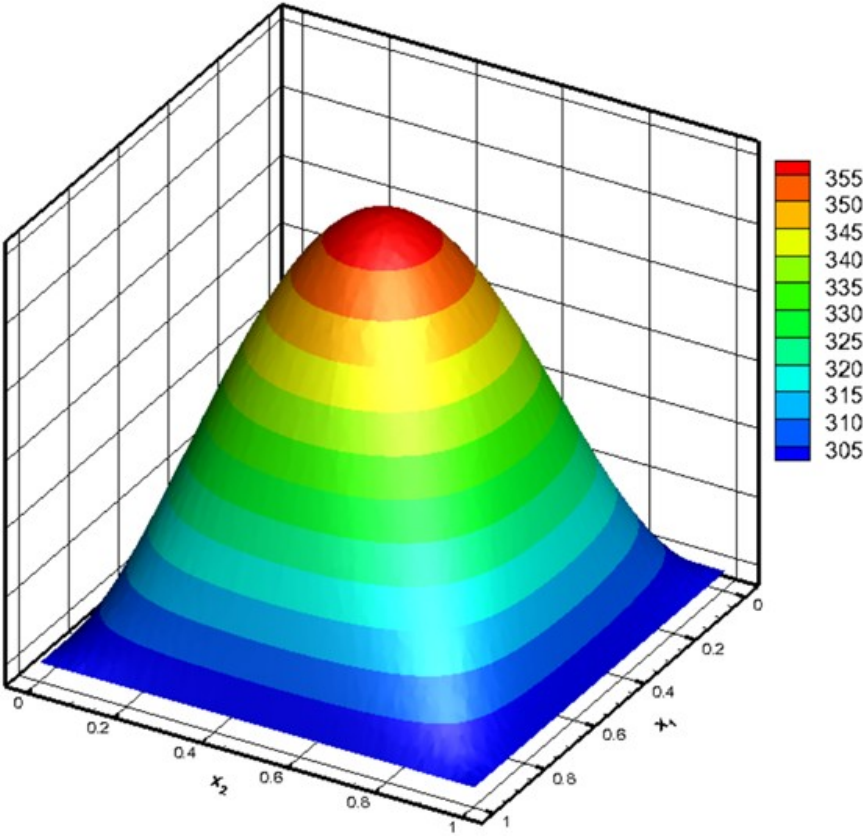}\\
  (a)
\end{minipage}
\begin{minipage}[c]{0.24\textwidth}
  \centering
  \includegraphics[width=31mm]{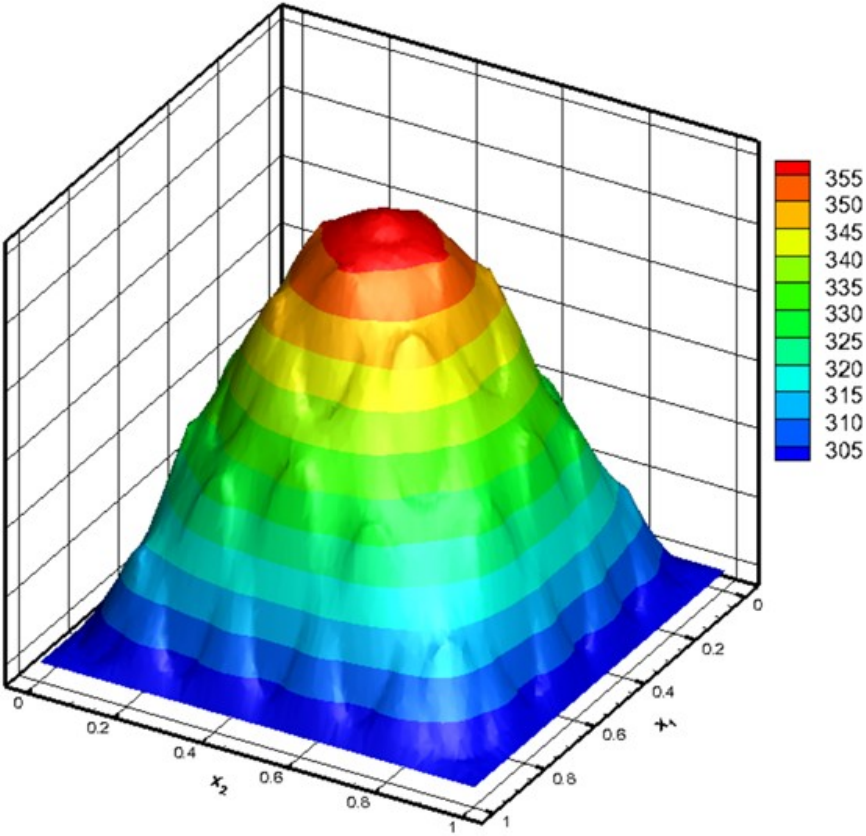}\\
  (b)
\end{minipage}
\begin{minipage}[c]{0.24\textwidth}
  \centering
  \includegraphics[width=31mm]{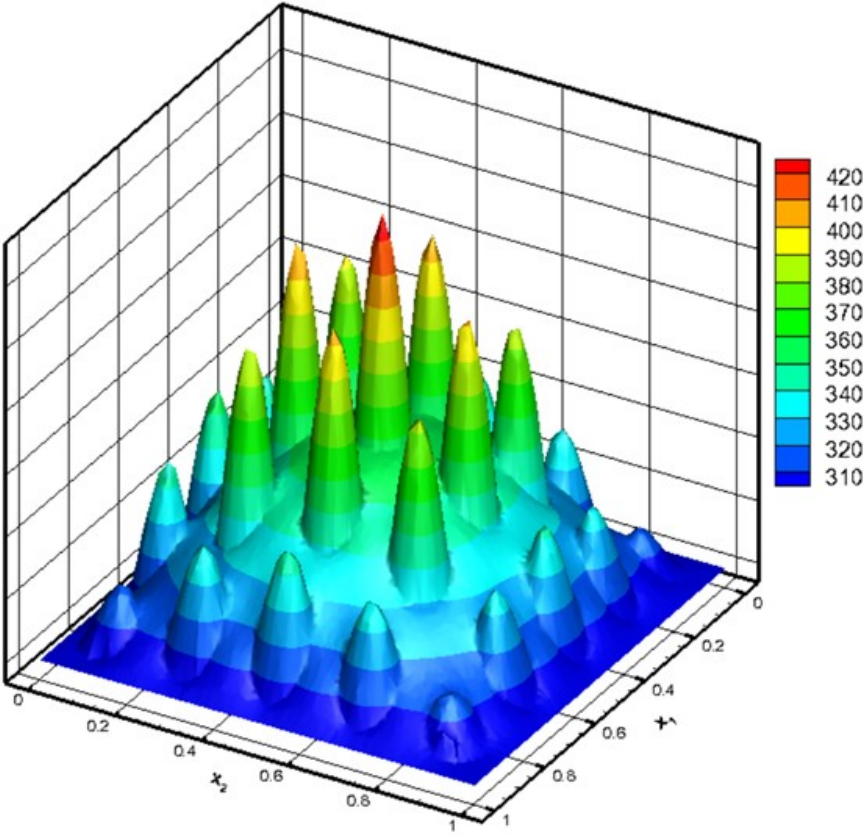}\\
  (c)
\end{minipage}
\begin{minipage}[c]{0.24\textwidth}
  \centering
  \includegraphics[width=31mm]{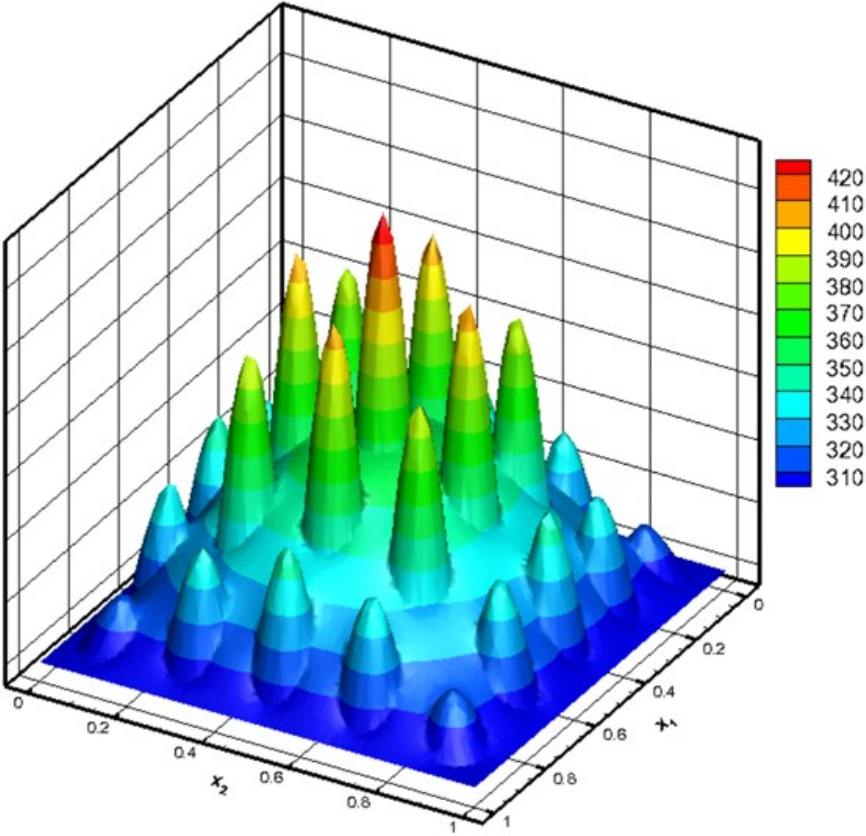}\\
  (d)
\end{minipage}
\caption{The temperature field in cross section $x_3=0.5$ at $t=1.0$: (a) $T_{0}$; (b) $T^{[1\varepsilon]}$; (c) $T^{[2\varepsilon]}$; (d) $T_{\rm{DNS}}^\varepsilon$.}
\end{figure}
\begin{figure}[!htb]
\centering
\begin{minipage}[c]{0.24\textwidth}
  \centering
  \includegraphics[width=30mm]{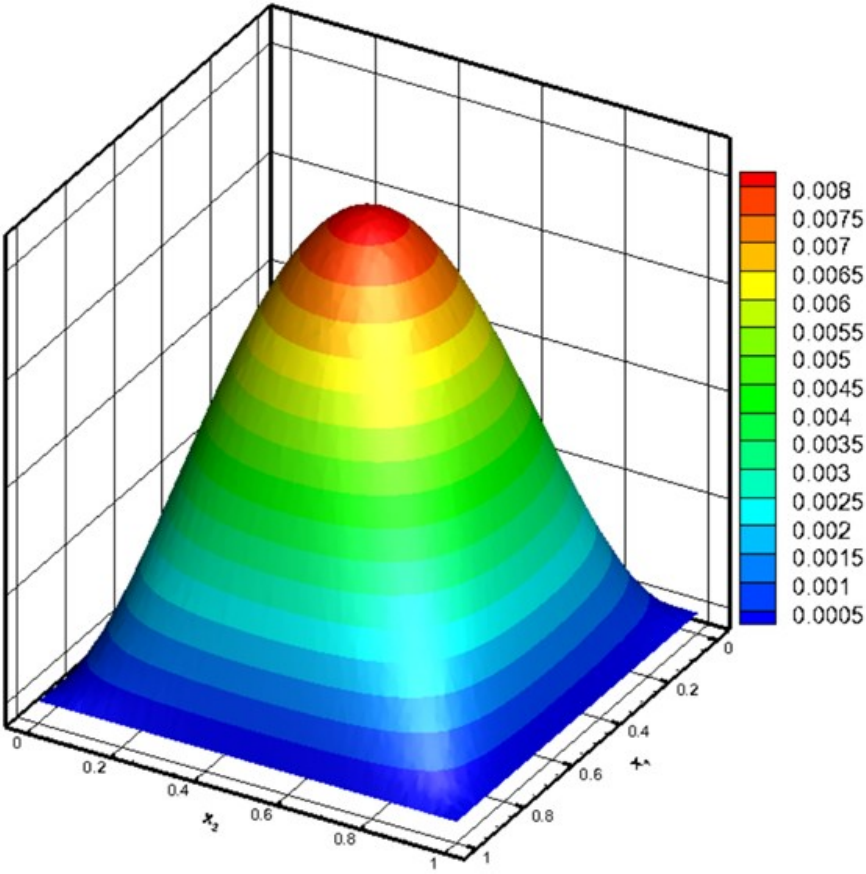}\\
  (a)
\end{minipage}
\begin{minipage}[c]{0.24\textwidth}
  \centering
  \includegraphics[width=30mm]{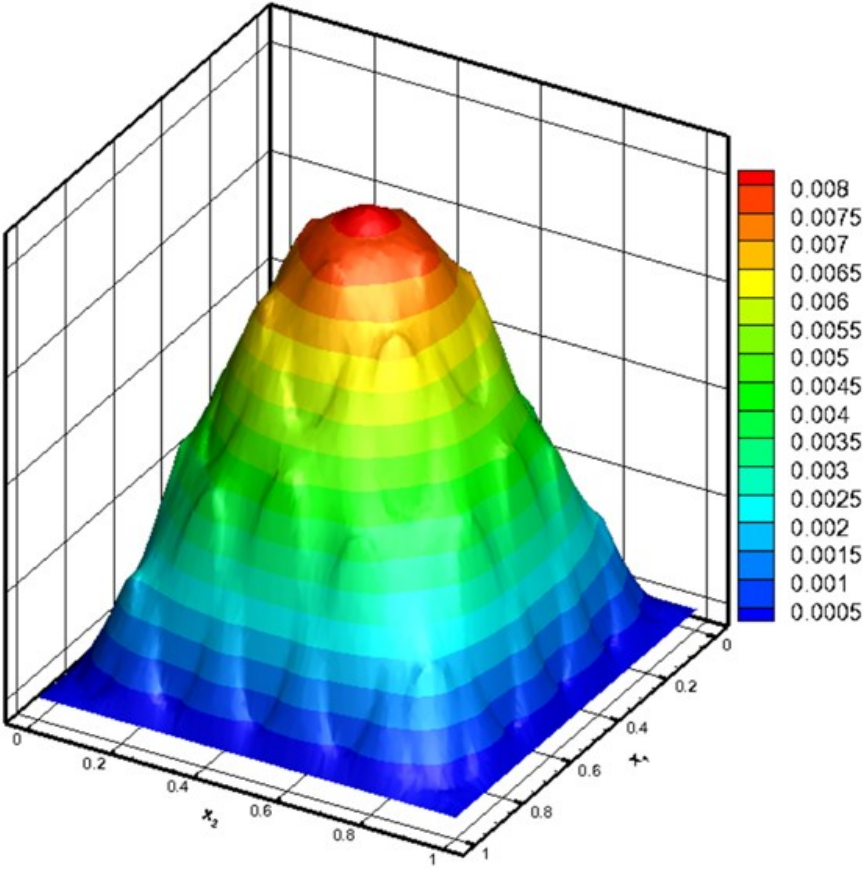}\\
  (b)
\end{minipage}
\begin{minipage}[c]{0.24\textwidth}
  \centering
  \includegraphics[width=30mm]{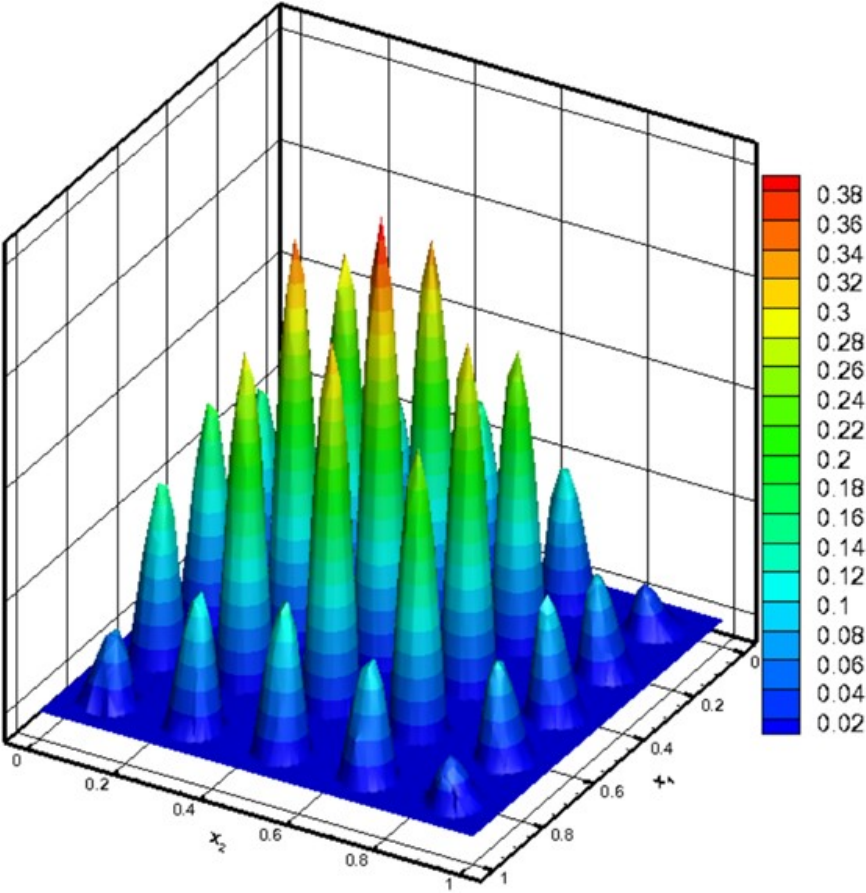}\\
  (c)
\end{minipage}
\begin{minipage}[c]{0.24\textwidth}
  \centering
  \includegraphics[width=30mm]{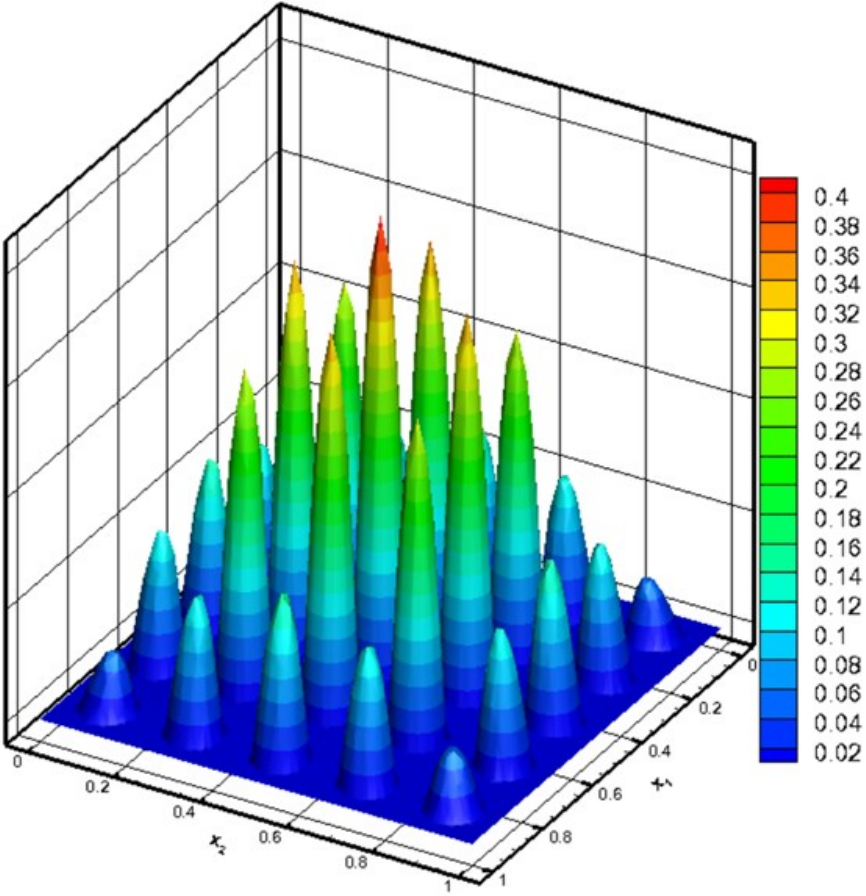}\\
  (d)
\end{minipage}
\caption{The electric potential field in cross section $x_3=0.5$ at $t=1.0$: (a) $\mathit{\Phi}_{0}$; (b) $\mathit{\Phi}^{[1\varepsilon]}$; (c) $\mathit{\Phi}^{[2\varepsilon]}$; (d) $\mathit{\Phi}_{{\rm{DNS}}}^\varepsilon$.}
\end{figure}
\begin{figure}[!htb]
\centering
\begin{minipage}[c]{0.24\textwidth}
  \centering
  \includegraphics[width=30mm]{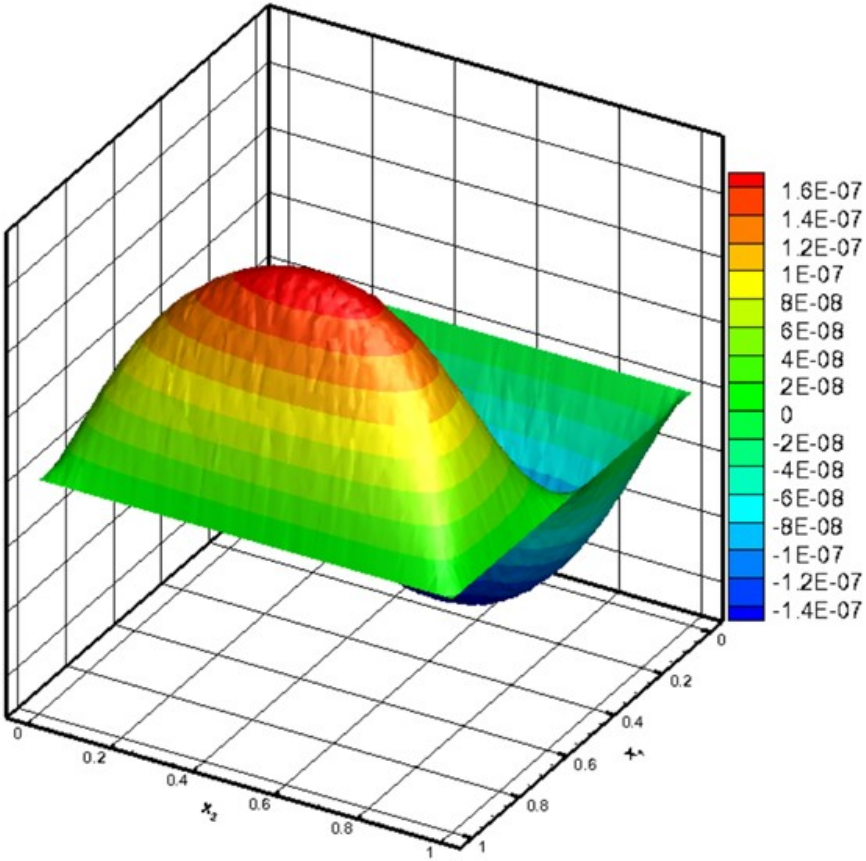}\\
  (a)
\end{minipage}
\begin{minipage}[c]{0.24\textwidth}
  \centering
  \includegraphics[width=30mm]{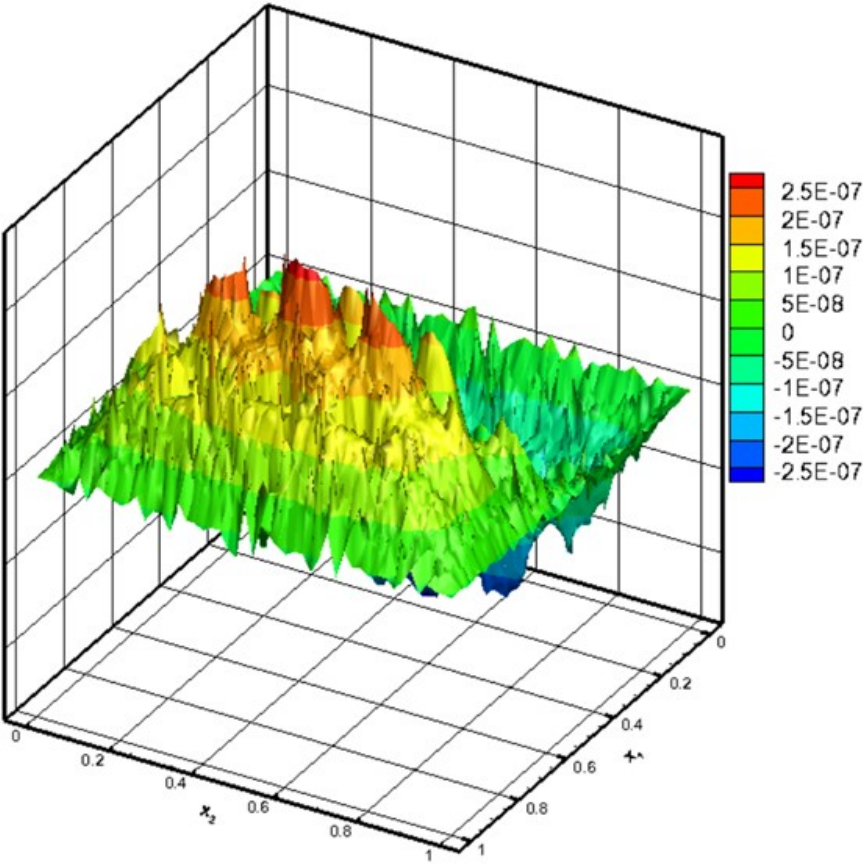}\\
  (b)
\end{minipage}
\begin{minipage}[c]{0.24\textwidth}
  \centering
  \includegraphics[width=30mm]{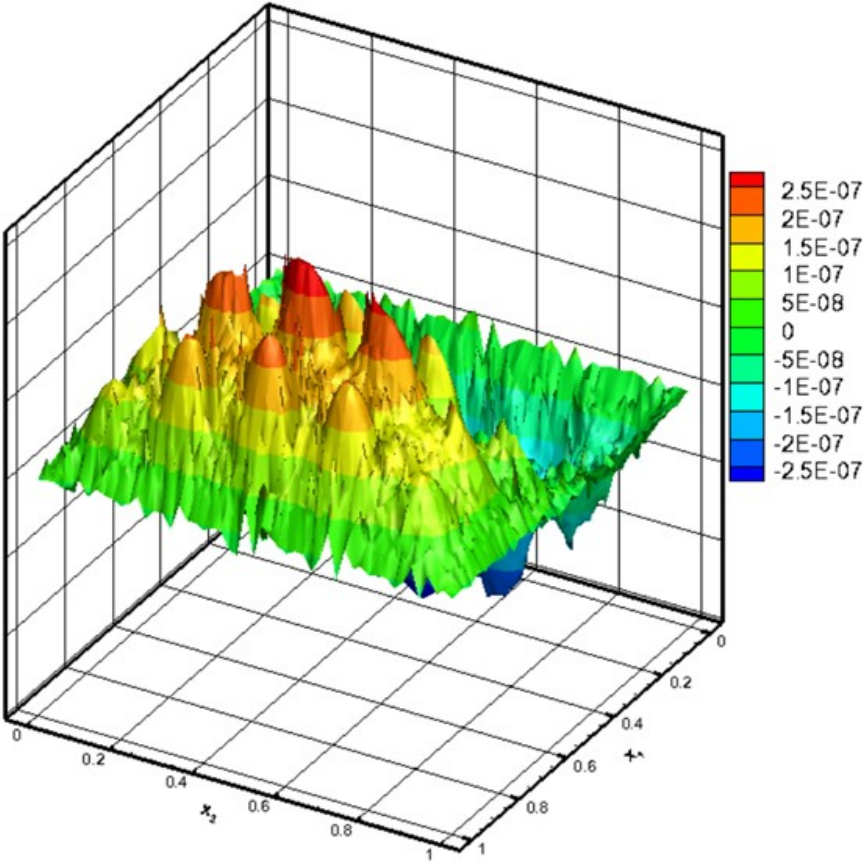}\\
  (c)
\end{minipage}
\begin{minipage}[c]{0.24\textwidth}
  \centering
  \includegraphics[width=30mm]{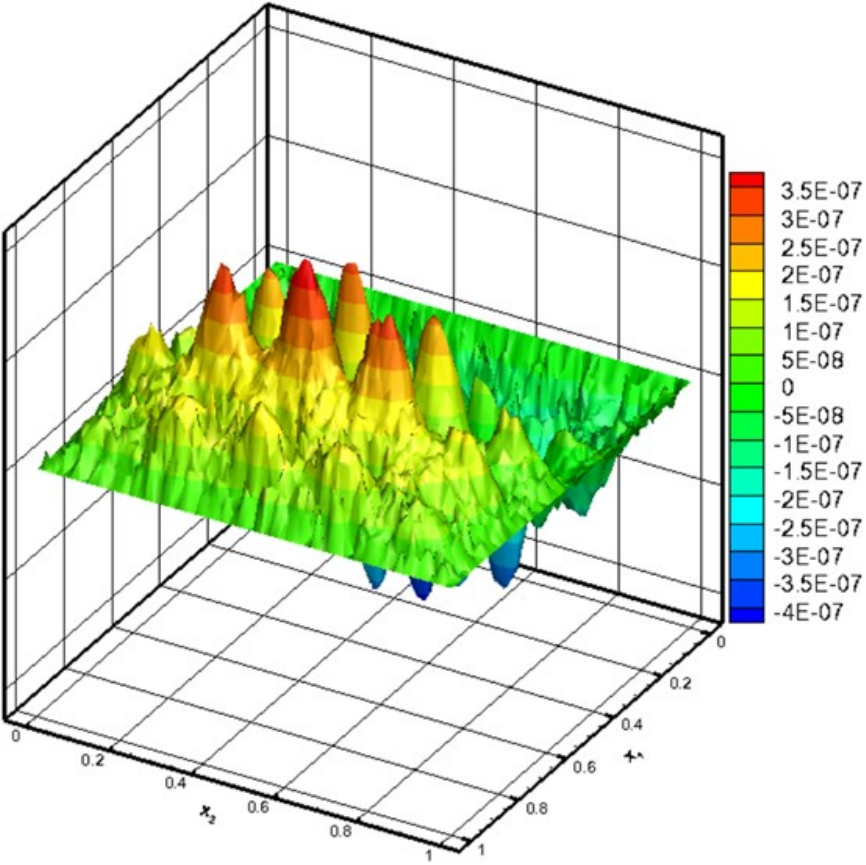}\\
  (d)
\end{minipage}
\caption{The displacement field in cross section $x_3=0.5$ at $t=1.0$: (a) $U_{10}$; (b) $U_1^{[1\varepsilon]}$; (c) $U_1^{[2\varepsilon]}$; (d) $U_{{1\rm{DNS}}}^\varepsilon$.}
\end{figure}
\begin{figure}[!htb]
\centering
\begin{minipage}[c]{0.24\textwidth}
  \centering
  \includegraphics[width=30mm]{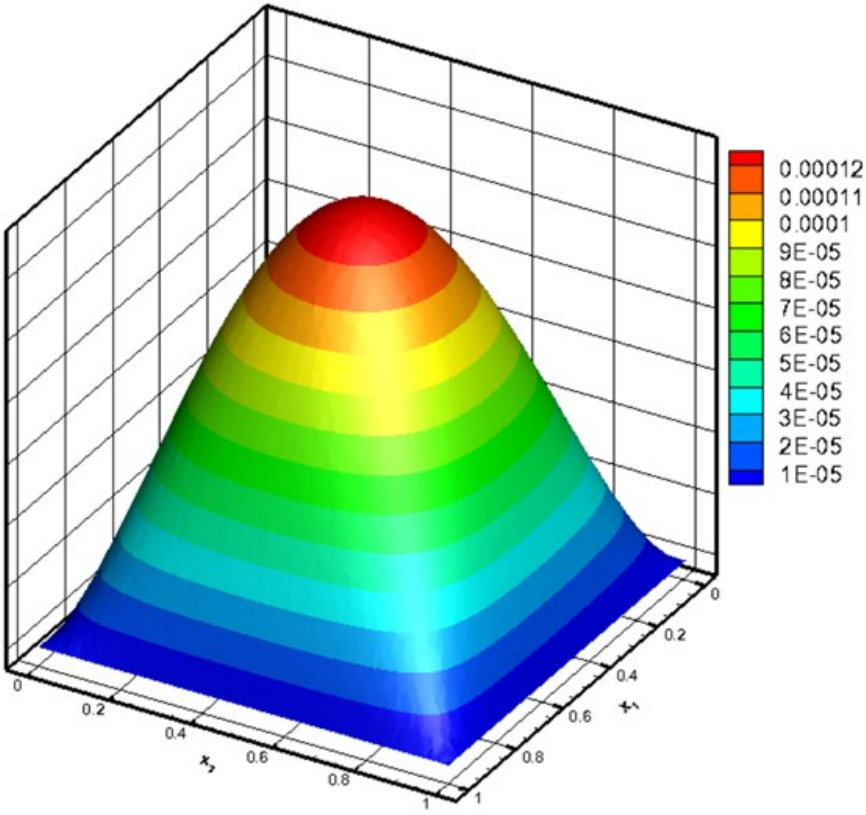}\\
  (a)
\end{minipage}
\begin{minipage}[c]{0.24\textwidth}
  \centering
  \includegraphics[width=30mm]{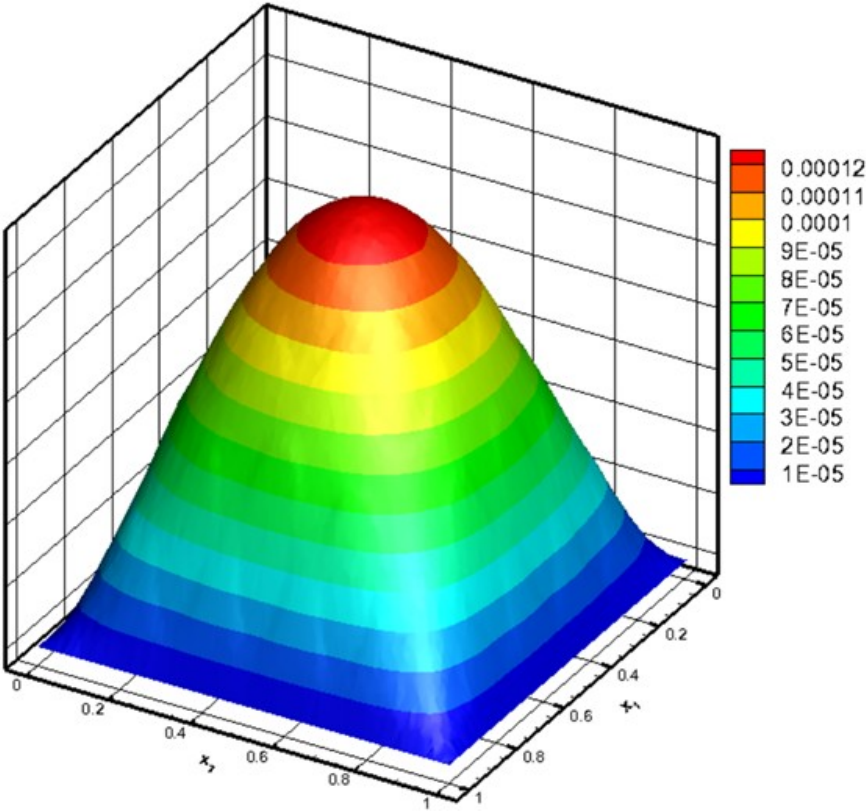}\\
  (b)
\end{minipage}
\begin{minipage}[c]{0.24\textwidth}
  \centering
  \includegraphics[width=30mm]{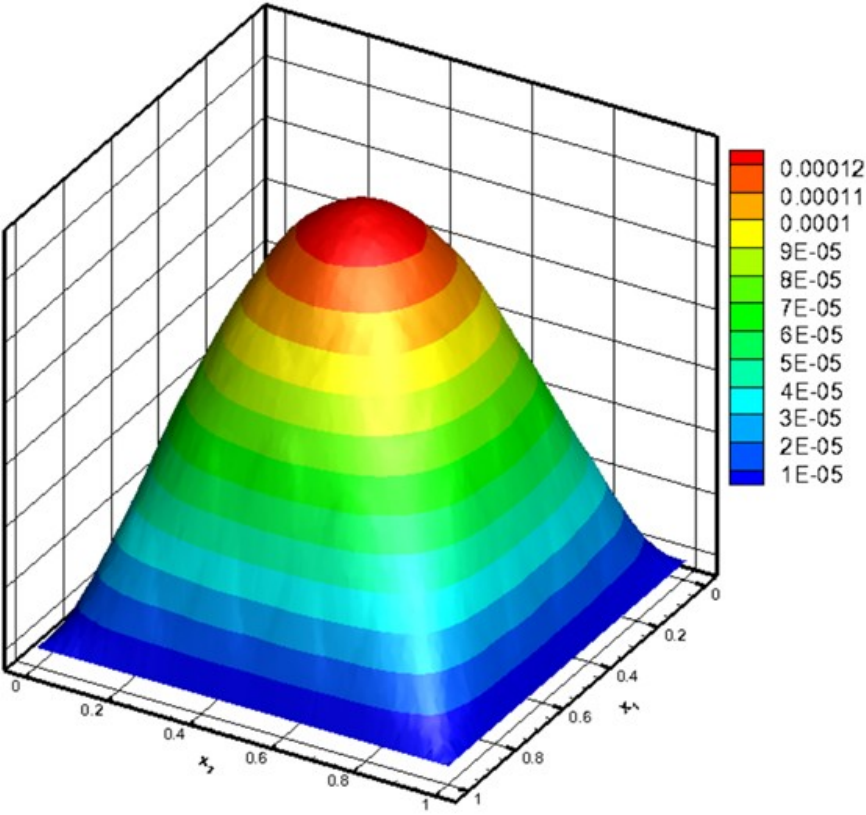}\\
  (c)
\end{minipage}
\begin{minipage}[c]{0.24\textwidth}
  \centering
  \includegraphics[width=30mm]{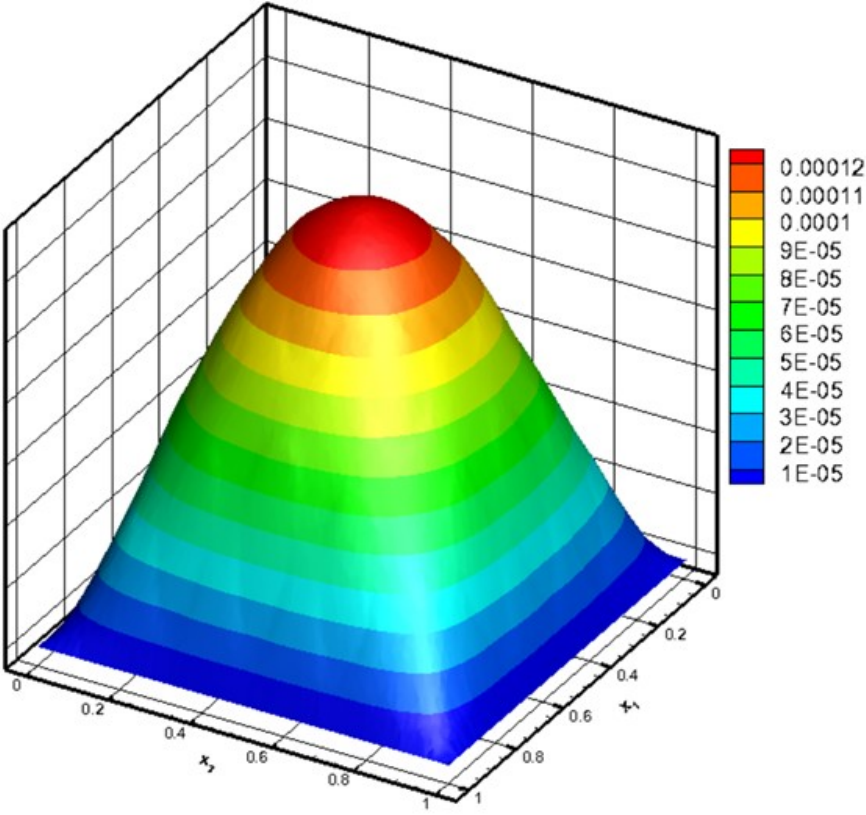}\\
  (d)
\end{minipage}
\caption{The displacement field in cross section $x_3=0.5$ at $t=1.0$: (a) $U_{30}$; (b) $U_3^{[1\varepsilon]}$; (c) $U_3^{[2\varepsilon]}$; (d) $U_{{3\rm{DNS}}}^\varepsilon$.}
\end{figure}
\begin{figure}[!htb]
\centering
\begin{minipage}[c]{0.42\textwidth}
  \centering
  \includegraphics[width=48mm]{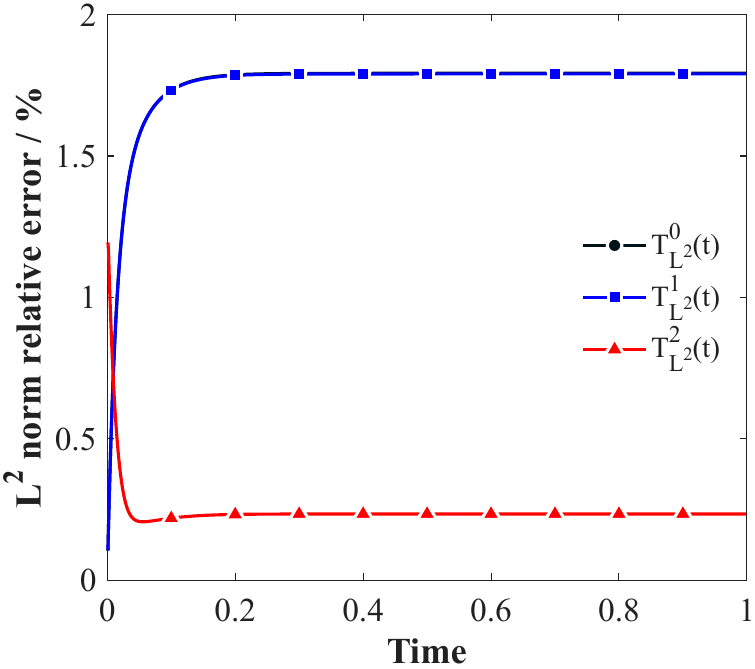}\\
  (a)
\end{minipage}
\begin{minipage}[c]{0.42\textwidth}
  \centering
  \includegraphics[width=48mm]{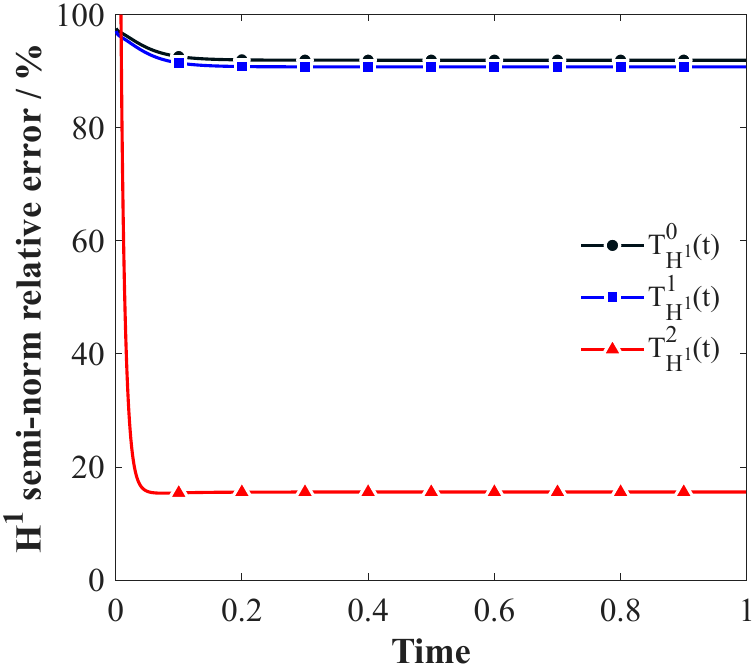}\\
  (b)
\end{minipage}
\begin{minipage}[c]{0.42\textwidth}
  \centering
  \includegraphics[width=48mm]{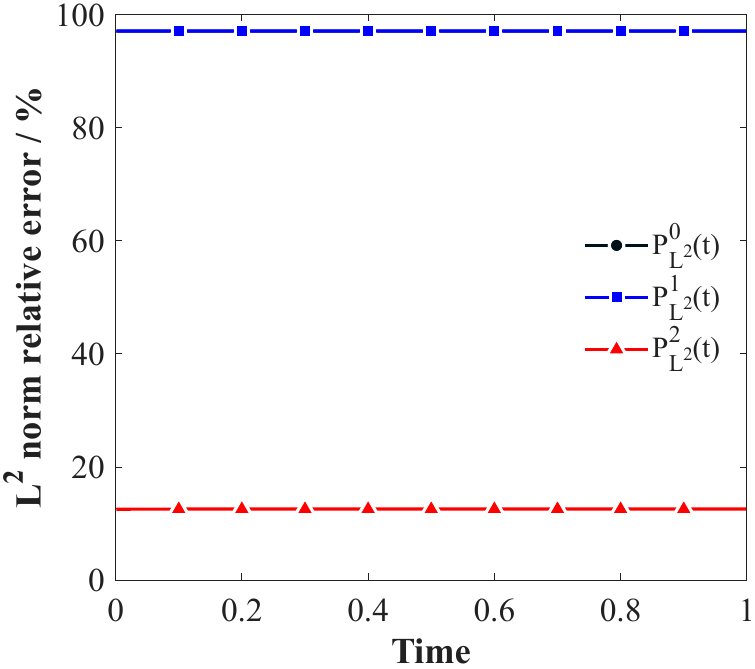}\\
  (c)
\end{minipage}
\begin{minipage}[c]{0.42\textwidth}
  \centering
  \includegraphics[width=48mm]{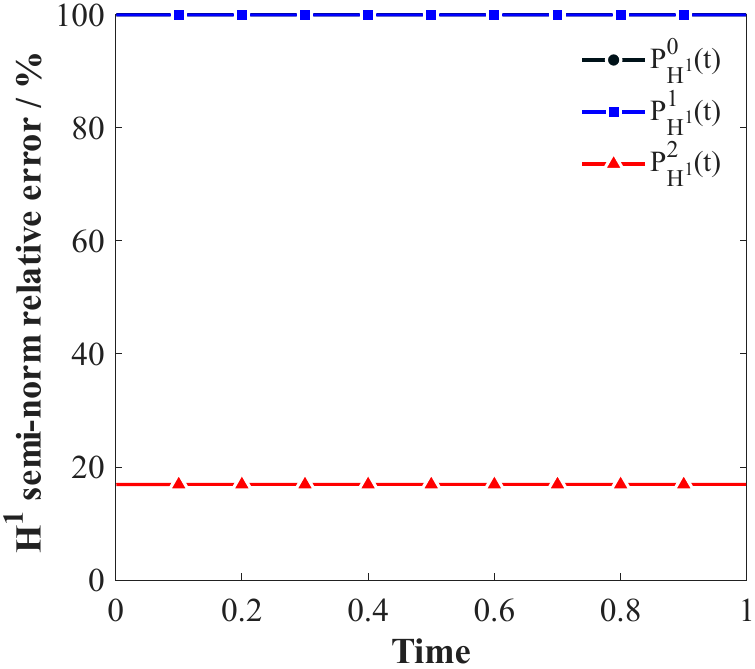}\\
  (d)
\end{minipage}
\begin{minipage}[c]{0.42\textwidth}
  \centering
  \includegraphics[width=48mm]{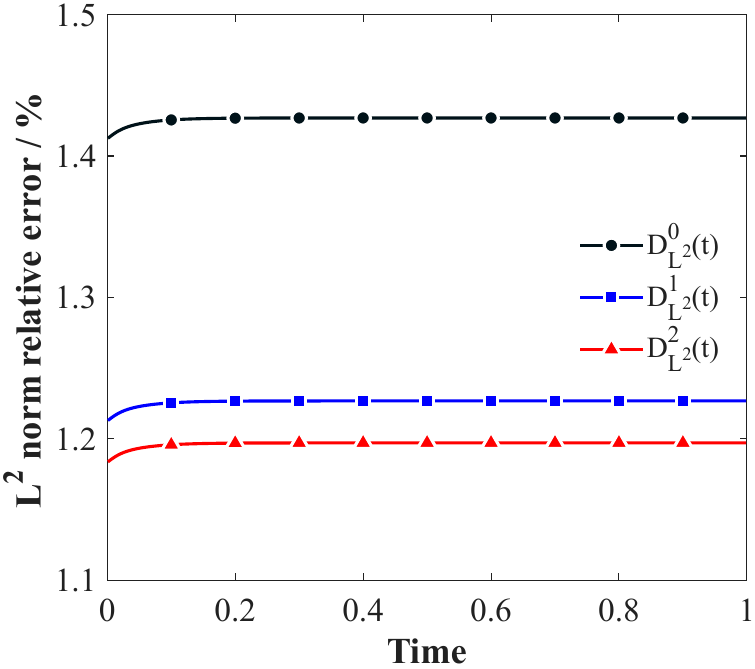}\\
  (e)
\end{minipage}
\begin{minipage}[c]{0.42\textwidth}
  \centering
  \includegraphics[width=48mm]{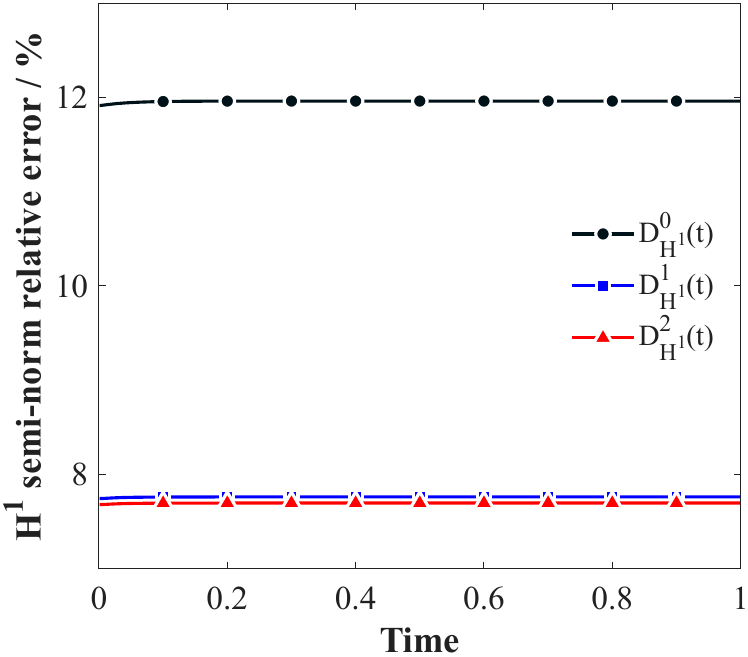}\\
  (f)
\end{minipage}
\caption{The evolutive relative errors of physical fields: (a) $\text{T}_{L^2}(t)$; (b) $\text{T}_{H^1}(t)$; (c) $\text{P}_{L^2}(t)$; (d) $\text{P}_{H^1}(t)$; (e) $\text{D}_{L^2}(t)$; (f) $\text{D}_{H^1}(t)$.}
\end{figure}

Table 3 shows that the HOMS method significantly reduces computational memory and time requirements relative to high-resolution DNS, yielding a 29.75$\%$ decrease in runtime. This efficiency gain scales favorably with the simulation time. In principle, greater reductions in computational cost could be realized by the proposed method. Nevertheless, to ensure a fair comparison with the reference solution, the time step size was chosen to be as small as that utilized in the high-resolution DNS. Besides, the numerical results in Figs.\hspace{1mm}7-10 demonstrate that only the HOMS solution can accurately capture the nonlinear thermo-electro-mechanical coupling behaviors of the 3D heterogeneous structure, whereas the standard homogenized method resolves only the macroscopic response, and the LOMS method provides insufficient microscopic corrections. In terms of temporal robustness, Fig.\hspace{1mm}11 effectively validate that the proposed HOMS approach can remain entirely stable without blow-up in the numerical solutions until $t=1.0$. These attributes make the HOMS approach particularly well-suited for complex 3D composite structures with numerous unit cells, circumventing the prohibitive computational cost of direct FEM simulations in large-scale multi-physics applications.
\section{Conclusions}
This work proposes a low-computational-cost high-order multi-scale method for high-accuracy simulation of time-dependent nonlinear thermo-electro-mechanical coupling problems in composite structures with highly microscopic heterogeneities. The pivotal contributions of this work are threefold: First, the novel macro-micro correlative formulations with the high-order correction terms are established for composite structures with periodically microscopic configurations. Second, rigorous local and global error analyses for multi-scale solutions of composite structures are derived in detail. Third, an efficient two-stage algorithm with off-line/on-line paradigm is designed to alleviate the intractable computational cost, and its convergence analysis is also derived. Numerical results confirmed that the proposed HOMS framework is high-accuracy and efficient for time-variant nonlinear three-field coupling simulation of large-scale composite structures, and simultaneously economize the computational cost for real-world applications. Furthermore, it is particularly noteworthy that simulating nonlinear multi-physics multi-scale problems of large-scale composite structures requires highly refined computational meshes to accurately resolve the micro-scale features. Such stringent mesh requirements render most conventional engineering solvers computationally prohibitive, failing to yield viable numerical solutions. The proposed method circumvents this bottleneck by delivering highly accurate numerical approximations with significantly reduced computational overhead, thereby offering significant practical value for practical engineering applications. As a future objective, we will focus on the extension of intricate nonlinear multi-physics problems including thermal radiation and convection effects under high-temperature environment, magnetic effect under strong electromagnetic environment. Furthermore, the computational efficiency of the proposed HOMS method is further accelerated in the on-line stage by exploiting the integrated spatiotemporal computational capability of randomized neural network approach.
\bibliographystyle{siamplain}
\bibliography{references1}
\end{document}